\documentclass[a4paper,10pt]{article}
\usepackage[margin = 0.9in,top=0.5in,headsep=9mm]{geometry}

\usepackage{amsmath,amsfonts,amsthm}
\usepackage{comment}
\usepackage{geometry}[1cm]

\usepackage{amsfonts}
\usepackage{amssymb}
\usepackage{multicol}
\usepackage{graphicx}
\usepackage{caption}
\usepackage{setspace}
\usepackage[titletoc,title]{appendix}
\usepackage{float}
\usepackage{lettrine}
\usepackage{changepage}
\usepackage[english]{babel}
\usepackage{csquotes}
\usepackage{textcomp, gensymb}
\usepackage{gensymb}
\usepackage{amsmath,amsthm}
\usepackage{bm}
\usepackage{subfig}
\usepackage{hyperref}
\usepackage{enumerate}
\usepackage{polynom}
\usepackage{enumerate}
\usepackage{mathtools}
\usepackage{algorithm}
\usepackage{algpseudocode}
\usepackage[mathscr]{euscript}
\usepackage{tikz}
\usepackage{graphicx}
\usepackage{hhline}
\usepackage[backend=bibtex]{biblatex}
\usepackage{nicefrac}
\usepackage{epsfig}
\usepackage{ifthen}
\usepackage{etoolbox}

\newcommand{\R}{\mathbb{R}}

\newcommand{\N}{\mathbb{N}}

\newcommand{\I}{\mathcal{I}}

\renewcommand{\L}{\mathcal{L}}

\DeclareMathOperator*{\argmin}{arg\,min}
\newcommand{\abs}[1]{\left| #1 \right|}
\newcommand{\pars}[1]{\left( #1 \right)}

\newcommand{\seq}[1]{\pars{#1_n}_{n\in\N}}

\newcommand{\Ker}[1]{\text{Ker}\pars{#1}}
\newcommand{\set}[1]{\left\{#1 \right\}}

\newcommand{\norm}[1]{\left\Vert #1 \right\Vert }

\renewcommand{\div}{\text{div} \,}

\newcommand{\inner}[2]{\left\langle #1, #2 \right\rangle}

\newcommand{\dom}{\text{dom}\,}
\newcommand{\ran}{\text{Ran}\,}

\definecolor{unibluedark}{RGB}{68, 111, 128}
\definecolor{unibluelight}{RGB}{84, 159, 198}
\def\mathcolor#1#{\@mathcolor{#1}}
\def\@mathcolor#1#2#3{%
  \protect\leavevmode
  \begingroup
    \color#1{#2}#3%
  \endgroup
}
\newcommand{\tv}[1]{\abs{#1}_{\text{TV}}}

\newtheorem{lemma}{Lemma}[section]
\newtheorem{corollary}{Corollary}[section]
\newtheorem{remark}{Remark}[section]
\newtheorem{definition}{Definition}[section]
\newtheorem{proposition}{Proposition}[section]
\newtheorem{theorem}{Theorem}[section]

\newtheorem{example}{Example}[section]

\newcommand{\bu}{\bar u}
\newcommand{\bv}{\bar v}

\renewcommand{\tv}[1]{\norm{D #1}_\mathcal{M}}

\definecolor{myblue}{rgb}{0.18,0.25,0.63}

\newbool{includeImages}
\setbool{includeImages}{true}

\title{Nested Bregman Iterations for Decomposition Problems}

\author{T. Wolf\,\footnote{Institute of Mathematics, University of Klagenfurt, Austria, email: {\tt tobias.wolf@aau.at} }, D. Driggs \footnote{Department of Applied Mathematics and Theoretical Physics (DAMTP)
University of Cambridge, UK, email: {\tt dtd24@cantab.ac.uk} }, 
K. Papafitsoros\,\footnote{School of Mathematical Sciences, Queen Mary University of London, UK, email: {\tt k.papafitsoros@qmul.ac.uk} }\;, 
E. Resmerita\,\footnote{Institute of Mathematics, University of Klagenfurt, Austria, email: {\tt elena.resmerita@aau.at} }\;, and
C.-B. Schönlieb \footnote{Department of Applied Mathematics and Theoretical Physics (DAMTP)
University of Cambridge, UK, email: {\tt cbs31@cam.ac.uk} } \\ }

\begin{document}

\maketitle
\begin{abstract}
We consider the task of image reconstruction while simultaneously decomposing the reconstructed image into components with different features. A commonly used tool for this is a variational approach with an infimal convolution of appropriate functions as a regularizer. Especially for noise corrupted observations, incorporating these functionals into the classical method of Bregman iterations provides a robust method for obtaining an overall good approximation of the true image, by stopping  the iteration early according to a discrepancy principle. However, crucially, the quality of the separate components  depends further on the proper choice of the regularization weights associated to the infimally convoluted functionals. Here, we propose the method of Nested Bregman iterations to improve a decomposition in a structured way. This allows to transform the task of choosing the weights into the problem of stopping the iteration according to a meaningful criterion based on normalized cross-correlation. We discuss the well-definedness and the convergence behavior of the proposed method, and illustrate its strength numerically with various image decomposition tasks employing infimal convolution functionals.
\end{abstract}

\section{Introduction}
 
Image decomposition is an image processing task in which, given a possibly degraded version $f\in Y$ of some ground truth image $x^{\dagger}\in X$, with $X, Y$ being some Hilbert spaces,  one seeks a decomposition of $x^{\dagger}$ into two or more components. Let us assume that $x^{\dagger}$ and $f$ are related through the equation
$Ax^{\dagger}=f$, where $A:X \to Y$ is a bounded, linear operator modeling some degradation process (the \emph{forward} operator), e.g.\  a convolution operator.  One is interested in finding ideal $x_1^\dagger,\dots, x_n^\dagger$ such that $x^\dagger  := \sum\limits_{i = 1}^nx_i^\dagger$ satisfies $Ax^\dagger  = f$, where each of the components $x_1^\dagger,\dots, x_n^\dagger$ is characterized by some specific features. These features can be geometric in nature, with examples including piecewise constant, piecewise smooth or   periodic structures. On the other hand, images also have  less-structured, highly oscillatory components like texture and other fine-scaled features. 
Decomposing an image in such components is a challenging problem, and it is often exacerbated by the addition of a highly oscillatory function $\eta\in Y$ to $f$, modelling some random noise component and resulting in noisy data $f^{\delta}$. Here we assume that the noise is additive and satisfies $\eta:=f^{\delta}-f$, with $\|f^{\delta}-f\|_{Y}\le \delta$. For simplicity, we focus on the case of two components of the image, denoted $u$ and $v$, satisfying $x^\dagger=u^\dagger+v^\dagger$. Thus, one ends up with the following relations:
\begin{equation}\label{eq:ill_posed_problem}
    Ax^{\dagger}+\eta=f^{\delta},
\end{equation}
\begin{equation}\label{eq:ill-posed_problem_components}
A(u^\dagger + v^\dagger) + \eta=f^{\delta}.
\end{equation}
Notably, the classical task of denoising can be regarded as a special instance of image decomposition where $v\equiv 0$. In that case, $x=u$ is obtained as a solution of a variational problem (Tikhonov regularization) 
\begin{equation}\label{intro:Tiknonov}
    \min_{x\in X}\; \frac{\lambda}{2} F(Ax,f^{\delta}) + J(x),
\end{equation}
or alternatively via solving the corresponding Morozov regularization problem
\begin{equation}\label{intro:Morozov}
    \min_{F(Ax,f^{\delta})\le \delta}  J(x).
\end{equation}
Here $F$ is the data-fidelity term which enforces data consistency, while  $J$ is a regularization functional which  makes the problem well-posed and imposes some prior structure on the reconstruction. For instance, $H^{1}$ regularization,  i.e., $J(x)=\frac{1}{2}\|\nabla x\|_{L^{2}}^{2}$ promotes smooth reconstructions \cite{Tikhonov}, while the total variation  seminorm $J(x)=\tv{x}$ \cite{rudin1992nonlinear, ChambolleLions} and its higher order extension, the total generalized variation (TGV) \cite{TGV}, $J(x)=TGV(x)$, promote piecewise constant and piecewise affine reconstructions, respectively.   The regularization parameter $\lambda>0$ in \eqref{intro:Tiknonov}  balances the effects of the data-fidelity and the regularization term.  Using the terminology of image decomposition, $\lambda$ influences how much of the structure imposed by $J$ will go to the component $x=u$ and how much of the noise will be incorporated in $\eta$. Hence, its numerical value has to be chosen carefully to ensure a meaningful reconstruction, which in this case is the decomposition into the estimation of ground truth and noise. 

While  the choice of $F$ in denoising is typically guided by the statistics of the noise, e.g. the $\frac{1}{2}\|\cdot\|_{L^{2}}^{2}$ or $\|\cdot\|_{L^{1}}$ norms for Gaussian or impulse noise, respectively, \cite{nikolova2002minimizers, ChambolleLions}, more sophisticated fidelity terms can be used in the case where \eqref{intro:Tiknonov} is used for decomposing $x^{\dagger}$ into cartoon/structure and texture. A prominent example is Meyer's $G$-norm  \cite{meyer2001oscillating}, for which various approximations have been proposed \cite{vese2003modeling, Osher_Sole_Vese_2003, Vese_negative_sobolev2008}. For our purposes,  we will make use of $\frac{1}{2}\|\cdot\|_{Y}^{2}$ as a fidelity term, often omitting the subscript, while the decomposition will be achieved via a proper choice of $J$. 
With regard to image decomposition where the aim is to recover two components $u$ and $v$, it has been shown that it can be quite advantageous if one employs regularization functionals $J$ which are defined via an infimal convolution of the form 
\begin{equation}\label{intro:inf_conv}
J(x)=(\alpha g \square \beta h)(x):= \min_{\substack{u,v\in X\\ x=u+v}} \{\alpha g(u) +\beta h(v) \},\quad \alpha, \beta>0.
\end{equation}
Here, the structures of $u$ and $v$ are shaped by  $g$ and $h$, respectively. A plethora of such examples can be found in the literature. For instance, a TV-TV$^{2}$ infimal convolution was suggested in \cite{ChambolleLions}, where $g(\cdot)=\alpha|\cdot|_{\mathrm{TV}}$ and $h(\cdot)=\beta|\cdot|_{\mathrm{TV}^{2}}$, the second order total variation, in order to promote reconstructions that consist of piecewise constant and affine parts. This type of infimal convolution  yields approximations of similar structure with TGV, where the two functionals coincide for dimension one \cite{papafitsoros2013study, TGV_learning2}. An analogous functional, employing the $L^{p}$ norm $\|\cdot\|_{L^{p}}$ instead of $|\cdot|_{\mathrm{TV}^{2}}$ was explored in \cite{journal_tvlp, tv_linf}, whose particular case of $p=2$ corresponds to the well-known Huber TV functional \cite{huber1973robust}. It turns out that in the latter case, the resulting functional has also a tight connection with the TV-$H^{1}$ infimal convolution, see Example \ref{ex:H1-TV} in the present work for more details. Infimal convolutions of TGV type functionals have also been considered in the literature, mainly for dynamic imaging \cite{TGV_dynamic_MRI, HollerKunischIC2014}, but also in relation to image decomposition. One such example is the infimal convolution of TGV and the \emph{oscillating} TGV functional, $TGV^{osci}$, \cite{oscICTV}, aiming at decomposing images into piecewise smooth and oscillating components.
A general framework for constructing  infimal convolution regularizers from frequently used functionals was introduced in \cite{Bredies_2022}, where we also refer for further references. 

Using an infimal convolution functional in a variational setting that employs a general forward operator yields the following minimization problem
\begin{equation}\label{intro:inf_conv_variational}
    \min_{u,v\in X}\;\frac{1}{2}\|A(u+v)-f^{\delta}\|^{2}+ \alpha g(u) + \beta h(v),
\end{equation}
where the parameter $\lambda$ in \eqref{intro:Tiknonov} has been absorbed into $\alpha$ and $\beta$. Similarly to our previous remark about the proper selection of $\lambda$, here the values of the regularization weights $\alpha$ and $\beta$ have to be chosen wisely. Large values of $\alpha$, respectively $\beta$, will enforce the component {$v$, respectively $u$}, to be overrepresented in the decomposition. A significant amount of recent research has been focusing on the automatic selection of these parameters. A prominent tool used for this is bilevel optimization. Guided by quality indicators such as the PSNR, SSIM or some norm, many of these approaches aim at optimizing the quality of the overall reconstruction $x=u+v$, but \emph{not} solely the one of $u$ or $v$ \cite{bilevellearning, DelosReyes2021, handbook_2019}. 
However, optimizing the parameters with respect to the quality of the decomposition is a very delicate problem. To the best of our knowledge, there is no consensus on how to quantify the quality of a decomposition. While earlier approaches manually select parameters \cite{AABC_decomposition2005, Vese_divBMO2005}, more recent works use additional information of the expected components to construct method dependent parameter choice rules, often enforcing a balancing principle for the individual regularizers \cite{HuskaDecompositionOnSurfaces_2019, HuskaKangMorigiLanza2021, SignalDecomposition2023}. Another method to assess the quality of a decomposition is to consider the dissimilarity of the individual components, rather than their similarity with the true ones. Based on the work \cite{AujaolGilboaChanOsher2006}, the (cross-)correlation of individual components has been used successfully for automatizing the parameter selection via bilevel methods \cite{TernaryImageDecomposition, QuaternaryImageDecomposition}.

Looking back at problems of type \eqref{intro:Tiknonov}, one approach that avoids the a priori selection of the regularization parameter $\lambda$ is the Bregman iteration \cite{bur-osh, Benning-Burger}. There, the solutions of a series of adaptive variational problems converge to a $J$-minimizing solution of $Ax=f$, in the noiseless case and  under certain assumptions, see Section \ref{sec:Bregman} for details. In practice, the iterations approximate the noisy data. Therefore, employing some discrepancy principle to stop the iteration early typically results in a good reconstruction.

\subsection*{Our contribution}
In this work, we propose a strategy for achieving a desirable decomposition based on appropriately terminating an iterative procedure rather than  making a suitable a priori choice of the regularization weights $\alpha$ and $\beta$ in an infimal convolution regularization functional. We achieve this by introducing the idea of Nested Bregman iterations combined with Morozov regularization. The starting point is the Morozov regularization problem
\begin{equation}\label{intro:starting_point}
    \min_{\substack{u,v\in X \\  \|A(u+v)-f^{\delta}\|\le \delta}} \{\alpha g(u) +\beta h(v)\},
\end{equation}
where we also allow the case $\delta=0$ for the noiseless case. The Nested Bregman iteration scheme consists of an \emph{outer} iteration, which corresponds to running Bregman iterations for \eqref{intro:starting_point} by considering $g$ as a regularizer and $h$ the fidelity term, with the aim of decreasing the share of the component $v$ in the decomposition. This avoids an overrepresentation of $h$, which occurs when $x \simeq v$. At every outer iteration, a minimization problem of the type
\begin{equation}\label{intro:outer}
      (u_l,v_l) \in \argmin\limits_{\substack{u,v\in X \\  \|A(u+v)-f^{\delta}\|\le \delta}} D_{\alpha g}^{p_{l-1}}(u,u_{l-1}) + \beta h(v) 
\end{equation}
needs to be solved, where $D_{\alpha g}^{p_{l-1}}(u,u_{l-1})$ denotes the Bregman distance between $u$ and $u_{l-1}$ at some appropriate $p_{l-1}\in \partial (\alpha g) (u_{l-1})$, see Section \ref{sec:preliminaries} for the definition of Bregman distances. The inner loop of this iteration is an iterative method for solving \eqref{intro:outer}. For the case of a noise-free observation, we propose to use Bregman iterations with data-fitting term $\frac{1}{2}\norm{A(u+v) -f}^2$ and penalty $D_{\alpha g}^{p_{l-1}}(u,u_{l-1}) + \beta h(v)$. This procedure motivates the notion of Nested Bregman iterations. For noise corrupted data, \eqref{intro:outer} constitutes a regularization of the ill-posed equation \eqref{eq:ill-posed_problem_components} and can be solved efficiently by methods from convex optimization. Conceptually, Nested Bregman iterations generate a sequence of regularized solutions to \eqref{eq:ill-posed_problem_components}, in which the contribution of component $v$ to the decomposition is iteratively decreased. We perform a detailed analysis of the scheme, in which we show the existence of the required subgradients that define the corresponding Bregman distances, examine the scheme's convergence properties, and justify the well-definedness of the associated minimization problems for selected infimal convolution regularizers: $L^{1}-H^{1}$ for separating peaks from an otherwise smooth signal, $TV-H^{1}$ for piecewise constant - smooth decompositions, as well as $TGV-TGV^{osci}$ for decomposing an image into piecewise smooth and oscillatory components. We then conduct a series of numerical experiments illustrating that the proposed Nested Bregman iteration scheme equipped with a cross-correlation based stopping criterion can recover   $x=u+v$ {at least as good with respect to PSNR  as related variational methods with optimal parameter choice obtained via either grid search or a  bilevel method}. {Furthermore, in the same time the}   {correct decomposition is attained}. To the best of our knowledge, this is the first time such a rule is used to early terminate an iterative method for an ill-posed problem. Crucially, this can be achieved by setting a suboptimal initialization for the regularizing weights, requiring minimal experimentation and manual tuning. 

\subsection*{Outline of the paper}
In Section \ref{sec:preliminaries}, necessary preliminary notions are mentioned and  notation is fixed. Section \ref{sec:regularization} recalls basic facts about  variational  regularization and classic results on the Bregman iterations. We introduce the Nested Bregman iterations in Section \ref{sec:nested_Bregman} and analyze separately the case of noiseless and noisy data. We discuss the well-posedness of the nested procedure for the aforementioned infimal convolution regularizers in Section \ref{sec:selected_regularizers}, and present the corresponding numerical experiments in Section \ref{sec:numerical_results}.

\section{Preliminaries}\label{sec:preliminaries}
Notation-wise, unless it is ambiguous, we will avoid using subscripts for the inner products and norms of the Hilbert spaces $X$ and $Y$. 
We first recall some notions and properties of convex functionals, subdifferentials, Bregman distances and infimal convolutions.

\begin{definition}
    Let $J:X \to \R\cup\set{+\infty}$ be a convex functional. 
    \begin{itemize}
        \item The domain of $J$ is defined by $\dom J = \set{x \in X : J(x) < +\infty}$. One  says that $J$ is proper if $\dom J \neq \emptyset$.
        \item The subdifferential of $J$ at $x_0 \in \dom J$ is
        \begin{equation*}
        \partial J(x_0) = \set{x^* \in X^*: \inner{x^*}{x-x_0}\le J(x) - J(x_0) \text{ for all } x \in X}.
    \end{equation*}
    \item The Bregman distance between points $x_0, x \in \dom J$ at $x^*\in \partial J(x_0)$ is given by
    \begin{equation*}
        D_J^{x^{\ast}}(x,x_0) = J(x) -J(x_0) -\inner{x^*}{x-x_0}.
    \end{equation*}
    \end{itemize}
\end{definition}

Several subdifferentials calculus rules are mentioned below.

\begin{proposition}\label{prop:properties_subdifferentials}
  \begin{enumerate} \item Let $J:X \to \R\cup \set{+\infty}$ be a proper, convex and lower semicontinuous functional and let $x_0 \in \dom J$. Then 
         $\partial(\lambda J)(x_0) = \lambda \partial J(x_0)$, for any $\lambda>0$.
    \item Let $J_1,J_2:X \to \R\cup\set{+\infty}$ be proper, convex and lower semicontinuous and let $x_0 \in \dom J_1\cap \dom J_2$.  If $(\dom J_1)^\circ\cap \dom J_2 \neq \emptyset$, then \begin{equation*}
        \partial(J_1+J_2)(x_0) = \partial J_1(x_0) + \partial J_2(x_0).
    \end{equation*}
\item Let $A:X \to Y$ be a bounded linear operator, $J: Y \to \R\cup\set{+\infty}$ be proper, convex and lower semicontinuous and $x_0 \in X$ such that $Ax_0 \in \dom J$. 
If {the interior of $\dom J$ and the image of $A$ have non-empty intersection}, then
 \begin{equation*}
    \partial(J\circ A)(x_0) = A^*\partial J(Ax_0),
\end{equation*}
where $A^*$ denotes the adjoint of $A$.
    \end{enumerate}
    \begin{proof}
             {See Proposition 16.6, Theorem 16.47(i) and Proposition 6.9(vii) in \cite{BauschkeCombettes}.} 
    \end{proof}
\end{proposition}

In order to enforce structures in a decomposition, we will use proper, convex and lower semicontinuous functionals $g,h: X \to \R\cup\set{+\infty}$. Finding a good decomposition of a solution to the ill-posed problem means solving the constrained minimization problem
\begin{align}\label{eq:decomposition_problem}
     &\min_{{u,v\in X}}\; \alpha g(u) +\beta h(v)\notag\\
     &s.t. \quad A(u+v)  =f,
 \end{align}
  where $\alpha,\beta >0$ are weights that balance how much the approximate solution is represented in each component. Often, it is useful to consider the two variable problem \eqref{eq:decomposition_problem} as a single variable one. For this, we recall the concept of infimal convolution. 
 \begin{definition}
The infimal convolution of two functions $g,h:X \to\R\cup\set{+\infty}$ is defined as
\begin{equation}\label{eq:def_inf_conv}
J(x):=(g\square h)(x)=\inf_{\substack{{u,v\in X}\\ x=u+v}}\{g(u)+h(v)\}.
\end{equation}
It is called exact at $x\in X$ if the infimum is attained, that is, there exists $\bu\in X$ such that 
\begin{equation}\label{eq:exact}
 J(x)=\inf_{{u\in X}}\;\{g(u)+h(x-u)\}=g(\bu)+h(x-\bu)=g(\bu)+h(\bv),
\end{equation}
with $\bv=x-\bu$.
 Note that $g\square h=h\square g$.
 \end{definition}
This means, if the infimal convolution $J = \alpha g\square \beta h$ is exact, \eqref{eq:decomposition_problem} becomes 
\begin{align}\label{eq:decomposition_problem_J}
    &\min_{{x\in X}}\; J(x) \notag \\ &s.t.\quad Ax = f.
\end{align}

We finally recall the following properties regarding  infimal convolution of proper, convex and lower semicontinuous functions: 

\begin{proposition}\label{prop:inf_conv_properties} Let $g, h:X \to\R\cup\set{+\infty}$ be proper, convex, lower semicontinuous functions, and $J=g\square h$. Then: 
\begin{enumerate}
\item $J$ is convex.
\item  If one of the following holds: 
\begin{enumerate}
\item $g$ is supercoercive, that is $\displaystyle{\lim_{\|u\|\to\infty}\frac{g(u)}{\|u\|}=\infty}$,
\item $h$ is bounded from below and $g$ is coercive, that is $\displaystyle{\lim_{\|u\|\to\infty}g(u)=\infty}$, 
\end{enumerate}
then $J$ is coercive, proper, lower semicontinuous and consequently exact on $X$.
\item[(iii)] $\dom(g\square h) = \dom g + \dom h$.
\item[(iv)] $\partial(g\square h)(x)=\partial g(\bu)\cap\partial h(\bv)$, for any  $x\in dom\,J$ and $\bu,\bv\in X$ such that \eqref{eq:exact} holds. 
\item[(v)]   If $h$ is bounded from below, and both $g$ and $h$ have (weakly) compact sublevel sets, then $J$ has (weakly) compact sublevel sets. In particular, if $g$ and $h$ are coercive, then $J$ is coercive.

\end{enumerate}
\begin{proof}
        For $(i)-(iii)$, see Propositions 12.11, 12.14 and 12.6 in \cite{BauschkeCombettes}, for $(iv)$, see Corollary 2.4.7 in \cite{zalinescu} and for $(v)$, see Theorems 2.3 and 2.9 in \cite{stromberg}.
\end{proof}
\end{proposition}

{For the remainder of this work, we will assume that the functionals $g$ and $h$ are proper, convex, lower-semicontinuous and non-negative. Further, we assume that there is $u_0 \in X$ with $g(u_0) = 0$. This implies $0 \in \partial g(u_0)$, so that by setting $p_0 = 0$, we can write $g(u) = D_g^{p_0}(u,u_0)$ for all $u\in \dom g$. These assumptions will be crucial for showing the well-definedness and convergence properties of the suggested method.}

\section{Regularization methods}\label{sec:regularization}
 Before presenting our proposed method for  decomposition problems, we review classical regularization methods and describe their numerical  performance for decomposition problems. To this end, we fix a $TV-H^{1}$ infimal convolution regularization as detailed next. \\

 \begin{example}\label{ex:H1-TV}
 Let $\Omega \subset \R^2$ be a square and consider  $x^\dagger = u^\dagger+v^\dagger$, where $u^\dagger$ is a smooth component and $v^\dagger$ is the characteristic function of a smaller square $S\subsetneq \Omega $. The operator $A\in \mathcal{L}(L^2(\Omega))$ is given as a convolution operator with a Gaussian  kernel $k$, that is, $x \mapsto Ax=k*x$.  In order to illustrate the behavior of multiple regularization methods, we will use this example for numerical experiments with a kernel of mean $0$ and variance $1$. The noise $\eta$ is realized by a centered normally distributed random variable with mean $0$ and variance $0.05$. We display the corresponding images in Figure \ref{fig:data}. Taking into account the decomposition structure  of $x^\dagger$,  we employ the squared $H^1-$seminorm as a regularizer $g$ for the $u$ component, and the total variation as $h$ for the $v$ component. We will properly define and analyze the functional arising from the infimal convolution of those regularizers in Section \ref{sec:selected_regularizers}. For now, we just formally treat them as \begin{equation}\label{eq:functionals_TV_H^1}
     g(u) = \frac{\alpha}{2}\int\limits_\Omega \abs{\nabla u}^2 , \quad 
     h(v) = \beta \int\limits_\Omega \abs{\nabla v}.
 \end{equation}

\begin{figure}[H]

\centering
\begin{minipage}[c]{0.24\textwidth}
    \includegraphics[scale = 0.4]{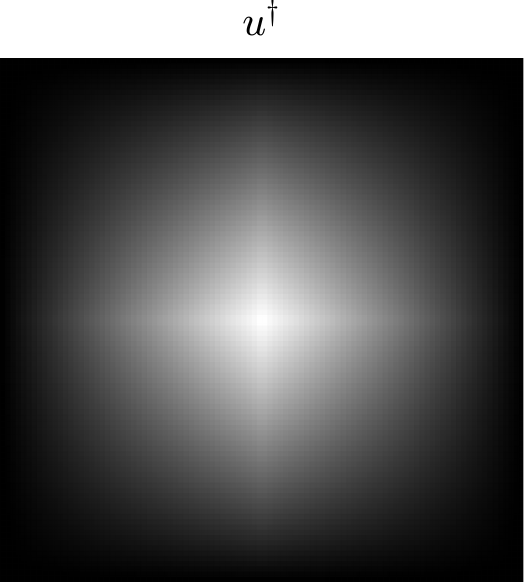}
\end{minipage}
\begin{minipage}[c]{0.24\textwidth}
    \includegraphics[scale = 0.4]{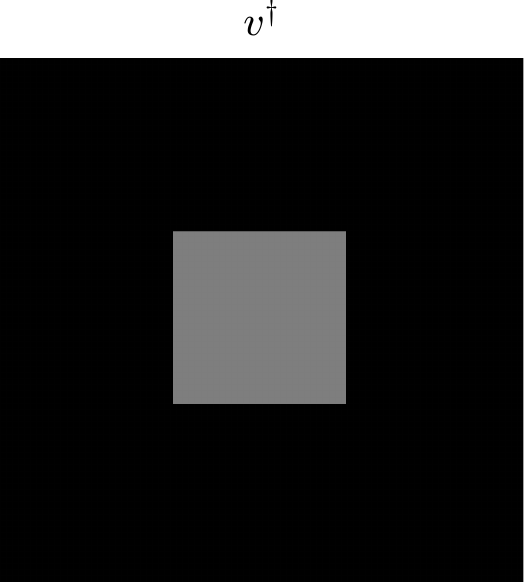}
\end{minipage}
\begin{minipage}[c]{0.24\textwidth}
    \includegraphics[scale = 0.4]{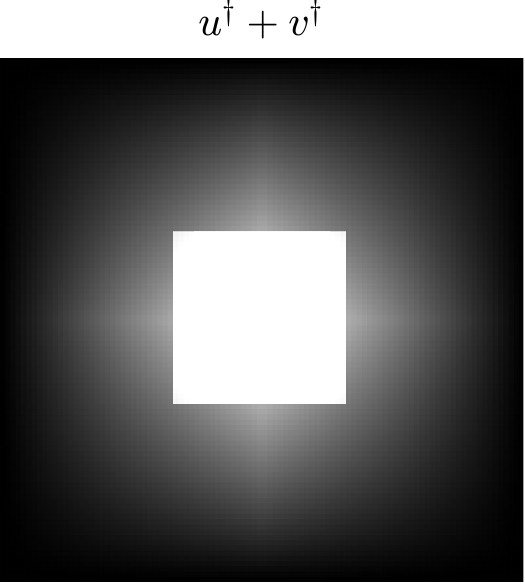}
\end{minipage}
\begin{minipage}[c]{0.24\textwidth}
    \includegraphics[scale = 0.4]{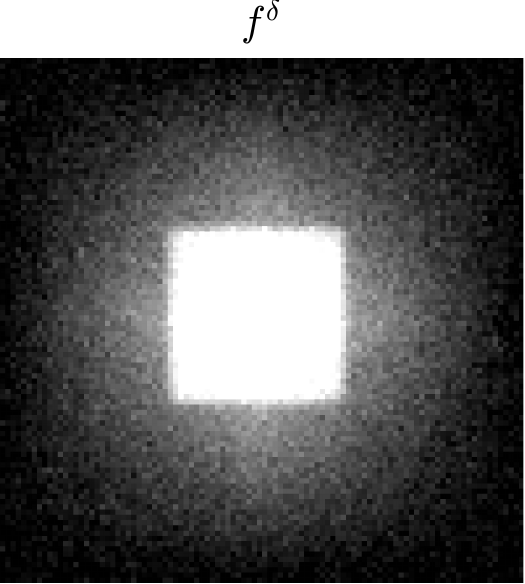}
\end{minipage}

\caption{Components $u^\dagger$ and $v^\dagger$, true image $x^\dagger = u^\dagger+v^\dagger$, as well as the blurred and noisy observation $f^\delta$.}
\label{fig:data}
\end{figure}
 \end{example}

\subsection{Tikhonov Regularization}
The idea of Tikhonov regularization is to approach \eqref{eq:decomposition_problem} by choosing an additional parameter $\lambda >0$ and solving
\begin{equation}\label{eq:variational_problem_J}
    \min\limits_{x \in X} \;\frac{\lambda}{2} \norm{Ax-f}^2 + J(x),
\end{equation}
where $J(x) = \pars{\alpha g \square \beta h}(x)$. If the infimal convolution is exact, this simplifies to the more intuitive problem of solving \begin{equation}\label{eq:variational_problem_decomposition}
    \min\limits_{u,v \in X} \;\frac{\lambda}{2}\norm{A(u+v) -f}^2 + \alpha g(u) +\beta h(v).
\end{equation}
Even though computing solutions of \eqref{eq:variational_problem_decomposition} can often be  relatively cheap, a single-step variational method often comes with some downsides. First,  minimizing $\frac{\lambda}{2}\norm{A(u+v) -f}^2$ instead of solving  $A(u+v) =f$ means that we are no longer guaranteed to obtain a solution of the ill-posed problem. In applications, only perturbed data $f^\delta$ instead of the true image $f$ are available. Thus, if some information about the noise (for instance an estimate of the form $\norm{f-f^\delta}\le \delta$ in the case of additive noise) is available, one chooses the regularization weights accordingly, so that $\norm{Ax-f^\delta}$ is in the range of the noise level $\delta$. Finding optimal parameters for these weights can for example be done using a bilevel approach \cite{HuskaKangMorigiLanza2021, BilevelHuberTV}. However, solving the corresponding bilevel problems can be computationally expensive or unfeasible for highly sophisticated choices of regularizing functionals.

Note that division by $\lambda$ does not affect the minimizers of \eqref{eq:variational_problem_decomposition}, meaning that for practical purposes only $\alpha$ and $\beta$ need to be tuned, while $\lambda$ can be kept constant. Figure \ref{fig:variational_deblurring} shows the decomposition results for  deblurring via Tikhonov regularization for different choices of $\alpha$ and $\beta$, where we set $\lambda = 1$ for convenience. Then, in order to determine appropriate parameters for the infimal convolution, we fix the ratio $\frac{\alpha}{\beta}$, and use a bisection method to find $\alpha$ such that the discrepancy $\delta - \norm{A(u+v)-f}$ is within a $0.1\%$ margin of the noise level $\delta$. Note that this is a successive tuning of the involved parameters and that finding a proper parameter choice for a given ratio is equivalent to determining $\lambda$ in \eqref{eq:variational_problem_decomposition} with $\alpha$ and $\beta$ fixed.  

We observe that the optimal PSNR value of the reconstruction $u+v$ under the performed parameter search strategy is obtained  for $\frac{\alpha}{\beta} = 47$ (Line 3 in Figure \ref{fig:variational_deblurring}). Notably, with this ratio, the reconstructed components $u$ and $v$ are also the closest to the original ones. As expected, if $\alpha$ is chosen too large (relatively to $\beta$), the piecewise constant component $v$ is overrepresented (Line 2 in Figure \ref{fig:variational_deblurring}), whereas a too small choice of $\alpha$ yields an overrepresentation of the smooth component $u$ (Line  5 in Figure \ref{fig:variational_deblurring}). 
\\ 
In conclusion, we observe that the variational approach produces suitable reconstructions with the right choice of parameters. However, finding such parameters is computationally expensive. For our strategy, the tuning of the parameters amounts to  finding the optimal parameters for a given ratio $\frac{\alpha}{\beta}$ and then improving this ratio. Therefore, our aim is to construct a method  which yields a good reconstruction with a meaningful decomposition in a more systematic and, consequently, more efficient manner.  \\

{
\begin{figure}[H]
\centering
\hspace{-0cm}\begin{minipage}[c]{0.24\textwidth}
    \includegraphics[scale = 0.4]{Images/Images_H1_TV/Images_Regularization/u.pdf}
\end{minipage}
\begin{minipage}[c]{0.24\textwidth}
    \includegraphics[scale = 0.4]{Images/Images_H1_TV/Images_Regularization/v.pdf}
\end{minipage}
\begin{minipage}[c]{0.24\textwidth}
    \includegraphics[scale = 0.4]{Images/Images_H1_TV/Images_Regularization/u+v.pdf}
\end{minipage}
\begin{minipage}[c]{0.24\textwidth}
    \includegraphics[scale = 0.4]{Images/Images_H1_TV/Images_Regularization/f.pdf}
\end{minipage}

\vspace{0.5cm}
\centering
\hspace{-0cm}\begin{minipage}[c]{0.24\textwidth}
   \vspace{-0.42cm} \includegraphics[scale = 0.435]{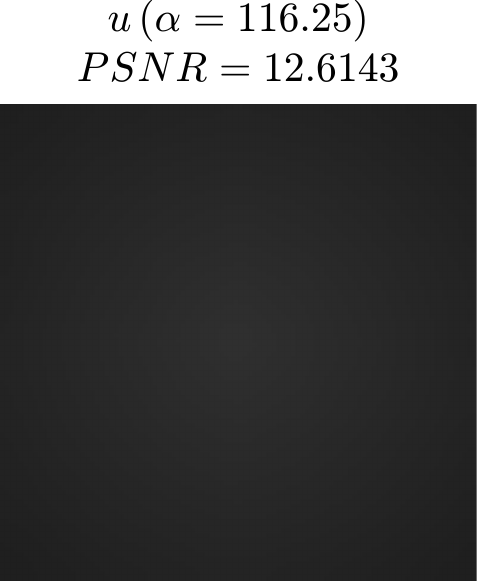}
\end{minipage}
\begin{minipage}[c]{0.24\textwidth}
   \vspace{-0.42cm} \includegraphics[scale = 0.435]{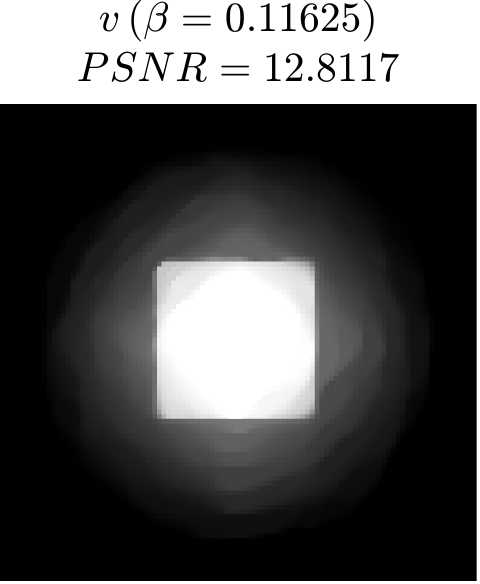}
\end{minipage}
\begin{minipage}[c]{0.24\textwidth}
    \includegraphics[scale = 0.4]{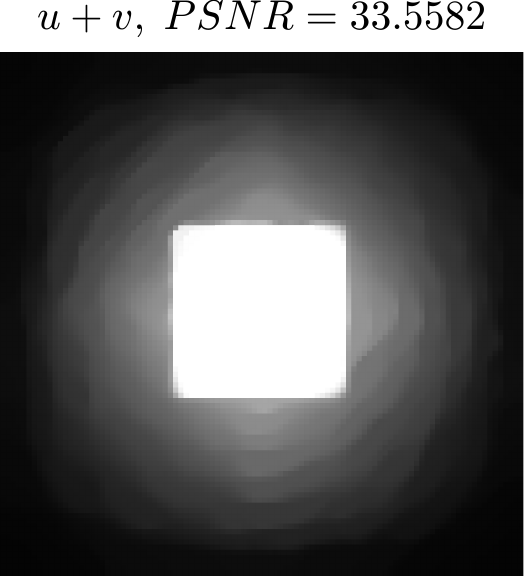}
\end{minipage}
\begin{minipage}[c]{0.24\textwidth}
    \includegraphics[scale = 0.4]{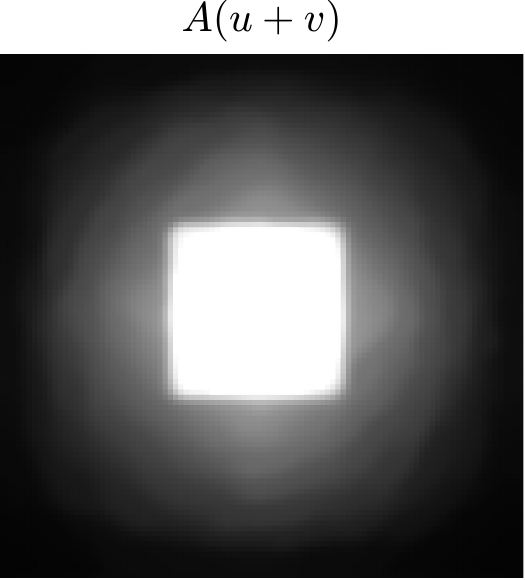}
\end{minipage}

\vspace{0.5cm}
\centering
\hspace{-0cm}\begin{minipage}[c]{0.24\textwidth}
   \vspace{-0.42cm} \includegraphics[scale = 0.435]{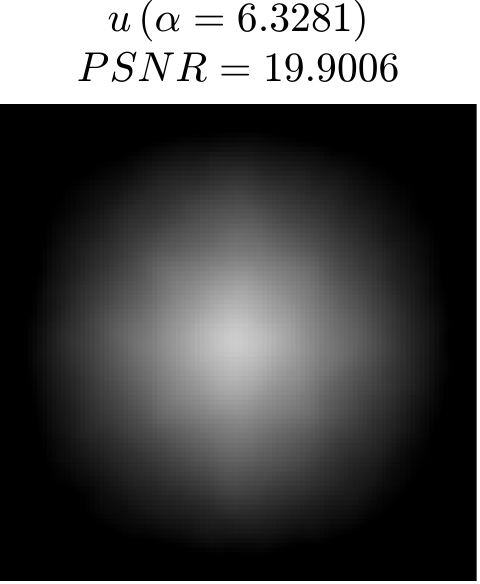}
\end{minipage}
\begin{minipage}[c]{0.24\textwidth}
   \vspace{-0.42cm} \includegraphics[scale = 0.435]{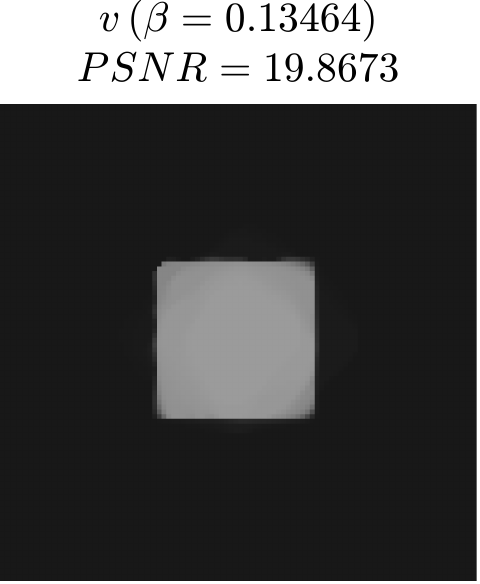}
\end{minipage}
\begin{minipage}[c]{0.24\textwidth}
    \includegraphics[scale = 0.4]{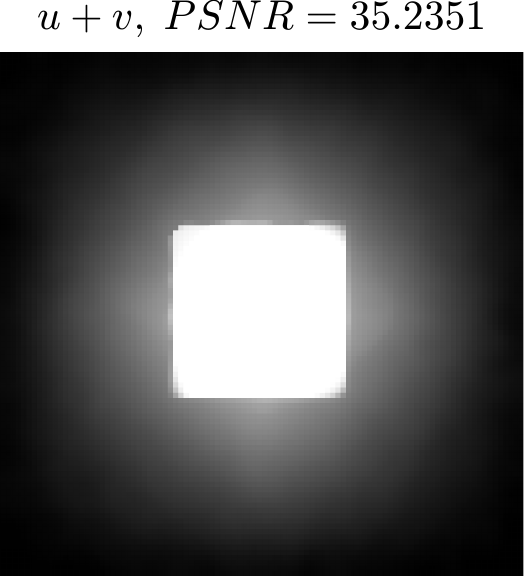}
\end{minipage}
\begin{minipage}[c]{0.24\textwidth}
    \includegraphics[scale = 0.4]{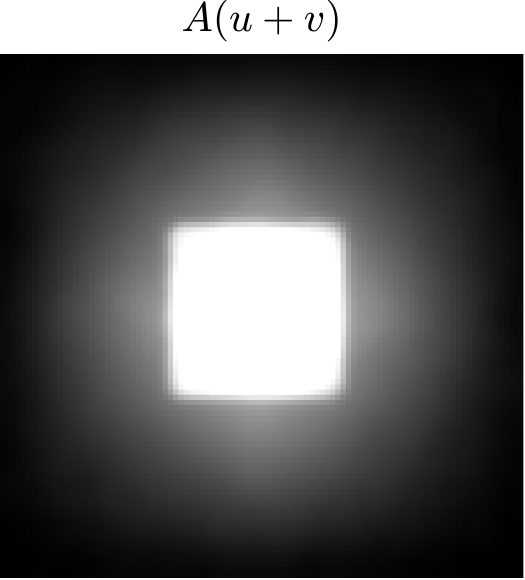}
\end{minipage}

\vspace{0.5cm}
\centering
\hspace{-0cm}\begin{minipage}[c]{0.24\textwidth}
  \vspace{-0.42cm}  \includegraphics[scale = 0.435]{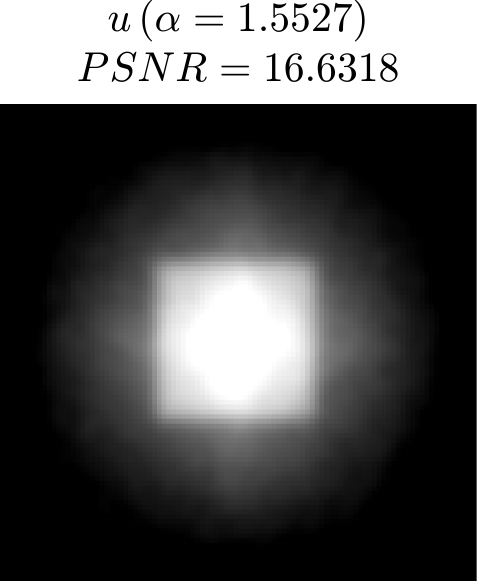}
\end{minipage}
\begin{minipage}[c]{0.24\textwidth}
   \vspace{-0.42cm} \includegraphics[scale = 0.435]{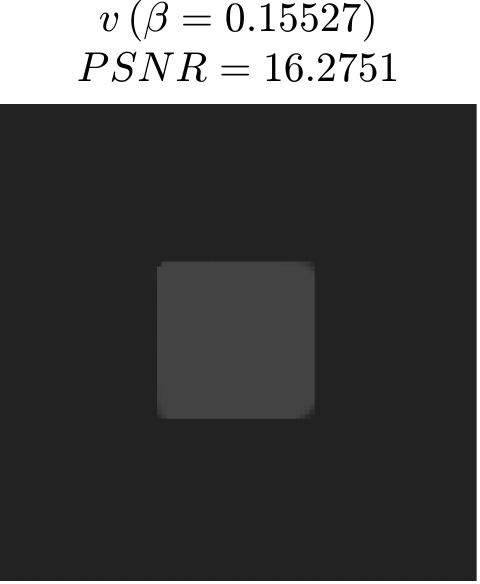}
\end{minipage}
\begin{minipage}[c]{0.24\textwidth}
    \includegraphics[scale = 0.4]{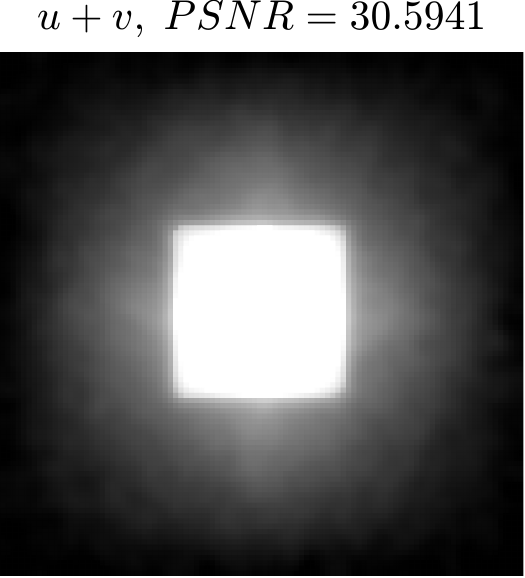}
\end{minipage}
\begin{minipage}[c]{0.24\textwidth}
    \includegraphics[scale = 0.4]{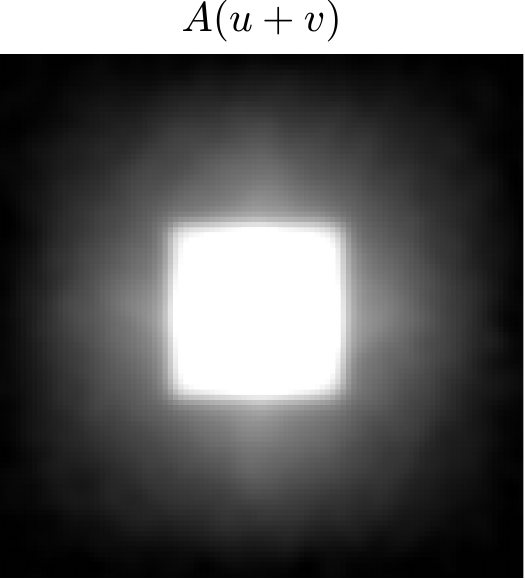}
\end{minipage}

\vspace{0.5cm}
\centering
\hspace{-0cm}\begin{minipage}[c]{0.24\textwidth}
   \vspace{-0.42cm} \includegraphics[scale = 0.435]{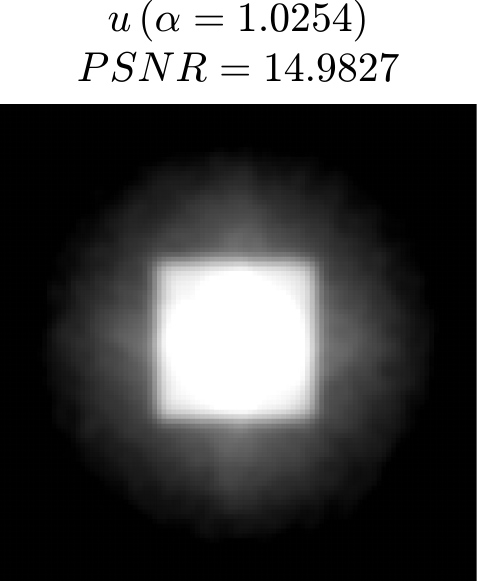}
\end{minipage}
\begin{minipage}[c]{0.24\textwidth}
   \vspace{-0.42cm} \includegraphics[scale = 0.435]{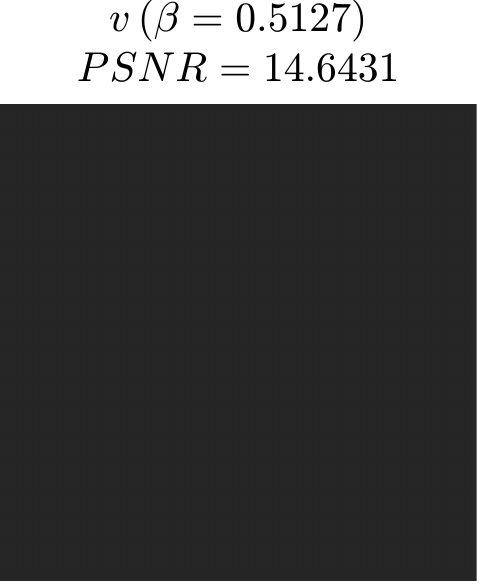}
\end{minipage}
\begin{minipage}[c]{0.24\textwidth}
    \includegraphics[scale = 0.4]{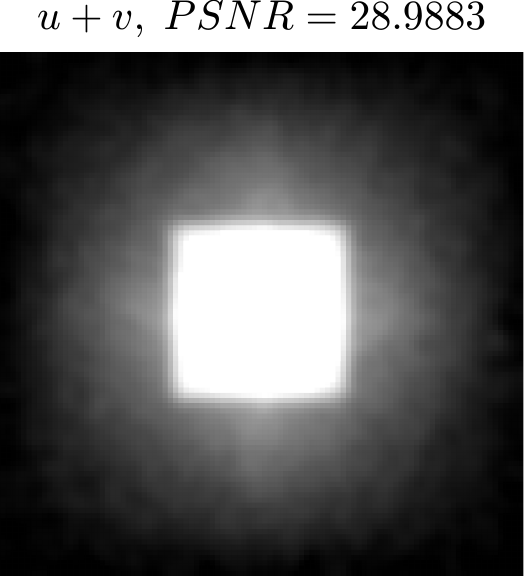}
\end{minipage}
\begin{minipage}[c]{0.24\textwidth}
    \includegraphics[scale = 0.4]{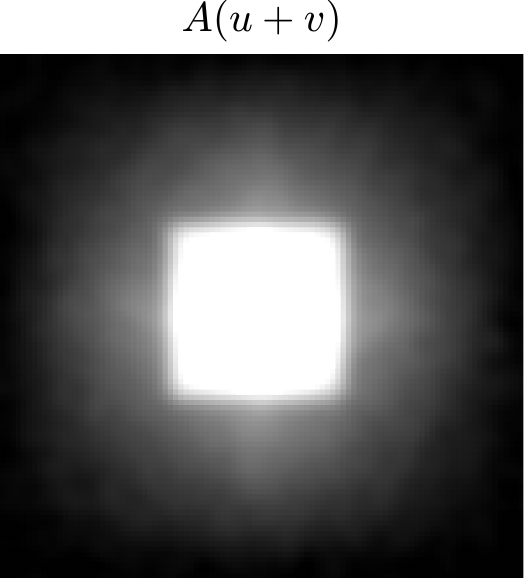}
\end{minipage}

\caption{Different decompositions for the deblurring problem, obtained by Tikhonov regularization with $J(u,v) = \frac{\alpha}{2}\int\limits_\Omega \abs{\nabla u}^2 + \beta \int\limits_\Omega \abs{\nabla v}$ for different ratios $\frac{\alpha}{\beta}$.  First line: true decomposition and observation. Lines $2-5$ (top to bottom): Reconstructions with ration $\frac{\alpha}{\beta} = 1000,47,10,2$}

\label{fig:variational_deblurring}

\end{figure}}{}

\subsection{Morozov Regularization}
Instead of choosing the regularization parameters in \eqref{eq:variational_problem_J} according to the discrepancy principle, one can consider the constraint problem  \begin{equation}\label{eq:Morozov_J}
    \begin{cases}
        &\min\limits_{x \in X} J(x)\\
        & s.t. \ \norm{Ax-f^\delta} \le \delta.
    \end{cases}
\end{equation}
This method is known as Morozov regularization or residual method, see e.g. \cite{IvanovVasinTanana,Lorenz_2013,GrasmairHaltmeierScherzer,Kaltenbacher_Morozov}. Assuming that $J = (\alpha g)\square (\beta h)$ and that the infimal convolution is exact, we can again simplify \eqref{eq:Morozov_J},  that is, we are interested in the solution of \begin{equation}\label{eq:Morozov_decomposition}
    \begin{cases}
        &\min\limits_{u,v \in X} \alpha g(u) + \beta h(v)\\
        & s.t. \ \norm{A(u+v)-f^\delta} \le \delta.
    \end{cases}
\end{equation}
In fact, under mild conditions, \eqref{eq:Morozov_J} is equivalent to \eqref{eq:variational_problem_J}. Namely, if there  is $x \in \dom J$ such that $\norm{Ax-f^\delta} < \delta$, then by \cite[Theorem 3.9]{barbu2012convexity}, there exists $\lambda$ 
such that the solution $x^*$ of \eqref{eq:Morozov_J} minimizes \eqref{eq:variational_problem_J}. Conversely, a solution $x_\lambda$ of \eqref{eq:variational_problem_J} solves \begin{equation*}
    \begin{cases}
        &\min\limits_{x \in X} J(x)\\
        &s.t. \ \norm{Ax-f^\delta} \le \norm{Ax_{\lambda}-f^\delta}.
    \end{cases}
\end{equation*}
The Morozov formulation of the variational regularization overcomes the issue of choosing the parameter $\lambda$. Thus, dividing by $\beta$ in \eqref{eq:Morozov_decomposition}, we only need to choose an appropriate value for the ratio $\frac{\alpha}{\beta}$ so that the decomposition problem can essentially be reduced to tuning a single parameter. For a more detailed analysis of the relation between Tikhonov and Morozov regularization in the case that $J(x) = \norm{Lx}^p$ for some linear map $L$ and $p>1$, the reader is referred to \cite[Section 3.5]{IvanovVasinTanana}.
\subsection{Bregman Iterations}\label{sec:Bregman}
 Another way of overcoming the need to choose the regularization parameter $\lambda$ in \eqref{eq:variational_problem_J} is performing Bregman iterations (see, for instance, \cite{orig_breg, Benning-Burger}). 

The initial step of Bregman iterations consists of computing \begin{equation}\label{eq:Bregman_initial_iterate}
    x_1 \in \argmin\limits_{x\in X}\set{\frac{{\lambda}}{2}\norm{Ax-f}^2 + J(x)}.
\end{equation}
Next,  {set $\xi_1= -\lambda A^*(Ax_1-f) \in\partial  J(x_1)$} and iterate for $k \ge 2$:
\begin{enumerate}
    \item Compute \begin{equation}\label{eq:Bregman_iteration_step}
        x_k\in \argmin\limits_{x \in X}\set{\frac{\lambda}{2}\norm{Ax-f}^2 + D_J^{\xi_{k-1}}(x,x_{k-1})}
    \end{equation}
    \item  {Set $\xi_k = \xi_{k-1}-\lambda  A^*(Ax_k-f)\in \partial J(x_k)$}.
\end{enumerate}

\begin{remark}\label{rem:Bregman_well_definedness}
The variational problem \eqref{eq:Bregman_iteration_step} is equivalent to \begin{equation*}
    x_k\in \argmin\limits_{x\in X}\set{\frac{\lambda}{2}\norm{Ax - (f+{\zeta_{k-1})}}^2 + J(x)},
\end{equation*}
where $\zeta_{k-1} = \sum\limits_{i = 1}^{k-1}(f-Ax_{i})$. Therefore,  Bregman iterations are well-defined, if \eqref{eq:variational_problem_J} is well-defined for all possible data.
\end{remark}

{We recall}  the following convergence results {for the Bregman iterates $x_k$, where $k \in \N= \set{1,2,\dots}$}.

\begin{theorem}\label{thm:Convergence_Bregman_Noisefree}
   \cite{bur-osh,fri_sch,fri_lor_res,Benning-Burger} Let $(x_k)_{k\in \N}$ be a sequence generated by \eqref{eq:Bregman_initial_iterate} and \eqref{eq:Bregman_iteration_step}. Assume that there is $x^\dagger \in \dom J$ that minimizes $\norm{Ax-f}^2$. Then one has \begin{equation}\label{eq:ConvergenceResidual_Bregman_noisefree}
        \norm{Ax_k-f}^2-\norm{Ax^\dagger-f}^2 \le \frac{J(x^\dagger)}{{\lambda}k},
    \end{equation}
    and
    \begin{equation}\label{eq:BoundIterate_Bregman_Noisefree}
        J(x_k) \le \frac{5}{2} J(x^\dagger),
    \end{equation}
    for all $k \in \N$. If, additionally, $Ax^\dagger = f$ and there is $q\in Y $ such that $A^*q\in\partial J(x^\dagger)$, then each weak cluster point  of $(x_k)_{k\in \N}$ is a $J$-minimizing solution of $Ax =f$.
\end{theorem}

In the case of noisy data $f^\delta$ with $\norm{f-f^\delta}\le \delta$, one can use the discrepancy principle to determine a stopping index of the iteration. That is, one chooses $\tau >1$ and stops the iteration at index $k^\delta$ defined by \begin{equation}\label{eq:discrepancy_principle}
    k^\delta = \min\set{k\in \N : \norm{Ax_k-f^\delta}<\tau\delta}.  
\end{equation}
With this stopping rule, one obtains the following version of Theorem \ref{thm:Convergence_Bregman_Noisefree}.
\begin{theorem}\label{thm:Convergence_Bregman_Noisy} \begin{enumerate}
    \item     Let $\left(x_k^\delta\right)$ be the sequence generated by \eqref{eq:Bregman_initial_iterate} and \eqref{eq:Bregman_iteration_step} with $f$ replaced by $f^\delta$. Assume that $\norm{f-f^\delta} \le \delta$ and that there is $x^\dagger \in \dom J$ which minimizes $\norm{Ax-f}^2$. Then one has\begin{equation}\label{eq:ConvergenceResidual_Bregman_noisy}
        \norm{Ax_k^\delta-f}^2-\norm{Ax^\dagger-f}^2 \le \delta^2 +\frac{J(x^\dagger)}{k}
    \end{equation}
        and
    \begin{equation}\label{eq:BoundIterate_Bregman_noisy}
        J(x_k^\delta) \le k\delta +\frac{5}{2} J(x^\dagger),
    \end{equation}
    for all $k \in \N$. 
    \item Additionally, assume that $Ax^\dagger = f$ and there is $q\in Y $ such that $A^*q\in\partial J(x^\dagger)$. Let   $(\delta_n)_{n\in \N}$ be a sequence with $\lim\limits_{n\to\infty}\delta_n = 0$, and let observations $f^{\delta_n}$  satisfy $\norm{f-f^{\delta_n}}\le \delta_n$. If $(x_{k^{\delta_{n}}})_{n \in \N}$ is a sequence obtained by stopping the iteration according to \eqref{eq:discrepancy_principle} for each $\delta_n$, then every weak cluster point of $\left(x_{k^{\delta_{n}}}\right)$ is a $J$-minimizing solution of $Ax = f$. 
\end{enumerate}
\end{theorem}
Theorem \ref{thm:Convergence_Bregman_Noisy} illustrates one major advantage of Bregman iterations. Instead of choosing a parameter $\lambda$, the iteration is stopped according to \eqref{eq:discrepancy_principle}. Therefore, only an estimate of the noise is necessary to obtain a good approximation of the ground truth. Additionally, in the case of exact data, the limit points of the iterates give a solution of \eqref{eq:decomposition_problem} if the data satisfies  the source condition $A^*q \in \partial J(x^\dagger)$. On the other hand, if the source condition is not satisfied, we might be looking for something different from a $J$-minimizing solution (cf \cite{Benning-Burger}).\\
As compared to Morozov regularization, the regularized solutions obtained from Bregman iterations usually yield better reconstructions, possibly due to the fact that they have less bias towards $J$. Namely, while the solutions $x^*$ of \eqref{eq:Morozov_J} satisfy $J(x_\lambda) \le J(x^\dagger)$, the  Bregman iterations verify $J(x_k) \le 2J(x^\dagger)$ (compare to \eqref{eq:BoundIterate_Bregman_noisy} in the noisy data case). Thus, the Morozov regularizers tend to be over-regularized. For a more detailed discussion, see Section 6.1 in \cite{Benning-Burger}. 
\\ Coming back to the decomposition problem, we note that for exact infimal convolutions, Bregman distances can be decomposed as shown below  in \eqref{eq:ecomposition_Bregman_distances}. 
Indeed, let $J = (\alpha g)\square (\beta h)$ and assume that the infimal convolution is exact. Let further $x,\hat x\in dom\,J$, $\hat \xi\in\partial J(\hat x)$. Due to the exactness of the infimal convolution, there exist $u,v,\hat u,\hat v\in X$ such that
\begin{equation*}
x=u+v,\quad \hat x=\hat u+\hat v,\quad J(x)=\alpha g(u)+\beta h(v),\quad J(\hat x)=\alpha g(\hat u)+\beta h(\hat v).
\end{equation*}
Then   Proposition \ref{prop:inf_conv_properties} (iv) yields 
\begin{eqnarray}\label{eq:ecomposition_Bregman_distances}
\notag D_J^{\hat\xi}(x,\hat x)&=& J(x)-J(\hat x)-\langle\hat\xi,x-\hat x\rangle\\\notag
&=& \alpha g(u)+\beta h(v)-\alpha g(\hat u)-\beta h(\hat v)-\langle\hat \xi,u+v-\hat u-\hat v\rangle\\
&=&  D_g^{\hat\xi}(u,\hat u)+D_h^{\hat\xi}(v,\hat v).
\end{eqnarray}
Hence, in this context, \eqref{eq:Bregman_iteration_step} becomes \begin{equation}\label{eq:Bregman_iteration_step_decomposition}
            (u^k,v^k)\in \argmin\limits_{u,v \in X}\set{\frac{1}{2}\norm{A(u+v)-f}^2 + D_{\alpha g}^{\xi_{k-1}}(u,u_{k-1}) + D_{\beta h}^{\xi_{k-1}}(v,v_{k-1})},
\end{equation}
with $\xi_k = \xi_{k-1}+ A^*(f-A(u_k+v_k))\in \partial(\alpha g)(u_k) \cap \partial (\beta h)(v_k) $.

Note that this formulation is precisely the one proposed in \cite{li_et_al}. However, no connection to infimal convolutions is mentioned or analyzed there.
\\
Figure \ref{fig:bregman_deblurring} shows the results of running Bregman iterations for Example \ref{ex:H1-TV} until the discrepancy principle \eqref{eq:discrepancy_principle} with $\tau = 1.001$ is met. As a consequence, the residual at the stopping index of Bregman iterations is not larger than the residual obtained earlier with Tikhonov regularization. In practice, we also observe that the residuals of the Bregman iteration and the ones obtained in the Tikhonov approach with a bisection for the parameter choice do not differ too much. The parameters $\alpha$ and $\beta$ were chosen to be $4$ times the parameters in the corresponding line of Figure \ref{fig:variational_deblurring}. That way, the ratios $\frac{\alpha}{\beta}$ are the same for both methods. The multiplication with the factor $4$ ensures that the iteration does not terminate after the initial step \eqref{eq:Bregman_initial_iterate}, so that the difference between Bregman iterations and a single step approach becomes clear. The experiments also suggest that, as long as the ratio $\frac{\alpha}{\beta}$ is kept constant, the reconstructions do not change much depending on the specific choices of $\alpha$ and $\beta$. We observe the largest PSNR again for the case $\frac{\alpha}{\beta} = 47$ (Line 3). As expected, the bias reduction towards the regularizer results in  the PSNR obtained by using Bregman iterations is consistently larger than the one obtained with variational regularization.  Furthermore, in all cases, the algorithm was terminated after only $3$ iterations. Of course, this is also due to  the parameters already being in a reasonable range. However, Bregman iterations in our experiments required less minimization steps than the variational approach with a bisection method to find suitable parameters. This poses another advantage, as no expensive search for suitable parameters is necessary. However, we observe again that only for a specific ratio $\frac{\alpha}{\beta}$ the components are accurately reconstructed. Thus, we conclude that Bregman iterations are not suitable for obtaining a good reconstruction of the individual components.

{
\begin{figure}[H]
\centering
\hspace{-0cm}\begin{minipage}[c]{0.24\textwidth}
    \includegraphics[scale = 0.4]{Images/Images_H1_TV/Images_Regularization/u.pdf}
\end{minipage}
\begin{minipage}[c]{0.24\textwidth}
    \includegraphics[scale = 0.4]{Images/Images_H1_TV/Images_Regularization/v.pdf}
\end{minipage}
\begin{minipage}[c]{0.24\textwidth}
    \includegraphics[scale = 0.4]{Images/Images_H1_TV/Images_Regularization/u+v.pdf}
\end{minipage}
\begin{minipage}[c]{0.24\textwidth}
    \includegraphics[scale = 0.4]{Images/Images_H1_TV/Images_Regularization/f.pdf}
\end{minipage}

\vspace{0.5cm}
\centering
\hspace{-0cm}\begin{minipage}[c]{0.24\textwidth}
    \vspace{-0.42cm}\includegraphics[scale = 0.435]{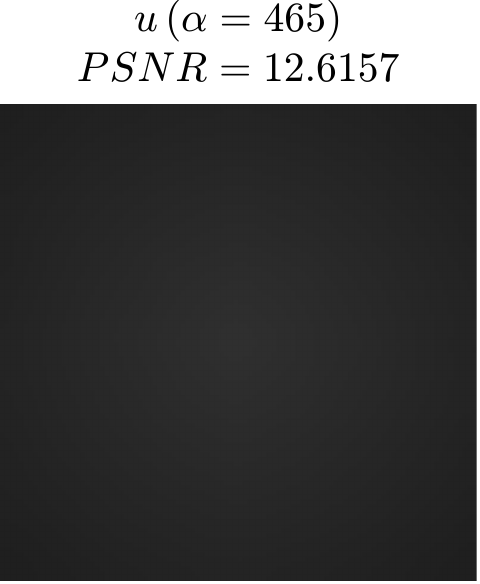}
\end{minipage}
\begin{minipage}[c]{0.24\textwidth}
  \vspace{-0.42cm}  \includegraphics[scale = 0.435]{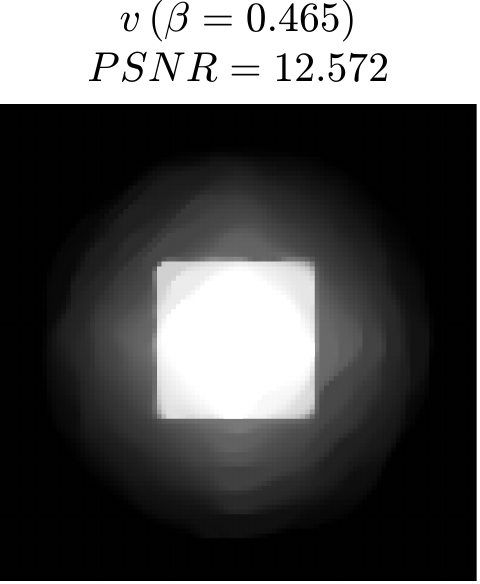}
\end{minipage}
\begin{minipage}[c]{0.24\textwidth}
    \includegraphics[scale = 0.4]{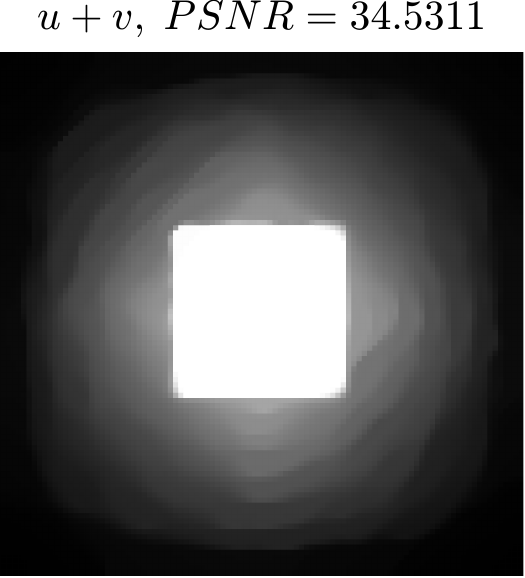}
\end{minipage}
\begin{minipage}[c]{0.24\textwidth}
    \includegraphics[scale = 0.4]{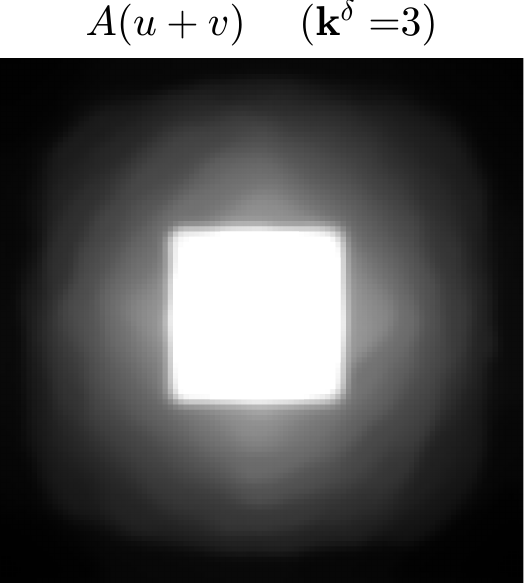}
\end{minipage}

\vspace{0.5cm}
\centering
\hspace{-0cm}\begin{minipage}[c]{0.24\textwidth}
  \vspace{-0.42cm}  \includegraphics[scale = 0.435]{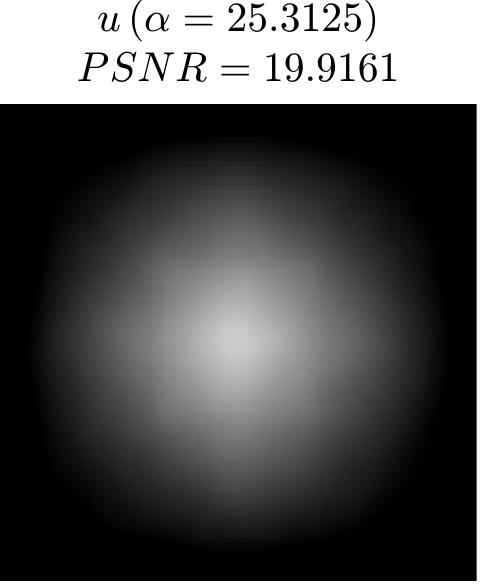}
\end{minipage}
\begin{minipage}[c]{0.24\textwidth}
  \vspace{-0.42cm}  \includegraphics[scale = 0.435]{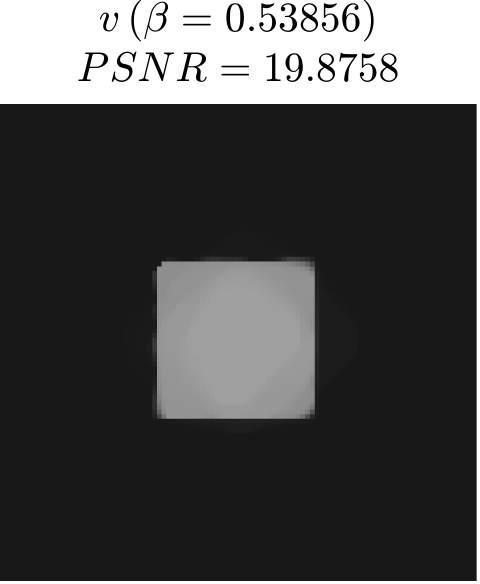}
\end{minipage}
\begin{minipage}[c]{0.24\textwidth}
    \includegraphics[scale = 0.4]{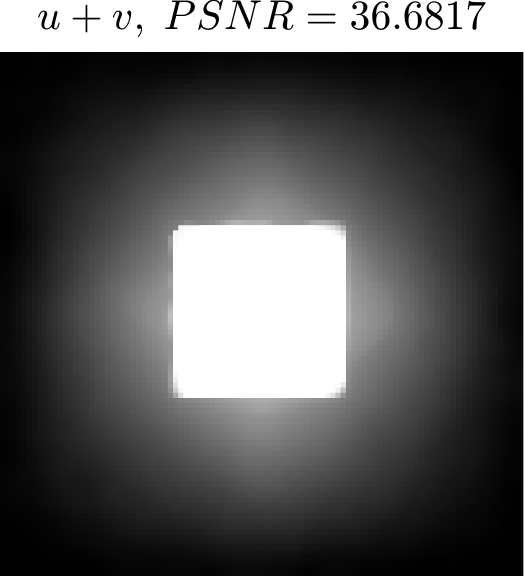}
\end{minipage}
\begin{minipage}[c]{0.24\textwidth}
    \includegraphics[scale = 0.4]{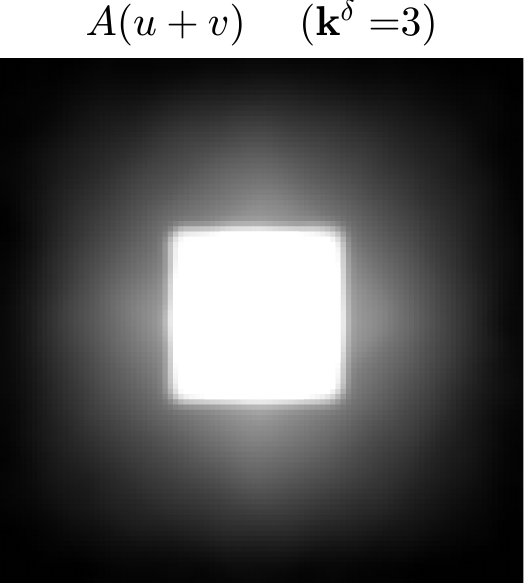}
\end{minipage}

\vspace{0.5cm}
\centering
\hspace{-0cm}\begin{minipage}[c]{0.24\textwidth}
   \vspace{-0.42cm} \includegraphics[scale = 0.435]{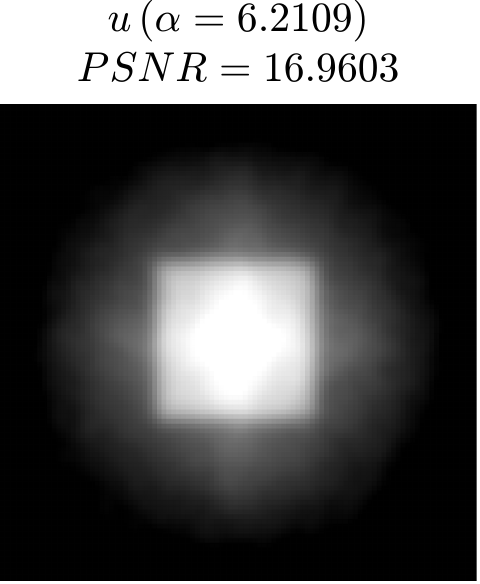}
\end{minipage}
\begin{minipage}[c]{0.24\textwidth}
   \vspace{-0.42cm} \includegraphics[scale = 0.435]{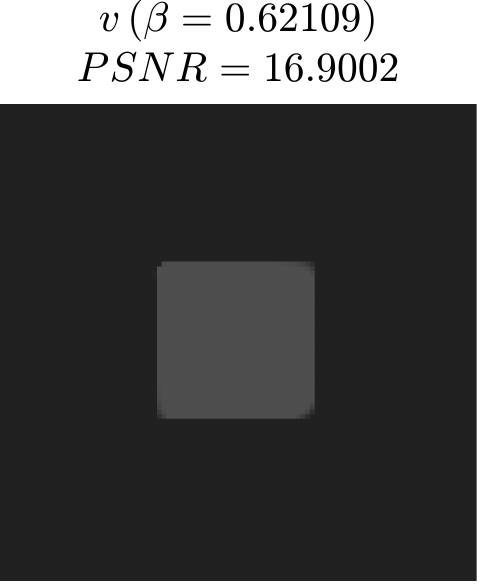}
\end{minipage}
\begin{minipage}[c]{0.24\textwidth}
    \includegraphics[scale = 0.4]{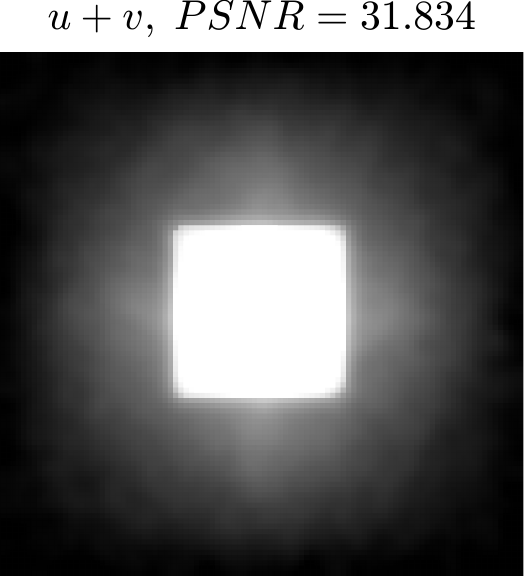}
\end{minipage}
\begin{minipage}[c]{0.24\textwidth}
    \includegraphics[scale = 0.4]{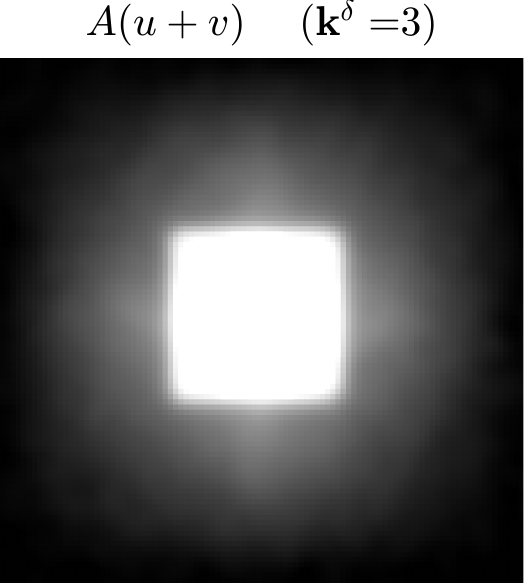}
\end{minipage}

\vspace{0.5cm}
\centering
\hspace{-0cm}\begin{minipage}[c]{0.24\textwidth}
   \vspace{-0.42cm} \includegraphics[scale = 0.435]{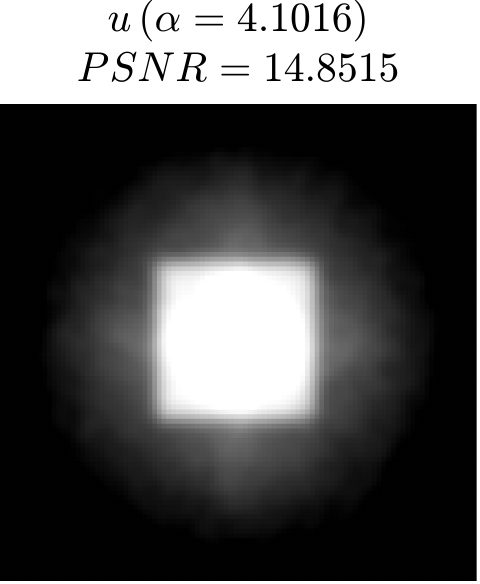}
\end{minipage}
\begin{minipage}[c]{0.24\textwidth}
  \vspace{-0.42cm}  \includegraphics[scale = 0.435]{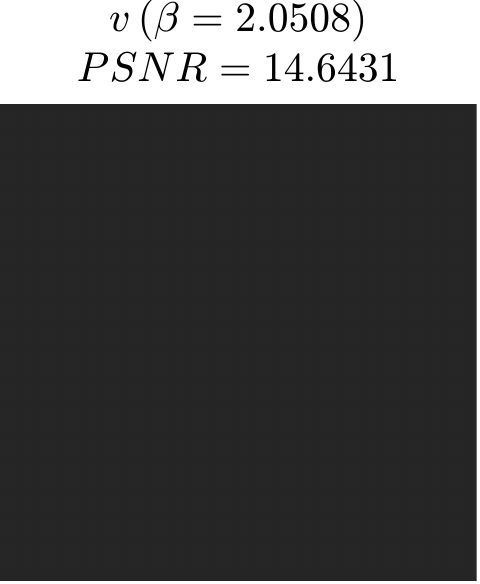}
\end{minipage}
\begin{minipage}[c]{0.24\textwidth}
    \includegraphics[scale = 0.4]{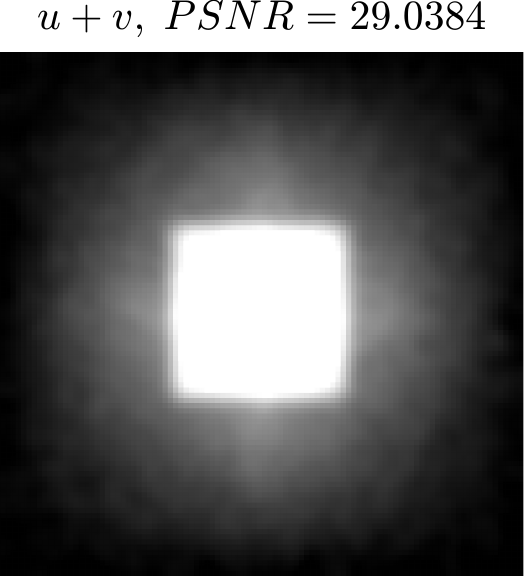}
\end{minipage}
\begin{minipage}[c]{0.24\textwidth}
    \includegraphics[scale = 0.4]{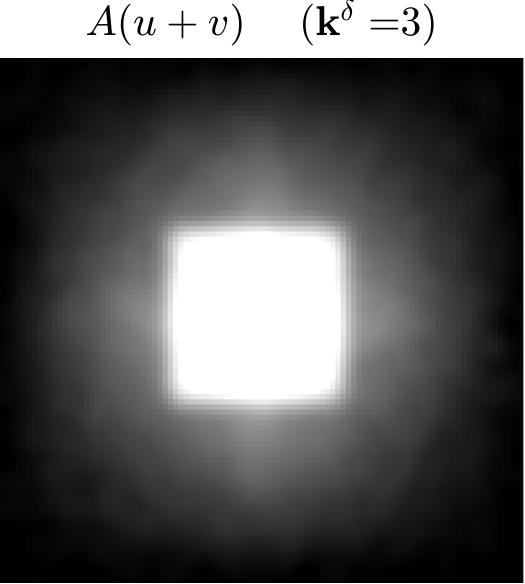}
\end{minipage}

\caption{Different decompositions of the deblurring problem, obtained by Bregman iterations  with $J(u,v) = \frac{\alpha}{2}\norm{\nabla u}_{L^2}^2 + \beta \tv{v}$ stopped via discrepancy principle $(\tau =1.001)$ First line: true decomposition and observation. Lines $2-5$ (top to bottom): Reconstructions with ration $\frac{\alpha}{\beta} = 1000,47,10,2$}

\label{fig:bregman_deblurring}

\end{figure}}{}

\section{Nested Bregman iterations}\label{sec:nested_Bregman}
As motivated in the seminal paper \cite{orig_breg}, by using Bregman iterations, one can avoid choosing the regularization parameter $\lambda$ in the Tikhonov approach. Recall that upon division by $\alpha$ or $\beta$, problem \eqref{eq:decomposition_problem} can be seen as a variational problem with only one parameter to choose. The idea of Nested Bregman iterations is to use Bregman iterations to transform the problem of choosing this parameter into the problem of stopping an iterative procedure. This means the following: Assume that the component $v$ is overrepresented in the decomposition, i.e.\ $x \simeq v$. We therefore aim to decrease its share in the decomposition. However, instead of increasing the value of $\beta$ and solving \eqref{eq:decomposition_problem} again with the adjusted weight, we  replace $\alpha g$ by a suitable Bregman distance and minimize again. That is, we run Bregman iterations with data fidelity term $h$ and regularizer $g$. Hence, we iteratively solve constraint problems of the form \eqref{eq:decomposition_problem}, where instead of $\alpha g$ we use a Bregman distance associated with this function at each step. Since the minimizers of these constraint problems are precisely the $J$-minimizing solutions of $A(u+v) = f$, they can be computed as the weak limit points of Bregman iterations (see Theorem \ref{thm:Convergence_Bregman_Noisefree}). Hence, we obtain a Nesting of Bregman iterations. However, when dealing with noise corrupted observations, multiple ways to define Nested Bregman iterations are thinkable. In the following, we will present two options in this respect, and outline their advantages. 
Since the algorithms we propose do not involve any choice of parameters, we simplify notation and write $g$ and $h$ instead of $\alpha g$ and $\beta h$. That is, we assume that an initial guess for the weightings is implicitly contained in the functionals $g$ and $h$. 

\subsection{The case of noise-free data}
We propose the following algorithm for the case of exact data $f$ (i.e. $\delta = 0$).

\begin{algorithm}[H]
  \begin{algorithmic}[1]
  \State Compute \begin{equation}\label{eq:NestedBregman_noiseless_initial}
      (u_{1},v_{1}) \in \argmin\limits_{A(u+v) =f} \set{g(u) +h(v)}
  \end{equation} using Bregman iterations.
  \State Choose $p_{1} \in \partial g(u_1)\cap \partial h(v_1)\cap \Ker A ^\perp$
 \For{$l =2,3,\dots$}
    \State Compute \begin{equation}\label{eq:NestedBregman_noiseless_step}
        (u_l,v_l) \in \argmin\limits_{A(u+v) = f} D_g^{p_{l-1}}(u,u_{l-1}) + h(v) 
    \end{equation} using Bregman iterations.
    \State Choose $p_{l}\in \partial g(u_l)\cap \Ker A ^\perp$ such that $p_l-p_{l-1} \in \partial h(v_l)$.
  \EndFor
  \end{algorithmic}
  \caption{Nested Bregman Iterations}\label{algo:NestedBregman_Noiseless}
\end{algorithm}

We want to stress that for the latter analysis of the proposed method, the choice of subdifferentials is crucial. The constraint $p_l \in \partial g(u_l)$ with $p_{l}- p_{l-1} \in \partial h(v_l)$ for $l\ge 2$ is necessary to establish convergence results for the sequence $(h(v_l))_{l\in\mathbb{N}}$, while the assumption $p_l \in \Ker A^\perp$ plays an important role in establishing the existence of minimizers for problem \ref{eq:NestedBregman_noiseless_step}.
It is straightforward to modify Algorithm \ref{algo:NestedBregman_Noiseless} if  $f$ is not attainable. In that case, we want to solve \eqref{eq:NestedBregman_noiseless_initial} and \eqref{eq:NestedBregman_noiseless_step} under the constraint that the distance of $A(u+v)$ and $f$ is minimal with respect to some similarity measure. For instance, we can replace the constraint $A(u+v)=  f$ by the normal equation $A^*A(u+v) =A^*f$ and solve the minimization problems \eqref{eq:NestedBregman_noiseless_initial} and \eqref{eq:NestedBregman_noiseless_step} among minimum-norm solutions instead. This is also consistent with solving problems \eqref{eq:NestedBregman_noiseless_initial} and \eqref{eq:NestedBregman_noiseless_step} using Bregman iterations, as defined in Section \ref{sec:Bregman}. However, Bregman iterations can also be used with other data fidelities than $\norm{\cdot}^2$. For more details, see \cite{Benning-Burger}.

\begin{remark}
    If $A = Id$, the constraint in \eqref{eq:NestedBregman_noiseless_initial} and \eqref{eq:NestedBregman_noiseless_step} becomes $u+v =f$. This means we are actually computing the minimizers in the infimal convolution of $D_g$ and $h$ if it is exact. In particular, for the choice $h = \frac{1}{2}\norm{\cdot}^2$, Algorithm \ref{algo:NestedBregman_Noiseless} coincides with classical Bregman iterations for the denoising problem. 
\end{remark}

For the remainder of this section, we will assume that all the minimization problems considered actually admit minimizers. The verification of this will be done for the concrete examples that we consider. However, we still need to verify the well-definedness of the choice of $p_l$. While the existence of $p_l \in \partial g(u_l)$ such that $p_l-p_{l-1} \in \partial h(v_l)$ follows under mild assumptions from the optimality condition of \eqref{eq:NestedBregman_noiseless_step}, the inclusion $p_l \in  \Ker A^\perp$ needs additional assumptions on the functions $g$ and $h$. As a first step, we compute the subgradient of the indicator function for the constraint $Ax = f$.  

\begin{proposition}\label{prop:subdifferential_affine_set}
    Let $C =\set{x \in X: Ax = f}$ and denote by \begin{equation}\label{eq:Indicator_function}
        \I_C(x) = \begin{cases}
            0 & \text{ if } x \in C\\
            \infty &\text{ if } x \notin C
        \end{cases}
    \end{equation} the indicator function of $C$. Then \begin{equation}
        \partial \I_C(x) = \Ker{A}^\perp 
    \end{equation} for all $x \in C$.
    \begin{proof}
       Let $\bar x \in C$. Since $C$ is convex, it holds \begin{equation*}
            \partial\I_C(\bar x) = \set{x^*\in X^*: \inner{x^*}{\tilde x-\bar x}\le 0 \text{ for all } \tilde x \in X \text{ with } A\tilde x = f}.
        \end{equation*}
        Now first let $x^* \in \Ker{A}^\perp$ and  fix $\tilde x$ such that $A\tilde x = f$. Now, $\bar x \in C$ implies $A(\tilde x-\bar x) = 0$, i.e. $\tilde x-\bar x \in \Ker{A}$. Thus, it must be $\inner{x^*}{ \tilde x-\bar x} = 0$ and therefore $x^*\in \partial \I_C(\bar x)$, which proves ${\Ker{A}^{\perp}} \subseteq \partial \I_C(\bar x)$. Conversely, let $x^*\in \partial \I_C(\bar x)$ and $\tilde x \in X$ such that $A \tilde x = f$, this means $\inner{x^*}{\tilde x-\bar x} \le 0$. But since $A(2\bar x-\tilde  x) = f$, it must also be $\inner{x^*}{\bar x-\tilde  x}=\inner{x^*}{(2\bar x-\tilde x)-{\bar{x}}} \le 0 $. That yields $\inner{x^*}{\bar x-\tilde x} = 0$, hence we have \begin{equation*}
            \partial \I_C(\bar x) = \set{x^*\in X : \inner{x^*}{\tilde x- \bar x} = 0\text{ for all } \tilde x \in {X} \text{ with } A\tilde x = f} 
        \end{equation*}
        Thus, the equality $C  =\set{x \in X: Ax = f}= \Ker{A} + \set{\bar x}$ implies $\partial \I_C(\bar x) \subseteq  \Ker{A}^\perp$.
    \end{proof}
\end{proposition}

This allows to show well-definedness of $p_l$.

\begin{lemma}\label{lemma:existence_subgradient_noisefree}
   Assume  that there is a solution of $Ax =f$ in the interior of $\dom g\square h$.  Furthermore, assume the infimal convolution $J_l=D_g^{p_{l-1}}(\cdot,u_{l-1})\square h$ (with the convention $p_{0} =0$ and $D_g^{p_0}(u,u_0) = g(u)$) is exact for all indices $l\in \N$. Then for any $l\in\N$, there is $p_l \in \partial g(u_l)\cap \Ker A ^\perp$ such that $p_l-{p_{l-1}} \in  {\partial h(v_l)}$.
    \begin{proof}
Since $\dom D_g^{p_{l-1}}(\cdot,u_l) =\dom g$, there must be a solution of $Ax =f$ in the interior of $\dom J_l$. Recall that $(u_l,v_l)$ satisfies\begin{equation*}
    (u_l,v_l) \in \argmin\limits_{A(u+v) =f} D_g^{p_{l-1}}(u,u_{l-1}) + h(v).
\end{equation*} Setting $x_l = u_l+v_l$ and using the exactness of $J_l$, the optimality condition of the previous equation becomes  \begin{equation*}
    0 \in \partial\pars{J_l+\I_{{\{Ax = f\}}}}(x_l).
\end{equation*} 
Since {$\dom J_{l}$} contains an inner point solving $Ax =f$, we obtain by Proposition \ref{prop:properties_subdifferentials} (ii)  that 
\begin{equation*}
    \partial\pars{J_l+\I_{\{Ax = f\}}}(x_l)= \partial J_l(x_l) + \partial{\I_{\{Ax = f\}}}(x_l) = \partial J_l(x_l) +\Ker A ^\perp.
\end{equation*}
In particular, this means $\partial J_l(x_l)  \cap \Ker{A}^\perp \neq \emptyset$. By Proposition \ref{prop:inf_conv_properties} (iv) in this setting, one has $\partial J_{l}(x_l) = \partial D_g^{p_{l-1}}(\cdot,u_{l-1})(u_l) \cap \partial h(v_l)$. Therefore $\pars{\partial g(u_l)-\set{p_{l-1}}}\cap\partial h(v_l) \cap\Ker A^\perp \neq 0$. This means there is $p_l \in \partial g(u_l)$ such that ${p_l-p_{l-1}} \in \partial h(v_l)$ and $p_l-p_{l-1} \in \Ker A^\perp$. Hence, by induction it also follows $p_l \in \Ker A^\perp$.
    \end{proof}
\end{lemma}

\begin{remark}\label{rem:update_differrentiable}
The necessary subgradient $p_l$ can be computed in a simple way when $g$ or $h$ is differentiable. If $g$ is differentiable, then $\partial g(u) =\set{\nabla g(u)}$ for all $u\in \dom g$ and thus $p_l = \nabla g(u_l)$. Otherwise, for $h$ being differentiable, one has $p_l = \nabla h(v_l)+p_{l-1}$.
\end{remark}

Before we analyze the convergence of Algorithm {\ref{algo:NestedBregman_Noiseless}}, let us introduce the algorithm for the case of noisy data. We will then provide a combined analysis for both situations. 

\subsection{The case of noisy data}
Assume that instead of the true data $f$, we only have access to a noise corrupted measurement $f^\delta$. We consider additive noise $\eta := f^\delta-f$ that  satisfies an estimate of the form 
\begin{equation}\label{eq:noise_level}
    \norm{\eta} \le \delta
\end{equation} for some known $\delta >0$. In this case, it does not make sense to solve an optimization problem of the form \eqref{eq:NestedBregman_noiseless_step} with $f$ replaced by $f^\delta$. Thus, we present $2$ different approaches on how to adapt Algorithm \ref{algo:NestedBregman_Noiseless} to this scenario. 

\subsubsection{Bregman iterated Morozov regularization}
Recall that Bregman iterations for noisy data are stopped  according to the discrepancy principle \eqref{eq:discrepancy_principle}. This is due to the fact that it does not make sense to approximate $f^\delta$ any further than  the accuracy of the estimated noise level. We therefore relax the constraint $A(u+v) =f$ to $\norm{A(u+v)-f^\delta}\le \delta$, which means \eqref{eq:NestedBregman_noiseless_initial} and \eqref{eq:NestedBregman_noiseless_step} are replaced by a Morozov regularization  with regularizer $D_g^{p_{l-1}}(\cdot,u_{l-1})\square h$. Hence, we obtain the following procedure. 

\begin{algorithm}[H]
  \begin{algorithmic}[1]
  \State Compute \begin{equation}\label{eq:NestedBregman_Morozov_initial}
      (u_{1}^\delta,v_{1}^\delta) \in \argmin\limits_{\norm{A(u+v)- f^\delta} \le  \delta} \set{g(u) +h(v)}
  \end{equation}
  \State Choose $p_{1}^\delta \in \partial g(u_1^\delta)\cap \partial h(v_1^\delta)\cap \ran{(A^*)}$ 
 \For{$l = 2,3,\dots$}
    \State Compute \begin{equation}\label{eq:NestedBregman_Morozov_step}
        (u_l^\delta,v_l^\delta) \in \argmin\limits_{\norm{A(u+v) - f^\delta}\le \delta} D_g^{p_{l-1}^\delta}(u,u_{l-1}^\delta) + h(v) 
    \end{equation}
    \State Choose $p_{l}^\delta\in \partial g(u_l^\delta)\cap \ran{(A^*)}$ such that $p_l^\delta-p_{l-1}^\delta \in \partial h(v_l^\delta)$.
  \EndFor
  \end{algorithmic}
  \caption{Nested Bregman with Morozov Regularization}\label{algo:NestedBregman_Noisy_Morozov}
\end{algorithm}

Again, we postpone the proof of existence for the minimizers in \eqref{eq:NestedBregman_Morozov_initial} and \eqref{eq:NestedBregman_Morozov_step} to the actual examples we consider. Thus, we begin by showing the existence of subgradients with the desired properties. Additionally, we assume that there is $x\in X$ with $\norm{Ax-f^\delta} < \delta$. If this is not the case, Algorithm \ref{algo:NestedBregman_Noisy_Morozov} reduces to Algorithm \ref{algo:NestedBregman_Noiseless} in the sense that the constraint $\norm{Ax -f^\delta}\le \delta$ would coincide with  $Ax  = f$. Indeed, assume $\norm{Ax-f^\delta}\ge \delta$ for all $x$ and there exists  $\tilde x$, $\tilde f \neq f$ with $A\tilde x = \tilde f$ such that $\norm{A\tilde x -f^\delta} = \delta$, then for any $t \in (0,1)$ it would be $\delta ^2 \le \norm{A\pars{tx^\dagger +(1-t)\tilde x} -f^\delta}^2 < \delta^2$. Therefore, the constraint $\norm{A(u+v) - f^\delta}\le \delta$ becomes $A(u+v) = f$. Furthermore, due to the regularizing nature of problems \eqref{eq:NestedBregman_Morozov_initial} and \eqref{eq:NestedBregman_Morozov_step}, no additional assumption that guarantees the sum-rule for subgradients is needed for the existence of $p_l^\delta$.

\begin{lemma}\label{lemma:existence_subdifferential_morozov}
    Let $l \in \N$ and assume the infimal convolution $J_l=D_g^{p_{l-1}^\delta}(\cdot,u_{l-1}^\delta)\square h$ (with the convention $p_{0}^\delta =0$) is exact. If there is $x \in \dom \pars{g\square h}$ with $\norm{Ax-f^\delta} < \delta$, then there exists $p_l^\delta \in \partial g(u_l^\delta)\cap\ran{(A^*)}$ such that $p_l^\delta-p_{l-1}^\delta\in \partial h(v_l^\delta)$. Additionally, for each $l \in \N$ it is \begin{equation}\label{eq:subdifferential_in_normal_cone}
        \inner{p_{l-1}^\delta-p_l^\delta}{v-v_l^\delta}\le 0
    \end{equation} for all $v \in X$ with $\norm{A(u_l^\delta+v)-f^\delta}\le \delta$.
    \begin{proof}
  
         Define $J_l(x) = (D_g^{p_{l-1}^{\delta}}(\cdot,u_{l-1}^\delta)\square h)(x)$. Since $\dom g = \dom D_g^{p_{l-1}^{\delta}}(\cdot,u_{l-1}^\delta)$, it follows that $x \in \dom J_l$ with $\norm{Ax-f^\delta }<\delta$. Hence, the Slater condition for problem \eqref{eq:NestedBregman_Morozov_initial}, respectively $\eqref{eq:NestedBregman_Morozov_step}$ is satisfied, and we can apply Theorem 3.9 in \cite{barbu2012convexity} to find that $x_{l}^\delta := u_l^\delta+v_l^\delta$ solves \begin{equation}\label{eq:Lagrange_multiplier}
            \min\limits_{x \in X}\set{ \frac{1}{2}\norm{Ax-f^\delta}^2 +\alpha_l J_l(x)}
        \end{equation}
        for some $\alpha_l>0$ (see also Remark 2.6 in \cite{Lorenz_2013}). Therefore $\frac{1}{\alpha_l}A^*(f^\delta-Ax_l^\delta) \in \partial J_l(x_l^\delta)$. Since $\partial J_l(x_l^\delta) = (\partial g(u_l^\delta) -\set{p_{l-1}^\delta})\cap \partial h(v_l^\delta)$, we inductively obtain that the choice \begin{equation*}
            p_l^\delta = p_{l-1}^\delta + \frac{1}{\alpha_l}A^*(f^\delta-A(u_l^\delta+v_l^\delta)) = \sum_{i = 1}^l \frac{1}{\alpha_i}A^*(f^\delta-A(u_i^\delta+v_i^\delta)) 
        \end{equation*} satisfies $p_l^\delta \in \partial g(u_l^\delta)$ and $p_l^\delta -p_{l-1}^\delta \in \partial h(v_l^\delta)$.
\\ In order to show \eqref{eq:subdifferential_in_normal_cone}, let  $v \in X$ with $\norm{A(u_l^\delta+v)-f^\delta}\le \delta$. We distinguish $2$ cases. If ${\norm{A(u_l^\delta+v_l^\delta)-f^{\delta}}} = \delta$, it is \begin{align*} 
            \inner{p_{l-1}^\delta-p_l^\delta}{v-v_l^\delta} &= \frac{1}{\alpha_l}\inner{A^*(A(u_l^\delta+v_l^\delta)-f^\delta)}{(v+u_l^\delta)-(v_l^\delta+u_l^\delta)}\\&=
            \frac{1}{\alpha_l}  \inner{A(u_l^\delta+v_l^\delta)-f^\delta}{A(v+u_l^\delta)-f^\delta} \\&-\frac{1}{\alpha_l} \inner{A(u_l^\delta+v_l^\delta)-f^\delta}{A(u_l^\delta+v_l^\delta)-f^\delta}           
            \\&\le\frac{1}{\alpha_l} \norm{A(u_l^\delta+v_l^\delta)-f^\delta}\pars{\norm{A(u_l^\delta+v)-f^\delta} - \norm{A(u_l^\delta+v_l^\delta)-f^\delta}} \\&\le 0.
        \end{align*}        
        Otherwise, if $\norm{A(u_l^\delta+v_l^\delta)-f^\delta}<\delta$, we must have that $x_l^\delta \in \argmin\limits_{\norm{Ax-f^\delta}\le \delta}J_l(x)$ is a local minimum of $J_l$. By convexity of $J_l$, it must also be a global minimum of $J_l$, so that $p_l^\delta = p_{l-1}^\delta$ is a feasible choice and the iteration becomes stationary. In this case, the claim follows inductively.
    \end{proof}
\end{lemma}

With this, we can  establish convergence properties of the Nested Bregman algorithm. The idea of the proof follows the convergence analysis in \cite{bur-osh} (see also Section 6.1 in \cite{Benning-Burger}). Further, we use the same notation for both Algorithm \ref{algo:NestedBregman_Noiseless} and Algorithm \ref{algo:NestedBregman_Noisy_Morozov}, where the noise-free case corresponds to $\delta = 0$. {For notational consistency, we denote $u_0^\delta = u_0$ and $p_0^\delta = p_0 = 0$, so that it is $g(u) = D_g^{p_0}(u,u_0)  = D_g^{p_0^\delta}(u,u_0^\delta) $ for all $u \in \dom g.$}
\begin{theorem}\label{thm:Convergence_Nested_Bregman-Morozov}
Let $(u_l^\delta,v_l^\delta)_{l \in \N}$ be the sequence generated by Algorithm \ref{algo:NestedBregman_Noiseless} or Algorithm \ref{algo:NestedBregman_Noisy_Morozov}. Under the assumptions of Lemma \ref{lemma:existence_subgradient_noisefree} for Algorithm \ref{algo:NestedBregman_Noiseless} or Lemma \ref{lemma:existence_subdifferential_morozov} for Algorithm \ref{algo:NestedBregman_Noisy_Morozov}, the following hold.
    \begin{enumerate}
                \item The sequence $(h(v_l^\delta))_{l \in \N}$ is decreasing.
        \item Assume that $g(\bar x) <\infty$ for a solution $\bar x$ of $Ax = f$. Then the estimate \begin{equation}\label{eq:Convergence_h_noisy}
            h(v_l^\delta) \le h(0) + \frac{g(\bar x)}{l}.  
        \end{equation}
        holds for any $l \in \N$. In particular, if ${h(0)= \inf h =\min h} $, this implies that $h(v_l^\delta)$ converges to $\min h$ {as $l \to \infty$}. 
    \end{enumerate}
    \begin{proof}
                \begin{enumerate}
        \item Let $l \ge 2$. Since the pair $(u_{l-1}^\delta,v_{l-1}^\delta)$ satisfies $\norm{A(u_{l-1}^\delta+v_{l-1}^\delta) -f^\delta}\le \delta$, it holds \begin{equation*}
            h(v_l^\delta) \le h(v_l^\delta) + D_g^{p_{l-1}^\delta}(u_l^\delta,u_{l-1}^\delta) \overset{\eqref{eq:NestedBregman_Morozov_step}}{\le} h(v_{l-1}^\delta) + D_g^{p_{l-1}^\delta}(u_{l-1}^\delta,u_{l-1}^\delta) = h(v_{l-1}^\delta).
        \end{equation*}
        \item Let $k \in \N$ and let $v = \bar x -u_k^\delta$. That means $\norm{A(u_k^\delta+v)-f^\delta} \le \delta$. Hence, for $\delta  = 0$, we have $A(u_k^\delta +v) = f$. Then $p_k^\delta \in  \Ker A^\perp$ implies $\inner{p_{k-1}^\delta-p_k^\delta}{\bar x -u_k^\delta -v_k^\delta}= 0$. Otherwise, for $\delta >0$ we use \eqref{eq:subdifferential_in_normal_cone} to see that
        \begin{equation}\label{eq:help_convergence_1}
            \inner{p_{k-1}^\delta-p_k^\delta}{\bar x -u_k^\delta -v_k^\delta} = \inner{p_{k-1}^\delta-p_k^\delta}{v-v_k^\delta}\le 0.
        \end{equation}
        Thus, using the three-point identity for Bregman distances, we have {for $k \ge 2$} \begin{align*}
            D_g^{p_{k}^\delta}(\bar x ,u_k^\delta)-D_g^{p_{k-1}^\delta}(\bar x ,u_{k-1}^\delta)&+D_g^{p_{k-1}^\delta}(u_k^\delta,u_{k-1}^\delta) = \inner{p_k^\delta-p_{k-1}^\delta}{u_k^\delta-\bar x} \\&= \inner{p_k^\delta-p_{k-1}^\delta}{0-v_k^\delta} + \inner{p_{k-1}^\delta-p_k^\delta}{\bar x -u_k^\delta -v_k^\delta}\\&\overset{\eqref{eq:help_convergence_1}}{\le} \inner{p_k^\delta-p_{k-1}^\delta}{0-v_k^\delta} \\&
            \le h(0)- h(v_k^\delta).
        \end{align*}
        Rearranging gives \begin{align}\label{eq:help_convergence_2}
            h(v_k^\delta)-h(0) &\le D_g^{p_{k-1}^\delta}(\bar x ,u_{k-1}^\delta)- D_g^{p_{k}^\delta}(\bar x ,u_k^\delta)-D_g^{p_{k-1}^\delta}(u_k^\delta,u_{k-1}^\delta)\notag\\&\le D_g^{p_{k-1}^\delta}(\bar x ,u_{k-1}^\delta)- D_g^{p_{k}^\delta}(\bar x ,u_k^\delta).
        \end{align}

        {For $k =1$, we analogously obtain \begin{align*}
    D_g^{p_1^\delta}(\bar x, u_1^\delta) -g(\bar x) + g(u_1^\delta) &= \inner{p_1^\delta}{u_1^\delta-\bar x} \\&= \inner{p_1^\delta}{0-v_1^\delta} - \inner{p_1^\delta}{\bar x - (u_1^\delta+v_1^\delta)} \\&\overset{\eqref{eq:help_convergence_1}}{\le} \inner{p_1^\delta}{0-v_1^\delta} \\&\le h(0) -h(v_1^\delta).
\end{align*}
Recall, that we denote $D_g^{p_0^\delta}(\cdot,u_0^\delta) = g$. Thus, by employing the non-negativity of $g$ we thus obtain \begin{equation}\label{eq:helo_convergnce_2_1}
            h(v_1^\delta)-h(0) \le g(\bar x) - D_g^{p_1^\delta}(\bar x,u_1^\delta) = D_g^{p_0^\delta}(\bar x,u_0^\delta) - D_g^{p_1^\delta}(\bar x,u_1^\delta).
        \end{equation} }
        Together with the monotonicity of $h(v_k^\delta)$, {summing \eqref{eq:help_convergence_2} and \eqref{eq:helo_convergnce_2_1}}  yields for any $l\in \N$\begin{align*}
            l\pars{h(v_l^\delta) -h(0)} &\overset{\eqref{eq:help_convergence_2}, \eqref{eq:helo_convergnce_2_1} }{\le} \sum\limits_{k = 1}^{l} \pars{D_g^{p_{k-1}^\delta}(\bar x,u_{k-1}^\delta) - D_g^{p_{k}^\delta}(\bar x ,u_{k}^\delta)} \\&{= D_g^{p_0^\delta}(\bar x,u_0^\delta) - D_g^{p_l^\delta}(\bar x,u_l^\delta)}  \\&= g(\bar x) -D_g^{p_l^\delta}(\bar x,u_l^\delta) \\&\le g(\bar x).
        \end{align*}
        Estimate \eqref{eq:Convergence_h_noisy} follows.
        \end{enumerate}
    \end{proof}
\end{theorem}

{\begin{remark}
    For the convergence of $(h(v_l^\delta))_{l \in \N}$ to $\inf h$, one does not need that the minimum of $h$ is attained at $0$. If there is a minimizing sequence $\seq{v}$ of $h$ and a sequence $\seq{x}$ of solutions to $Ax = f$ such that $g(x_n-v_n)$ is finite for all $n \in \N$, we also obtain $\lim\limits_{n \to \infty} h(v_l^\delta) = \inf h$. Indeed, define the functions $h_n := h(\cdot+v_n)$ and $g_n:=g(\cdot-v_n)$. Then the pair $(u_l^\delta+v_n,v_l^\delta-v_n)$ is a solution of \eqref{eq:NestedBregman_Morozov_initial} for $l = 1$ or \eqref{eq:NestedBregman_Morozov_step} for $l \ge 2$ with $g_n$ and $h_n$ instead of $g$ and $h$ respectively. Thus, by Theorem \ref{thm:Convergence_Nested_Bregman-Morozov} it holds \begin{equation*}
        h(v_l^\delta) = h_n(v_l^\delta-v_n) \le h_n(0) + \frac{g_n(x_n)}{l} = h(v_n) + \frac{g(x_n-v_n)}{l}.
    \end{equation*}
    Letting $l \to \infty$ gives \begin{equation*}
        \limsup\limits_{l \to \infty}h(v_l^\delta) \le h(v_n)
    \end{equation*} for all $n \in \N$ and thus $\lim\limits_{l \to \infty}h(v_l^\delta) = \inf h$.
\end{remark}}

Having established the convergence of $h(v_n^\delta)$, one can see how Algorithm \ref{algo:NestedBregman_Noisy_Morozov} can be used to find good components $u$ and $v$ which also solve $\norm{A(u+v)-f^\delta} \le \delta$. Assume that in \eqref{eq:decomposition_problem} the parameter $\beta$ is chosen too small, i.e. that the component $v$ is overrepresented in the decomposition. In this case,  we expect the share of the component $v$ to decrease within the iteration. Since the constraint $\norm{A(u+v)-f^\delta}\le \delta$ ensures that $u_l^\delta+v_l^\delta$ is a reasonable approximation to the true solution of \eqref{eq:ill_posed_problem}, the share of the component $u$ should be increasing. This means the decomposition should change from $v$ being over-regularized to $u$ being over-regularized. Hence, instead of choosing appropriate parameters in \eqref{eq:variational_problem_J}, we now only need to stop the iteration according to some meaningful rule. As one example for such rule, we will use a cross-correlation based stopping criterion in the numerical experiments in Section \ref{sec:numerical_results}.

\subsubsection{Nested Bregman iterations stopped by the discrepancy principle}
Instead of replacing the minimization problems \eqref{eq:NestedBregman_noiseless_initial} and \eqref{eq:NestedBregman_noiseless_step} by Morozov regularization, we also propose to run Bregman iterations until the discrepancy principle is satisfied. Hence, we obtain the following algorithm.

\begin{algorithm}[H]
  \begin{algorithmic}[1]
    \State Set $\tilde u_0^\delta = 0$, $\tilde v_0^\delta = 0$,  $\tilde p_0^\delta = 0$, $k = 0$.
    \State Define $J_1(u,v) = g(u) +h(v)$.
    \While {$\norm{f^\delta - A(\tilde u_k^\delta+\tilde v_k^\delta)} > \tau \delta$}
        \State Set $k \leftarrow k+1$.
        \State Compute $(\tilde u_k^\delta,\tilde v_k^\delta) \in \argmin\limits_{u,v \in X} \frac{1}{2}\norm{f^\delta-A(u+v)}^2 + D_{J_1}^{\tilde p_{k-1}^\delta}\pars{(u,v),(\tilde u_{k-1}^\delta,\tilde v_{k-1}^\delta)} $
        \State Set $\tilde p_k^\delta = \tilde p_{k-1}^\delta + A^*(f^\delta - A(\tilde u_k^\delta+v_k^\delta))$.
    \EndWhile
    \State Set $u_1^\delta = \tilde u_k^\delta$, $v_1^\delta = \tilde v_k^\delta$, $p_1^\delta = \tilde p_k^\delta$.
\For{$l= 2,3,\dots$}
        \State Set $\tilde u_0^\delta = 0$, $\tilde v_0^\delta = 0$,  $\tilde p_0^\delta = 0$, $k = 0$.
        \State Define $J_l(u,v) = D_g^{p_{l-1}^\delta}(u,u_{l-1}^\delta) +h(v)$.
        \While {$\norm{f^\delta - A(\tilde u_k^\delta+\tilde v_k^\delta)} > \tau \delta$}
        \State Set $k \leftarrow k+1$.
        \State Compute $(\tilde u_k^\delta,\tilde v_k^\delta) \in \argmin\limits_{u,v \in X} \frac{1}{2}\norm{f^\delta-A(u+v)}^2 + D_{J_l}^{\tilde p_{k-1}^\delta}\pars{(u,v),(\tilde u_{k-1}^\delta,\tilde v_{k-1}^\delta)} $
        \State Set $\tilde p_k^\delta = \tilde p_{k-1}^\delta + A^*(f^\delta - A(\tilde u_k^\delta+v_k^\delta))$.
    \EndWhile
    \State Set $u_l^\delta = \tilde u_k^\delta$, $v_l^\delta = \tilde v_k^\delta$, $p_l^\delta = p_{l-1}^\delta + \tilde p_k^\delta$.
    \EndFor
  \end{algorithmic}
  \caption{Nested Bregman for noisy data}\label{algo:NestedBregman_Noisy}
\end{algorithm}

This algorithm consists of a sequence of inner loops within an outer loop. In the $l-$th inner loop, we run Bregman iterations with regularizer $J_l = D_g^{p_{l-1}^\delta}(\cdot,u_{l-1}^\delta)\square h$. At the beginning of each outer iteration, we then replace $g$ by the Bregman distance of $g$ at the iterate where the discrepancy principle was met. This comes with several advantages. The properties $p_l^\delta \in \partial g(u_l^\delta) \cap \ran{(A^*)}$ and $p_l^\delta -p_{l-1}^\delta \in \partial h(v_l^\delta)$ follow immediately from the algorithm, which avoids the computation of a suitable subdifferential as compared to Algorithms \ref{algo:NestedBregman_Noiseless} and \ref{algo:NestedBregman_Noisy_Morozov}. Additionally, this algorithm can easily be implemented for other data fidelity terms than $\norm{\cdot}^2$, as there is no constraint to be taken care of. However, we can not simply extend the convergence guarantees from Algorithm \ref{algo:NestedBregman_Noisy_Morozov} to this scenario, because the iterates of the inner loops at the discrepancy principle do not necessarily solve a minimization problem of the form \eqref{eq:decomposition_problem_J}. Thus, monotonicity of $h(v_l^\delta)$ may not be established. Moreover, we need to solve more minimization problems for this procedure than for Algorithm \ref{algo:NestedBregman_Noisy_Morozov}, resulting in a larger computational effort. However, because in general the iterates of Bregman iterations stopped via the discrepancy principle tend to be superior to the ones obtained from variational regularization (see the numerical examples in Section \ref{sec:Bregman}, as well as the discussion in \cite{Benning-Burger}), this method is still of interest. In terms of improving a decomposition with an over-regularized component, the method produces comparable results to Algorithm \ref{algo:NestedBregman_Noisy_Morozov}, as can be seen in Section \ref{sec:numerical_results}.

\section{Well-definedness for selected regularizers}\label{sec:selected_regularizers}
In this part, we give an overview of the infimally-convoluted regularizers which we use in the numerical section. We show exactness of the infimal convolutions and well-definedness of the minimization problems from Algorithms \ref{algo:NestedBregman_Noiseless} and \ref{algo:NestedBregman_Noisy_Morozov} in the infinite dimensional setting. As opposed to the situation of Remark \ref{rem:Bregman_well_definedness},  problems \eqref{eq:NestedBregman_noiseless_step} and \eqref{eq:NestedBregman_Morozov_step} are not equivalent to problems of the form \eqref{eq:NestedBregman_noiseless_initial} or \eqref{eq:NestedBregman_Morozov_initial}, because we update the  Bregman distance of one functional within an infimal convolution, rather than the Bregman distance of the entire infimal convolution.
 In fact, coercivity of an infimal convolution is not necessarily preserved if one of the functionals is replaced with a Bregman distance. As a simple counterexample, let $X = \R$ and $g = h = \abs{\cdot}$. Then for any $x \neq 0$ and $\lambda \in \R$, it is $g\square h(\lambda x) = \abs{\lambda}\abs{x}$. However, for $p = \frac{x}{\abs{x}} \in \partial g(x)$, we have $D_g^p(\lambda x,x) = 0$ for all $\lambda \ge 0$. Hence, $\lim\limits_{\lambda \to \infty}g\square h(\lambda x) = \infty$, while $\lim\limits_{\lambda \to \infty} (D_g^p(\cdot,x)\square h)(\lambda x) = 0$. However, it is possible to show that any coercivity result for $g\square h$ within the feasible sets $\set{(u,v): A(u+v) = f}$ and $\set{(u,v): \norm{A(u+v)-f^\delta}\le \delta}$ can be extended to $D_g^{p_l^\delta}(\cdot,u_l^\delta)\square h$. Thus, the constraints in the minimization problems ensure their well-definedness. For brevity, we give a combined proof for the cases of exact and noisy data, where the exact data case corresponds to $\delta = 0$.

\begin{lemma}\label{lemma:coercivity_bregmanized}
    For $l \in \N$, let $p_l^\delta$ be obtained as in Algorithm \ref{algo:NestedBregman_Noisy_Morozov}, or Algorithm \ref{algo:NestedBregman_Noiseless} in the case $\delta = 0$ (but still denoted by $p_l^\delta)$ and assume that $f \in \ran (A)$.   Assume further that $h$ satisfies 
    \begin{equation}\label{eq:assumptions_coercivity_h}
       h(-x) = h(x)  {\text{ for all }} x \in X.    \end{equation}
Let $((u_n,v_n))_{n\in \N}$ be a sequence in $X\times X$ such that $\norm{A(u_n+v_n)-f^\delta}\le \delta$ for all $n\in\N$ and $\lim\limits_{n \to \infty}g(u_n) +h(v_n) = \infty$, then $\lim\limits_{n \to \infty}D_g^{p_l^\delta}(u_n,u_l^\delta) +h(v_n) = \infty$.
\begin{proof}
{Recall that we assume $g(u_0) = 0$, so that upon setting $u_0^\delta = u_0$ and $p_0^\delta = 0$, it is $D_g^{p_0}(u,u_0^\delta) = g(u)$ for all $u \in \dom g$.}
 {Fix} $x \in X$ with $Ax = f$. If $\delta >0$, we have $p_l^\delta \in \ran({A^*})$ by Lemma \ref{lemma:existence_subdifferential_morozov}. Hence, there is $q_l^\delta \in Y$ with $A^* q_l^\delta = p_l^\delta$. Therefore, \begin{align*}
     \inner{p_l^\delta}{u_n} &= \inner{q_l^\delta}{Au_n} = \inner{q_l^\delta}{A(u_n+v_n)-f} + \inner{q_l^\delta}{f-Av_n} \\&= \inner{q_l^\delta}{A(u_n+v_n)-f} + \inner{p_l^\delta}{x-v_n}.
 \end{align*}
This implies \begin{equation}\label{eq:Estimate_Bregman_Distance_Coercivity_improved1}
      -\inner{p_l^\delta}{u_n} \ge \inner{p_l^\delta}{v_n-x}-\delta\norm{q_l^\delta}.
 \end{equation}
 Due to \eqref{eq:assumptions_coercivity_h} and $(p_l^\delta-p_{l-1}^\delta) \in \partial h(v_l^\delta)$, it is $(p_{l-1}-p_l)\in \partial h(-v_l^\delta)$ and we deduce that \begin{equation*}
     \inner{p_l^\delta}{v_n} = \inner{p_l^\delta-p_{l-1}^\delta}{v_n} + \inner{p_{l-1}^\delta}{v_n}\ge h(-v_l^\delta) -h(v_n) + \inner{p_{l-1}^\delta-p_l^\delta}{v_l^\delta} + \inner{p_{l-1}^\delta}{v_n}.
 \end{equation*} 
 Iterating this estimate, we obtain \begin{equation}\label{eq:Estimate_Bregman_Distance_Coercivity_improved2}
     \inner{p_l^\delta}{v_n} \ge \sum\limits_{i = 1}^l\pars{h(-v_i^\delta)+\inner{p_{i-1}^\delta-p_i^\delta}{v_i^\delta}} - lh(v_n).
 \end{equation}
  Denote \begin{equation*}
      C_l = -g(u_l^\delta) -\sum\limits_{i = 1}^l\pars{h(-v_i^\delta)+\inner{p_{i-1}^\delta-p_i^\delta}{v_i^\delta}} -\inner{p_l^\delta}{x}-\delta\norm{q_l^\delta}+\inner{p_l^\delta}{u_l^\delta}.
  \end{equation*} Then, by combining \eqref{eq:Estimate_Bregman_Distance_Coercivity_improved1} and \eqref{eq:Estimate_Bregman_Distance_Coercivity_improved2}, we obtain \begin{align*}
         D_g^{p_l^\delta}(u_n,u_l^\delta) +h(v_n) &\ge \frac{1}{l+1}D_g^{p_l^\delta}(u_n,u_l^\delta)+h(v_n)\\ &= \frac{1}{l+1}\pars{g(u_n) -g(u_l^\delta) -\inner{p_l^\delta}{u_n-u_l^\delta}} + h(v_n)\\&\ge \frac{1}{l+1}\pars{g(u_n) -lh(v_n)+ C_l} + h(v_n)\\&\ge \frac{1}{l+1}\pars{g(u_n) +h(v_n) +C_l}.
 \end{align*}
 Since $C_l$ depends only on $l$ (but not on $n$), the claim follows. For $\delta =  0$, we recall that the constraint means $A(u_n+v_n) = f$. Consequently, it is $x-(u_n+v_n)\in \Ker A$ for all $n\in \N$. Together with $p_l^\delta \in \Ker A^\perp$ we obtain \begin{equation*}
      \inner{p_l^\delta}{u_n} = \inner{p_l^\delta}{x-v_n}.
  \end{equation*} Hence, the previous estimates work analogously for the case $\delta = 0$ by just omitting the terms of the form $ \inner{q_l^\delta}{A(u_n+v_n)-f}$.
\end{proof}
\end{lemma}

\begin{remark}\label{rem:coercivity_for_other_f}
    Upon close inspection of the proof of Lemma \ref{lemma:coercivity_bregmanized}, we notice that the result also holds for any sequence $(u_n,v_n)$ satisfying $\norm{A(u_n+v_n)-\tilde f^\delta}\le \tilde \delta$ for another $\tilde f^\delta$ and $\tilde \delta$ than the ones used in \eqref{eq:NestedBregman_Morozov_initial} and \eqref{eq:NestedBregman_Morozov_step}. This means that, as long as we assume that $u,v$ are such that $A(u+v)$ is bounded, we can extend coercivity of $g\square h$ to the infimal convolution $D_g^{p_l^\delta}\square h$. This will be helpful in showing its exactness later.
\end{remark}

In what follows, we prove exactness of the infimal convolutions and well-definedness of the minimizing problems in Algorithm \ref{algo:NestedBregman_Noisy_Morozov} for selected regularizers in the infinite dimensional function space setting. We would like  to point out that, technically, we only need to show that the infimal convolutions are exact at the minimizers of the constraint problem. By definition, this is the case as soon as those minimizers exist. Nonetheless, extending exactness of the infimal convolution for all $x\in X$ is interesting in itself and helps simplify the proofs of the existence of minimizers. The proofs mostly use standard coercivity arguments for the corresponding problems with $g\square h$ instead of $D_g^{p_l^\delta}(\cdot,x_{l})\square h$. We then make use of Lemma \ref{lemma:coercivity_bregmanized} to extend the results to our setting. The regularizers we consider are merely exemplary and shall illustrate the variety of infimal convolution based decompositions and the applicability of our proposed methods.

\subsection{\texorpdfstring{$L^1-H^1$}{TEXT}}
 Let $\Omega \subseteq \R^d$ be a bounded Lipschitz domain. In order to separate peaks from an otherwise smooth signal, we use the infimal convolution of $g = \frac{\alpha}{2}\norm{\nabla\cdot}_{L^2}^2$  and $h = \beta\norm{\cdot}_{L^1}$.
 Hence, let $X = Y = L^2(\Omega)$ . We show that the minimization problems in Algorithm \ref{algo:NestedBregman_Noisy_Morozov} are well-defined for the case of denoising. Note that  the first step of our algorithm corresponds to denoising with regularizer $g\square h$, for which a proof of well-definedness might exist in literature. Since we could not find one, we outline the arguments below and illustrate how Lemma \ref{lemma:coercivity_bregmanized} helps to extend the well-definedness to problems of the form \eqref{eq:NestedBregman_Morozov_step}. We first recall the space of functions of bounded variation. \begin{definition}\cite{AmbrosioBV}
     The (isotropic) total variation of a function $u \in L^1(\Omega)$ is defined as \begin{equation*}
        \tv{u} = \sup\set{\int\limits_\Omega u\div  \varphi : \ \varphi \in C_0^\infty(\Omega,\R^2), \;\norm{\varphi(x)}_2 \le 1 \text{ for all }x\in \Omega}, 
\end{equation*}
 where $\norm{\cdot}_2$ denotes the Euclidean norm on $\R^d$. {It holds that \begin{equation*}
     \tv{u} = \norm{Du}_{\mathcal{M}},
 \end{equation*} where} $Du$ {denotes} the finite Radon measure  that represents the distributional derivative of $u$ and $\norm{\cdot}_\mathcal{M}$ the Radon norm. Recall that, for a given Radon measure $\mu$ on $\Omega$, the Radon norm is defined as 
 \begin{equation*}  
   \norm{\mu}_{\mathcal{M}}=\sup\set{\int\limits_\Omega \varphi \,d\mu: \ \varphi \in C_0^\infty(\Omega,\R^2),\; \norm{\varphi(x)}_2 \le 1 \text{ for all }x\in \Omega}.
 \end{equation*}
 The domain of the total variation is then the space of functions of bounded variation $BV(\Omega) = \set{u \in L^1(\Omega) : \tv{u} < \infty}$. If equipped with the norm $\norm{u}_{BV} = \norm{u}_{L^1} + \tv{u}$, this space  is complete and continuously embedded in $L^\frac{d}{d-1}(\Omega)$ (for the case $d =1$ we use the convention $\frac{1}{0} = \infty$). Furthermore, if $u \in W^{1,1}(\Omega)$, then \begin{equation*}
     \tv{u} = \int_\Omega \norm{\nabla u}_2,
 \end{equation*}
 where $\nabla u$ denotes the weak gradient of $u$.
 \end{definition}
 
 \begin{lemma}\label{lemma:existence_L1_H1}
       Let $g(u) = \frac{\alpha}{2}\norm{\nabla u}_{L^2}^2$, $h(v) = \beta\norm{v}_{L^1}$   and $A = Id$. Additionally, assume that $f \in L^1(\Omega)$. Then the infimal convolutions $D_g^{p_{l}}(\cdot,u_{l}^\delta)\square h$ are exact for every $x \in X$, and the minimization problems \eqref{eq:NestedBregman_Morozov_initial} and \eqref{eq:NestedBregman_Morozov_step} are well-defined for every $l\in \N_0$.
       \begin{proof}

           Let $l \in \mathbb{N}_0$. For simpler notation, denote $g_l(\cdot):=D_g^{p_l^\delta}(\cdot,u_{l}^\delta)$ and $J_l = g_l\square h$, with the convention $g_l = g$ and $p_l^\delta = 0$ if $l = 0$  throughout the proof. We show the results in the case $d \ge 2$. For $d = 1$ the results can be proved analogously by replacing weak convergence in $L^\frac{d}{d-1}(\Omega)$ with weak*-convergence in $L^\infty(\Omega)$.
                      \begin{enumerate}
               \item We start by showing the exactness of the infimal convolution. 
               Let $x\in X$  and consider a minimizing sequence $(u_n,v_n)$ such that $u_n+v_n = x$ and $\lim\limits_{n \to \infty} g_l(u_n) +h(v_n) = J_l(x)$. Hence, the sequences $(g_l(u_n))_{n \in \N}$ and $(h(v_n))_{n \in \N}$ are individually {bounded}, since $g_l$ and $h$ are both non-negative. Furthermore, {$\norm{\nabla u_n}_{L^{2}}$} must be bounded. If it was unbounded, then $(\frac{\alpha}{2}\norm{\nabla u_n}_{L^2}^2 + \beta\norm{v_n}_{L^1})_{n\in \N}$ would be unbounded. But since $u_n+v_n = x$ for all $n \in \N$, we can apply Lemma \ref{lemma:coercivity_bregmanized} with $f^\delta$ replaced by $x$ and an arbitrary $\delta >0$ to find that that $(g_l(u_n) +h(v_n))_{n\in \N}$ would be unbounded. Therefore, $(\tv{u_n})_{n\in \N}$ is bounded since \begin{equation*}
                  \tv{u_n} =  \int\limits_\Omega \abs{\nabla u_n} \le \abs{\Omega}^\frac{1}{2}\pars{ \int\limits_\Omega\abs{\nabla u_n}^2}^\frac{1}{2}.
               \end{equation*}  
               Since $u_n = x-v_n$ and $x \in L^1(\Omega)$, we obtain that $(u_n)_{n\in \N}$ is bounded in $L^1(\Omega)$. Therefore, $(\norm{u_n}_{BV})_{n\in \N}$ is bounded. This implies $(\norm{u_n}_{L^\frac{d}{d-1}})_{n\in \N}$ is bounded, so that we can assume (upon passing to a subsequence, denoted the same) that $(u_n)_{n\in \N}$ converges weakly in $L^\frac{d}{d-1}(\Omega)$ to some $u^*$. Consequently, $(v_n)_{n\in \N}$ converges weakly in $L^\frac{d}{d-1}(\Omega)$ (on a subsequence) to $v^* = x-u^*$. If we can  show that  both $g_l$ and $h$ are lower semicontinuous with respect to $L^\frac{d}{d-1}(\Omega)$, then the claim follows by standard arguments. Indeed, first note that $p_l^\delta = \sum_{i = 1}^l(p_i^\delta-p_{i-1}^\delta)$ with $p_i^\delta-p_{i-1}^\delta \in \partial h(v_i^\delta)$ for all $i \in \set{1,\dots l}$. Because $\partial h(v_i^\delta) \subset \set{v \in L^2(\Omega) : \norm{v}_{L ^\infty} \le \frac{1}{\beta}} \subset L^\infty(\Omega)$, we have $p_l^\delta \in L^\infty(\Omega)$. Therefore, the duality paring $\inner{p_l^\delta}{\cdot}$ is (weakly) continuous on $L^\frac{d}{d-1}(\Omega)$. Next, note that for any $u \in L^\frac{d}{d-1}(\Omega)$ it is \begin{equation*}
                   \norm{\nabla u}_{L^2} = \sup\limits_{\substack{\varphi \in C_0^\infty(\Omega,\R^2)\\ \norm{\varphi}_{L^2} \le 1}} \int\limits_\Omega  \nabla u \cdot \varphi = \sup\limits_{\substack{\varphi \in C_0^\infty(\Omega,\R^2)\\ \norm{\varphi}_{L^2} \le 1}} \int\limits_\Omega  u\, \div \varphi.
               \end{equation*}
               Thus, $\norm{\nabla \cdot}_{L^2}$ is a supremum of continuous functionals  on $L^\frac{d}{d-1}(\Omega)$, which implies lower-semicontinuity. In summary, $g_l$ is lower semicontinuous and convex on $L^\frac{d}{d-1}(\Omega)$. Similarly, we can write \begin{equation*}
                   \norm{v}_{L^1} = \sup\limits_{\substack{\varphi \in C_0^\infty(\Omega,\R)\\ \norm{\varphi}_{L^\infty}\le 1}} \int\limits_\Omega {v} \,  \varphi.
               \end{equation*}
                Again, this implies that $\norm{ \cdot}_{L^1}$ is lower semicontinuous with respect to $L^\frac{d}{d-1}(\Omega)$. {It remains to show $u^* \in L^2(\Omega)$. Due to $ \norm{\nabla u^*}_{L^2}< \infty$ and $u^* \in L^\frac{d}{d-1}(\Omega)$, it follows that $u\in L^2(\Omega)$ as outlined next. Since $\Omega$ is a bounded, Lipschitz domain, there is a sequence $\seq{\bar u} \subseteq C^\infty (\bar \Omega)$, such that $\seq{\bar u}$ converges to $u^*$ in $L^\frac{d}{d-1}(\Omega)$ and $(\nabla \bar u_n)_{n \in \N}$ converges to $\nabla u^*$ in $L^2(\Omega)$. In particular, we obtain that  $\int_\Omega \bar u_n$ converges to $\int_\Omega u^*$, so that the Poincaré-Wirtinger inequality in $W^{1,2}(\Omega)$ implies \begin{equation*}
    \norm{u^*-\int_\Omega u^*}_{L^2} \le \liminf\limits_{n \to \infty} \norm{\bar u_n -\int_\Omega \bar u_n}_{L^2} \le C\liminf\limits_{n \to \infty} \norm{\nabla \bar u_n}_{L^2} < \infty.
\end{equation*}
Since $x \in L^2(\Omega)$, this also implies $v^* = x-u^* \in L^2(\Omega)$.}

               \item Let $(x_n)_{n\in \N}$ be a sequence with $\norm{x_n -f^\delta}_{L^2} \le \delta$ and \begin{equation*}
                   \lim\limits_{n \to \infty}J_l(x_n) = \inf\limits_{\norm{x-f^\delta}_{L^2}\le \delta} J_l(x).
               \end{equation*}
               In particular, this implies that $(x_n-f^\delta)_{n\in \N}$ is bounded in $L^1(\Omega)$. Additionally, there exist $u_n$ and $v_n$ such that $u_n+v_n = x_n$ and $J_l(x_n) = g_l(u_n) +h(v_n)$ for all $n \in \N$. Thus, by using the same arguments as before, we can assume that $(u_n)_{n\in \N}$ and $(v_n)_{n\in \N}$ are weakly convergent (up to a subsequence) in $L^\frac{d}{d-1}(\Omega)$. Again, the claim follows by lower-semicontinuity of $g_l$ and $h$ and Lemma \ref{lemma:coercivity_bregmanized}.
           \end{enumerate}
       \end{proof}
\end{lemma}

\subsection{\texorpdfstring{$TV-H^1$}{TEXT}}\label{subsec:TV_H1}
We formalize the setting of Example \ref{ex:H1-TV} in a more general context. For this consider a bounded Lipschitz domain  $\Omega \subset \R^d$. Let $X = L^2(\Omega)$ and consider $A \in \L(X,Y)$, where $Y$ is some Hilbert space. In order to separate smooth and piecewise constant components of a function $x \in L^2(\Omega)$, we use the  squared $H^1$-seminorm as a regularizer for the smooth part:  \begin{equation*}
     g(u) = \frac{\alpha}{2}\norm{\nabla u}_{L^2}^2 = \int\limits_\Omega \norm{\nabla u}_2^2,
 \end{equation*}
for $u \in H^1(\Omega)$. For $u \notin H^1(\Omega)$, we set $g(u) = \infty$. The piecewise constant part will be penalized by the total variation:
\begin{equation*}
    h(v) = \beta\tv{v}.
    \end{equation*}
    Since  $H^1(\Omega) \subset BV(\Omega)$, we obtain $\dom (g\square h) = BV(\Omega)$  by Proposition \ref{prop:inf_conv_properties} (iii). \\

  \begin{remark}\label{remark_Huber}
 Let us point out the similarity of the regularizer $J =  g\square  h$  with the Huber-TV functional \cite{BilevelHuberTV,journal_tvlp,HintermuellerInfConv}. The Huber-TV functional is defined on $BV(\Omega)$ as a convex function of the measure $Du$ \cite{Dem}, i.e.
    \begin{equation*}
     TV_\gamma(x) :=  \int\limits_\Omega f_\gamma(\nabla x) + \norm{D^s x}_{\mathcal{M}},
 \end{equation*}
where $f_\gamma(y)$ is given for some $\gamma>0$ and  $y\in \mathbb{R}^{d}$ by
\begin{equation*}
    f_\gamma (y) := \left(\frac{1}{2\gamma}\norm{\cdot}_2^2 \square \norm{\cdot}_2\right)(y) = \begin{cases}
        \frac{\norm{y}_2^2}{2\gamma}  &\text{ if } \norm{y}_2 \le \gamma,\\
        \norm{y}_2 - \frac{\gamma}{2} &\text{ if } \norm{y}_2\ge \gamma.
    \end{cases}
\end{equation*}
Here $\nabla x$ and $D^s x$ denote the absolutely continuous and  the singular part of $Dx$ with regard to the Lebesgue measure, respectively. In particular, there holds $\nabla x \in L^1(\Omega)$. Setting $\gamma = \frac{\beta}{\alpha}$, we have that \begin{align}\label{eq:Huber_TV_integral}
     \beta TV_\gamma(x) &= \int\limits_\Omega \pars{\inf_{w(z)\in \R^d}\frac{\alpha}{2}\norm{ w(z)}_2^2 +\beta\norm{\nabla x(z) -w(z) }_2}dz + \beta \norm{D^s x}_{\mathcal{M}}\notag\\
     &= \inf\limits_{w \in L^2(\Omega,\R^d)} \int\limits_\Omega \pars{\frac{\alpha}{2}\norm{ w(z)}_2^2 +\beta\norm{\nabla x(z) -w(z) }_2}dz + \beta\norm{D^s x}_{\mathcal{M}}.
 \end{align}
 Note that the last minimization problem in \eqref{eq:Huber_TV_integral} admits indeed a solution $w_{\gamma}\in L^{\infty}(\Omega)$, see \cite{journal_tvlp} and also \eqref{w_gamma} below.
Next, since $\tv{x} = \norm{\nabla x}_{L^1} + \norm{D^sx}_{\mathcal{M}}$, for all $x \in BV(\Omega)$ we have 
\begin{align}\label{eq:Tv-H1_integral}
    (\alpha g\square \beta h)(x) &= \inf\limits_{u \in H^1(\Omega)} \pars{ \int\limits_\Omega \pars{\frac{\alpha}{2}\norm{\nabla u(z)}_2^2+ \beta\norm{\nabla(x(z)- u(z))}_2}dz +\beta\norm{D^s (x-u)}_{\mathcal{M}}}\notag \\
 &\overset{D^su = 0}{=} \inf\limits_{u \in H^1(\Omega)} \pars{\int\limits_\Omega \pars{\frac{\alpha}{2}\norm{\nabla u(z)}_2^2+ \beta\norm{\nabla x(z)- \nabla u(z)}_2}dz +\beta\norm{D^s x}_{\mathcal{M}}}\notag\\
    & = \inf\limits_{ \substack{w \in \text{grad } H^1}} \pars{\int\limits_\Omega \pars{\frac{\alpha}{2}\norm{ w(z)}_2^2 +\beta\norm{\nabla x (z)- w(z) }_2}dz + \beta\norm{D^s x}_{\mathcal{M}}}. 
\end{align}
where we have defined $\text{grad } H^1 = \set{w \in L^2(\Omega,\R^d) : w = \nabla u\text{ for some } u \in H^{1}(\Omega)}$.
Comparing \eqref{eq:Huber_TV_integral} and \eqref{eq:Tv-H1_integral}, one obtains $\beta TV_\gamma(x) \le (\alpha g \square \beta h)(x)$ for every $x\in BV(\Omega)$. Define   \begin{equation*}
    H_0(\div 0) = \set{w \in L^2(\Omega,\R^d) : \div w = 0 \text{ and } w\cdot n = 0\text{ on } \partial \Omega},
\end{equation*} where $n$ denotes the outer unit normal of $\Omega$ and all derivatives have to be understood in a weak sense. Then  $\text{grad } H^1$ and $H_0(\div 0)$ are both closed and orthogonal subspaces of $L^2(\Omega,\R^d)$, so that we obtain (see for instance \cite[p. 216]{BookSpectral}) \begin{equation*}
    L^2(\Omega,\R^d) = \text{grad } H^1 \oplus H_0(\div 0).
\end{equation*}   
This means that in general, we do not expect $(\alpha g \square \beta h) (x) = \beta TV_\gamma (x)$. However, for $d = 1$, one can easily see that $H_0(\div 0) = \set {0}$, so that equality holds. Furthermore, it is 
\begin{equation}\label{w_gamma}
    \argmin\limits_{w \in \R^d}\;\set{\frac{\alpha}{2}\norm{w}_2^2 +\beta\norm{\nabla x(z) -w}_2} = \min\set{1,\frac{\beta}{\alpha\norm{\nabla x(z)}_2}}\nabla x(z):= w_\gamma(z),
\end{equation}
for all $z \in \Omega$.
Thus, equality of $(\alpha g \square \beta h) (x)$ and $\beta TV_\gamma(x)$ also holds if $w_\gamma \in \text{grad } H^1$. For instance, this happens  for all $x$ such that $\norm{\nabla x(z)}_2$ is constant or $\norm{\nabla x(z)}_{2} \le \frac{\alpha}{\beta}$ almost everywhere.
 \end{remark}

Next, we give a formal proof of existence for the minimization problems in Algorithm \ref{algo:NestedBregman_Noisy_Morozov} with those regularizers. As a first step of our analysis, we establish a Poincaré-Wirtinger inequality for the infimal convolution of $g$ and $h$. To this end, let $p = \frac{d}{d-1}$ for $d \ge 2$ for the remainder of the section. In the case $d = 1$, we set $p = 2$ for simplicity. However, in the one-dimensional case all results remain valid for any $1\le p\le \infty$.
\begin{corollary}\label{cor:Poincare_Wirtinger_TV-H1}
    Let $J = g\square h$,  and define $\tilde J (\tilde{x}) = \min\limits_{u\in L^p(\Omega)} \frac{\alpha}{2}\norm{u}_{L^p}^2 + \beta\norm{\tilde{x}-u}_{L^p}$ for every $\tilde{x}\in L^{p}(\Omega)$. Furthermore, denote $\mu_u = \frac{1}{\abs{\Omega}}\int\limits_{\Omega} u$ the mean of a function $u\in L^{p}(\Omega)$. Then there is a constant $C>0$ (which only depends on $\Omega$) such that \begin{equation}\label{eq:Poincare_Wirtinger_TV-H1}
        \tilde J(x-\mu_x) \le C J(x)
    \end{equation}
    for all $x \in L^{p}(\Omega)$.
    \begin{proof}
        Since $\dom J=BV(\Omega)$, we can assume that  $x \in BV(\Omega)$. Let $u \in L^{p}(\Omega)$, then by the Poincaré-Wirtinger inequalities for $BV(\Omega)$ and $H^1(\Omega)$ and the boundedness of $\Omega$, there is $C >0$, which only depends on $\Omega$, such that \begin{equation*}
            \tilde J(x-\mu_x)   \le  \frac{\alpha}{2} \norm{u-\mu_u}_{L^p}^2 + \beta \norm{(x-\mu_x)-(u-\mu_u)}_{L^p} \le C \pars{\frac{\alpha}{2} \norm{\nabla u}_{L^2}^2 + \beta \tv{(x-u)}}.
        \end{equation*}
        Since $u$ was chosen arbitrarily, we can take the infimum over all $u$ on the right-hand side of the previous equation to obtain \eqref{eq:Poincare_Wirtinger_TV-H1}.
    \end{proof}
\end{corollary}

This allows to prove the existence result.

\begin{lemma}\label{lemma:existence_TV_H1}
   Let $g(u) = \frac{\alpha}{2}\norm{\nabla u}_{L^2}^2 $, $h(v) = \beta \tv{v}$ and $A \in \mathcal{L}(L^2(\Omega),Y)$. For $d \ge 3$ assume additionally that $A$ is continuous on $L^p(\Omega)$. Furthermore, assume that there is $x \in \dom g\square h$ with $Ax = f$. Then for any $l \in \N$, the infimal convolutions  in \eqref{eq:NestedBregman_Morozov_initial} and \eqref{eq:NestedBregman_Morozov_step} are exact, and the minimization problems are well-defined.
\begin{proof}

Let $l \in \N$. For simpler notation, denote again $g_l(\cdot):=D_g^{p_{l-1}^\delta}(\cdot,u_{l-1}^\delta) $ with the convention $g_l = g$ if $l = 1$ throughout the proof. Additionally define $J_l = g_l\square h$. Note that, as in the proof of Lemma \ref{lemma:existence_L1_H1}, the functionals $g_l$ and $h$ are lower semicontinuous in $L^p$. Furthermore, for any constant function $c(x) = c\in \R$,  it is $\inner{p_{l-1}^\delta}{c} = \inner{p_{l-1}^\delta}{u_{l-1}^\delta+c-u_{l-1}^\delta} \le \frac{\alpha}{2}\norm{\nabla (u_{l-1}^\delta+c)}_{L^2}^2 - \frac{\alpha}{2}\norm{\nabla u_{l-1}^\delta}_{L^2}^2 = 0$. By considering $-c$ instead, we conclude $\inner{p_{l-1}^\delta}{c} = 0$. Hence, $g_l(u+c) = g_l(u)$ for all $u\in \dom g$ and $c \in \R$. 
\begin{enumerate}
\item Let us first show the exactness of the infimal convolution. Let $x \in \dom J_l$, and note that clearly $J_l(x) \ge 0$. Therefore, let $\seq{u}$ be a minimizing sequence, that is $u_n\in\dom J_l$ for all $n \in \N$ and $\lim\limits_{n \to \infty} g_l(x-u_n) + h(u_n) = J_l(x)$. Since for any $c \in \R$, one has $g_l(x-(u_n+c))+h(u_n+c) = g_l(x-u_n) + h(u_n)$, we may assume that $\mu_{u_n} = 0$.  Since $x \in \dom J_l$, we must have that $(h(u_n))_{n \in \N}$ is bounded. Hence, the Poincaré-Wirtinger inequality for $TV$ implies that $(u_n)_{n \in \N}$ is bounded in $L^p(\Omega)$. Consequently, (by passing to a subsequence), we can assume that $(u_n)_{n \in \N}$ converges weakly  to some $u^*\in L^p(\Omega)$. Since $g_l$ and $h$ are lower semicontinuous with respect to $L^p(\Omega)$, we obtain \begin{equation*}
       g_l(x-u^*) +h(u^*) \le \liminf_{n \to \infty} \pars{g_l(x-u_n) + h(u_n)}   = J_l(x). 
    \end{equation*}
    Hence, $J_{l}$ is exact.

    \item Next, we show that the constraint problem  \begin{equation*}
       \begin{cases}
           &\min J_l(x) \\
           & s.t. \  \norm{Ax-f^\delta}\le \delta
       \end{cases} 
    \end{equation*} is well-defined. Let $\seq{x}$ be a sequence with $\norm{Ax_n - f^\delta}\le \delta$ for all $n\in \N$ and $\lim\limits_{n \to \infty} J_l(x_n) = \inf\limits_{\norm{Ax-f^{\delta}}\le \delta}J_l(x)$. We first show that $(x_n)_{n \in \N}$ can be assumed to have bounded mean. We claim that $(x_n-\mu_{x_n})_{n \in \N}$ is bounded in $L^p(\Omega)$. By contradiction, suppose that (on a subsequence denoted the same) $\lim\limits_{n \to \infty} \norm{x_n-\mu_{x_n}}_{L^p} = \infty$. Let $\tilde J$ as in Corollary \ref{cor:Poincare_Wirtinger_TV-H1}. One can easily see that $\tilde J$ is coercive with respect to $L^p(\Omega)$. Therefore $\lim\limits_{n \to \infty}\tilde J(x_n-\mu_{x_n}) = \infty$. This implies $\lim\limits_{n\to \infty}(g\square h) (x_n-\mu_{x_n}) = \infty$. But that would mean $\lim\limits_{n \to \infty}J_l(x_n-\mu_{x_n}) = \lim\limits_{n \to \infty}J_l(x_n) = \infty$ by Lemma \ref{lemma:coercivity_bregmanized}, which is a contradiction. Thus,  $(x_n-\mu_{x_n})_{n \in \N}$ and consequently $(A(x_n-\mu_{x_n}))_{n \in \N}$ must be bounded in $L^p(\Omega)$ by continuity of $A$. Consider two cases. First assume $A\chi_\Omega \neq 0$, where $\chi$ is the characteristic function of $\Omega$. If $(\mu_{x_n})_{n \in \N}$ was unbounded, then $(Ax_n)_{n \in \N}$ would have to be unbounded as well, which contradicts $\norm{Ax_n-f^\delta} \le \delta$ for all $n \in \N$. Otherwise, if $A\chi_\Omega = 0$, we can replace $(x_n)_{n\in \N}$ by a sequence with mean $0$, since $A(x_n+c) = Ax_n$ and $J_l(x_n+c) = J_l(x_n)$ for all $c \in \R$. Combining both cases, we can choose $(x_n)_{n\in \N}$ such that $(\mu_{x_n})_{n\in \N}$ is bounded. \\
    Next, since $J_l$ is exact, for any $n \in \N$ there exist $u_n\in \dom g_l,\ v_n\in \dom h$ with $u_n+v_n = x_n$ and $J_l(x_n)= g_l(u_n) +h(v_n) $. Since adding a constant to $u_n$ and subtracting it from $v_n$ does not change the value of $g_l(u_n) +h(v_n)$ (see part (i) of this proof), we may assume that $\mu_{v_n} = 0$ for all $n\in \N$. This means $(v_n)_{n\in \N}$ is bounded by non-negativity of $g$ and the Poincaré-Wirtinger inequality of $h$:
    \begin{equation*}
       \beta \norm{v_n}_{L^p} \le C\beta \tv{v_n}  \le CJ_l(x_n).
    \end{equation*}
 In summary, we have that both $(x_n)_{n\in \N}$ and $(v_n)_{n\in \N}$ must be bounded in $L^p(\Omega)$. Consequently, $(u_n)_{n\in \N}$ must also be bounded in $L^p(\Omega)$. Therefore, (by passing to subsequences) we can assume that $(u_n)_{n\in \N}$ and $(v_n)_{n\in \N}$ converge weakly in $L^p(\Omega)$ to some $u^*,\, v^* \in \L^p(\Omega)$. By continuity and linearity of $A$, the set $\set{x \in X :\norm{Ax-f^\delta}\le \delta}$ is $L^p-$weakly closed, therefore one has $\norm{A(u^*+v^*)-f^\delta}\le \delta$. Because $g_l$ and $h$ are lower-semicontinuous on $L^p(\Omega)$, it holds that \begin{equation*}
        g_l(u^*) +h(u^*) \le \liminf\limits_{n \to \infty} \pars{g_l(u_n) +h(v_n)} = \lim\limits_{n \to \infty}J_l(x_n) = \inf\limits_{\norm{Ax-f^\delta}\le \delta} J_l(x).
    \end{equation*} 
    This shows the existence of a minimizer.  
     \end{enumerate}
\end{proof} 
\end{lemma}

\subsection{Oscillatory TGV}
As a last example, we recall the infimally-convoluted oscillation TGV from \cite{oscICTV}. To be precise, let $\Omega\subset \R^d$ be a bounded domain with Lipschitz boundary and $Y$ be some Hilbert space. Consider the functionals \begin{equation}\label{eq:TGV}
    g(u) = TGV_{\alpha_1,\beta_1}^2(u)  = \min\limits_{{w}\in BD(\Omega)} \alpha_1 \norm{\nabla u-w}_{\mathcal{M}} + \beta_1\norm{\mathcal{E}w}_{\mathcal{M}},
\end{equation} 
and \begin{equation}\label{eq:TGV_osci}
    h(v) = TGV_{\alpha_2,\beta_2,C}^{osci}(v) = \min\limits_{w \in BD(\Omega)} \alpha_2\norm{\nabla v-w}_{\mathcal{M}} + \beta_2\norm{\mathcal{E}w + C{v}}_{\mathcal{M}},
\end{equation}
where $\mathcal{E}$ denotes the weak symmetrized derivative, $BD(\Omega)$ is the space of functions of bounded deformation, and $C\in \R^{d\times d}$ is defined as the matrix with entries $c_{ij} = \omega_i\omega_j$ for some $\omega \in \R^d$. That means using $g$ and $h$ as regularizers promotes piecewise affine and  oscillatory functions, respectively (compare to Lemma \ref{lem:properties_TGV_osci} (ii) below). Additionally, note that the total generalized variation $g$ is defined via an infimal convolution itself. Thus, this example serves to show that our proposed method can also be applied for problems involving more sophisticated regularizers that yield decompositions with components having quite different structures.\\
In order to show exactness of the infimal convolutions and  well-definedness of the minimization problems occurring in Algorithm \ref{algo:NestedBregman_Noisy_Morozov}, we use the arguments from \cite[Section 3]{oscICTV} and apply Lemma \ref{lemma:coercivity_bregmanized} where necessary. Let us first state some properties of $g$ and $h$. Again, we denote $p = \frac{d}{d-1}$ for $d\ge 2$ and $p = 2$ for $d = 1$.
\begin{lemma}\label{lem:properties_TGV_osci}
    Let $\alpha_1,\alpha_2,\beta_1,\beta_2\ge 0$ and $C = \omega^T\omega$ with $\omega \in \R^d\setminus\set{0}$. Define $g(u) = TGV_{\alpha_1,\beta_1}^2(u)$ and $ h(v) = TGV_{\alpha_2,\beta_2,C}^{osci}(v)$. Then the  following statements hold. \begin{enumerate}
        \item $g$ and $h$ are seminorms and lower semicontinuous with respect to $L^p$.
        \item It is \begin{equation*}
            \ker g = \set{u(z) = a\cdot z +b :  a\in \R^d, b\in \R}
        \end{equation*} and \begin{equation*}
            \ker h = \set{v(z) = a \sin(\omega\cdot z) + b\cos(\omega\cdot z): a,b \in \R}.
        \end{equation*} Furthermore, $g\square h$ is a lower semicontinuous seminorm on $L^p(\Omega)$ with kernel $\ker g\square h = \ker g + \ker h$.
        \item Let $P_g$ and $P_h$ be the $L^p-$projection on $\ker g$ and $\ker h$, respectively. Then there are constants $c_1,c_2 >0$ such that \begin{equation}\label{eq:Poincare_Wirtinger_TGV}
            \norm{u-P_gu}_{L^p} \le c_1 g(u)
        \end{equation} for all $u \in L^2(\Omega)$ and \begin{equation}\label{eq:Poincare_Wirtinger_TGV_osci}
            \norm{v-P_hv}_{L^p} \le c_2h(v)
        \end{equation} for all $v \in L^p(\Omega)$. Furthermore, there exists $K >0$ such that \begin{equation}\label{eq:seminorm_estimate_projections}
            \norm{x}_{L^p} \le K\pars{\norm{x-P_gx}_{L^p}+\norm{x-P_hx}_{L^p}}
        \end{equation} for all $x\in L^p(\Omega)$
    \end{enumerate}
    \begin{proof}
        All statements follow from Proposition $2$, Proposition $4$, Lemma $3$ and Theorem 1 in \cite{oscICTV} by taking into account that $\ker g \cap \ker h = \set{0}$. 
    \end{proof}
\end{lemma}

\begin{lemma}\label{lemma:existence_TGV_osci}
    Let $g(u) = TGV_{\alpha_1,\beta_1}^2(u)$, $ h(v) = TGV_{\alpha_2,\beta_2,C}^{osci}(v)$ and $A \in \L(L^2(\Omega,Y))$. Additionally, assume that there is $x \in \dom g\square h$ with $Ax =f$.Then for any $l \in \N$, the infimal convolutions  in \eqref{eq:NestedBregman_Morozov_initial} and \eqref{eq:NestedBregman_Morozov_step} are exact, and the minimization problems are well-defined.
    \begin{proof}
        Let $l \in N$. As before, denote $g_l(\cdot):=D_g^{p_{l-1}^\delta}(\cdot,u_{l-1}^\delta) $ and $J_l = g_l\square h$  with the convention $g_l = g$ if $l = 1$ throughout the proof.
        \begin{enumerate}
            \item For showing the exactness of the infimal convolution, let $x \in \dom g\square h$ and consider a sequence $(u_n,v_n)_{n \in \N}$ such that $u_n+v_n = x$ for all $n\in \N$ and $\lim\limits_{n \to \infty} g_l(u_n) +h(v_n) = J_l(x)$. By non-negativity of $g_l$ and $h$ and convergence, this implies that $(g_l(u_n))_{n\in \N}$ and $(h(v_n))_{n\in \N}$ must be bounded. By \eqref{eq:Poincare_Wirtinger_TGV_osci}, this means $(v_n -P_hv_n)_{n\in \N}$ is bounded in $L^p\Omega)$. Using Lemma \ref{lemma:coercivity_bregmanized} with $f$ replaced by $x$ and arbitrary $\delta >0$ (see also Remark \ref{rem:coercivity_for_other_f}), we also obtain that $(u_n-P_g u_n)_{n\in \N}$ is bounded in $L^p(\Omega)$ by \eqref{eq:Poincare_Wirtinger_TGV}. Therefore, the sequence $v_n -P_g v_n = x-u_n -P_g(x-u_n) = (x-P_gx)- (u_n-P_gu_n)$ must be bounded in $L^p(\Omega)$, so that by \ref{eq:seminorm_estimate_projections} $(v_n)_{n\in \N}$ must be bounded in the same space. We can thus pass to a subsequence (still denoted by $(v_n)_{n\in \N}$) which converges weakly in $L^p(\Omega)$ to some $v^*\in L^p(\Omega)$. Consequently, $(u_n)_{n\in \N}$ converges weakly in $L^p(\Omega)$ to $u^* = x-v^*$ (up to a subsequence). Using lower semicontinuity of $g_l$ and $h$, it must be $J_l(x) = g_l(u^*) +h(v^*)$.
            \item Note that, as in the proof of Lemma \ref{lemma:existence_TV_H1}, it is $\inner{p_{l-1}^\delta}{x} = 0$ for all $x \in \ker g$. Therefore, let $(u_n,v_n)_{n\in \N}$ be a minimizing sequence, i.e.  $\norm{A(u_n+v_n) -f^\delta}\le \delta$ and $\lim\limits_{n\to \infty} g_l(u_n) + h(u_n) = \inf\limits_{\norm{A(x)-f^\delta}\le \delta}J_l(x)$. Let $P_{g,A}$ be the $L^p$-projection onto $\ker g\cap \ker A$ and $P_{h,A}$ the $L^p$-projection on $\ker h\cap \ker A$. Now replace $u_n$ by $\tilde u_n = u_n -P_{g,A}u_n$ and $v_n$ by $\tilde v_n - P_{h,A}$ for each $n\in \N$. In particular, this means $A(\tilde u_n + \tilde v_n) = A(u_n+v_n)$, $g_l(\tilde u_n) = g_l(\tilde u_n)$, $h(\tilde v_n) = h(v_n)$, $\tilde u_n - P_g\tilde u_n = \tilde u_n$ and $\tilde v_n - P_h \tilde v_n = \tilde v_n$. Thus, by Lemma \ref{lemma:coercivity_bregmanized}, \eqref{eq:Poincare_Wirtinger_TGV} and \eqref{eq:Poincare_Wirtinger_TGV_osci}, we must have that $(\tilde u_n)_{n\in \N}$ and $(\tilde v_n)_{n\in \N}$ are bounded in $L^p(\Omega)$. Again, passing to weakly convergent subsequences and using lower semicontinuity finishes the proof.
        \end{enumerate}
    \end{proof}
\end{lemma}

\begin{remark}
 Lemmas \ref{lemma:existence_L1_H1},\ref{lemma:existence_TV_H1} and \ref{lemma:existence_TGV_osci} also hold with the roles of $g$ and $h$ in Algorithm \ref{algo:NestedBregman_Noisy_Morozov} interchanged. 
\end{remark}

\section[Numerical results]{Numerical results\footnote{The program code is available as ancillary file from the arXiv page of this paper.}}\label{sec:numerical_results}

The purpose of this section is to illustrate the general behavior of our proposed method, as well as its ability to produce approximate solutions of \eqref{eq:ill_posed_problem} with meaningful components $u$ and $v$ when equipped with a suitable stopping criterion. To this end, we apply Algorithms \ref{algo:NestedBregman_Noisy_Morozov} and \ref{algo:NestedBregman_Noisy} with the pairs of  regularizers introduced in Section \ref{sec:selected_regularizers} to synthetic images. Furthermore, we use a stopping criterion that does not depend on the regularizers, but rather implicitly assumes that the optimal components are structurally dissimilar. For more sophisticated applications on real data, we would recommend to use problem specific criteria that take into consideration a priori information about the expected components.

\subsection{Stopping criterion}
Following \cite{TernaryImageDecomposition}, we recall the discrete normalized cross-correlation of two non-zero signals, respectively matrices. \begin{definition}
    \begin{enumerate}
        \item Given $u,v \in \R^N\setminus\set{0}$, the sample normalized cross-correlation of $u$ and $v$ is defined as the vector $\rho(u,v) \in \R^{2N-1}$ with components
        \begin{equation}\label{eq:normalized_cross_correlation_1D}
            \rho_{k}(u,v) = \frac{1}{\norm{u}_2\norm{v}_2}\sum\limits_{i = 1}^N u_i v_{i+k},
        \end{equation}  
        for $k \in \set{-(N-1),\dots ,0,\dots ,N-1}$. Here we assume a periodic extension of $v$, that is $v_{i\pm N} := v_i$ for all $i \in \set{1,\dots,N}$.
        \item For $u,v \in \R^{N\times M}\setminus\set{0}$, we analogously define the sample normalized cross-correlation as the matrix $\rho(u,v) \in \R^{(2N-1)\times(2M-1)}$ with entries 
        \begin{equation}\label{eq:normalized_cross_correlation_2D}
            \rho_{k,l}(u,v) = \frac{1}{\norm{u}_2 \norm{v}_2} \sum\limits_{i = 1}^N \sum_{j = 1}^M u_{i,j} v_{i+k,j+l}
        \end{equation}
        for $(k,l) \in \set{-(N-1),\dots ,0,\dots ,N-1}\times \set{-(M-1),\dots ,0,\dots ,M-1}$, where $\norm{\cdot}_2$ denotes the Frobenius norm of a matrix. Again, we assume periodic boundary conditions: $v_{i\pm N, j\pm M} := v_{i,j}$.
    \end{enumerate}
\end{definition}

\begin{remark}
    Under the assumption of periodic boundary conditions, we can define $u'_j = u_{-j}$ for signals and $u'_{i,j} = u_{-i,-j}$ for images, so that \eqref{eq:normalized_cross_correlation_1D} and \eqref{eq:normalized_cross_correlation_2D} can be written as \begin{equation*}
        \rho(u,v)  = \frac{1}{\norm{u}_2\norm{v}_2} (u'*v),
    \end{equation*}
    where $*$ denotes the discrete convolution. Since  the components appearing as minimizers in Algorithms \ref{algo:NestedBregman_Noisy_Morozov} and \ref{algo:NestedBregman_Noisy} with the regularizers from  Section \ref{sec:selected_regularizers} are also in $L^2(\Omega)$ for dimension $d \in \set{1,2}$, this allows for a straightforward generalization of the normalized cross correlation in the function space setting. However, for working in a discretized setting, the formulations \eqref{eq:normalized_cross_correlation_1D} and \eqref{eq:normalized_cross_correlation_2D} are sufficient.
\end{remark}

The normalized cross correlation has been used to determine optimal parameter choices in Tikhonov regularization for decompositions \cite{AujaolGilboaChanOsher2006,TernaryImageDecomposition, QuaternaryImageDecomposition}, relying on the  idea that structurally different components should be uncorrelated or have small correlation. In \cite{AujaolGilboaChanOsher2006}, only the empirical correlation of two images (which corresponds to the entry $\rho_{0,0}$ in the notation \eqref{eq:normalized_cross_correlation_2D}) is considered for an a posteriori parameter choice principle. By using  increasing values of a regularization parameter $\lambda$, an observed image is decomposed in a cartoon part $u_\lambda$ and a texture part $v_\lambda$. Then the optimal parameter is determined as the first local minimum of $\rho_{0,0}(u_\lambda,v_\lambda)$ with respect to $\lambda$. The authors point out that the criterion is  expected to work well only for relatively simple decompositions. The paper \cite{TernaryImageDecomposition} employs a two-stage variational decomposition model. For the second step of the method, the authors introduce a more informative criterion by defining a scalar measure of correlation as \begin{equation}\label{eq:scalar_correlation}
   \mathcal{C}(u,v) := \frac{1}{MN} \norm{\rho(u,v)}_2^2
\end{equation}
 and by minimizing $\mathcal{C}$ in a bilevel approach, where the lower level problem is a minimization problem that decomposes a denoised image into a cartoon and a texture part. It turns out that combining these approaches allows us to define a stopping criterion for our experiments. Thus,  we stop the iteration at the first local minimum of $\mathcal{C}(u_n^\delta,v_n^\delta)$, taking care that the algorithm does not terminate after the first step. This means, even if $\mathcal{C}(u_1^\delta,v_1^\delta) \le \mathcal{C}(u_2^\delta,v_2^\delta)$, we continue iterating  until a local minimum is reached. We would like to stress that this criterion is heuristic, and the existence of a local minimum of $\mathcal{C}$ is not guaranteed by our theory.

\subsection{\texorpdfstring{Experiment 1: $L^1-H^1$-decompositions}{TEXT}}
We consider the denoising problem $(A= Id)$ for the $1-$dimensional signal, consisting of two components. The first component $u^\dagger$ is  a discretization of one period of a sine curve on $n = 300$ equidistant nodes. The second component $v^\dagger$ consists of multiple peaks. The noisy observation was created by adding Gaussian noise $\eta^\delta$ with mean $\mu = 0$ and variance $\sigma = 0.05$ to the sum of the components: $f^\delta  = u^\dagger+v^\dagger +\eta^\delta$ - see the first line in Figure \ref{fig:nested_morozov_l1_l2}. We  employ Algorithm \ref{algo:NestedBregman_Noisy_Morozov} with $g = \frac{\alpha}{2}\norm{D\cdot }_2$, $\alpha = 400$ and $h = \norm{\cdot}_1$, where $D$ is the following matrix for finite differences 
\begin{equation*}
    D = \begin{pmatrix*}
        &1&-1& &0\\
        &&\ddots&\ddots&&\\
        &0&&1&-1
    \end{pmatrix*}.
\end{equation*}

\begin{figure}[H]
    \centering
    \begin{minipage}[c]{0.49\textwidth}
        \includegraphics[scale = 0.5]{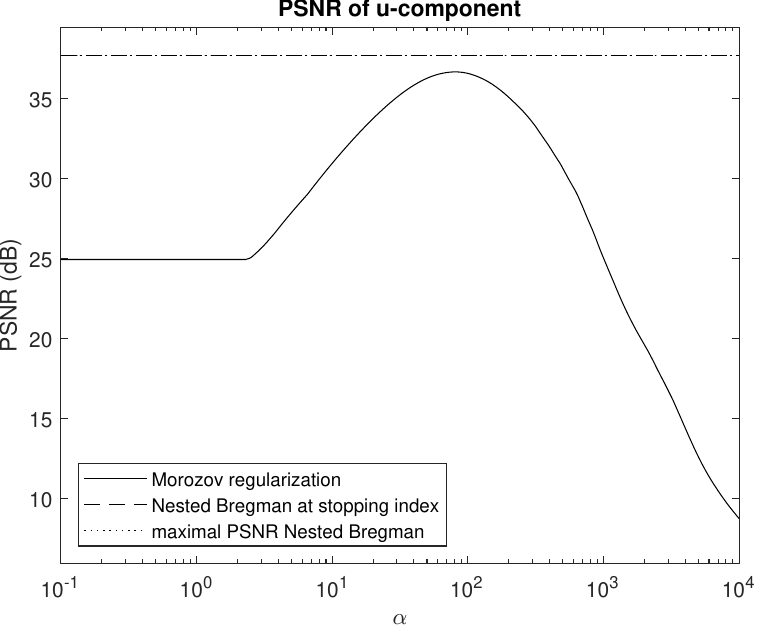}
    \end{minipage}
        \begin{minipage}[c]{0.49\textwidth}
        \includegraphics[scale = 0.5]{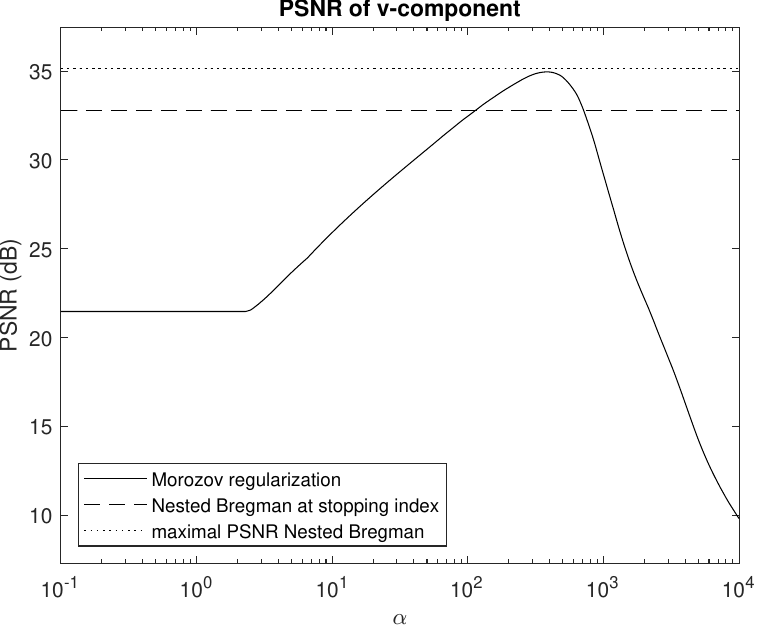}
    \end{minipage}
\vspace{0.7em}

        \begin{minipage}[c]{0.49\textwidth}
        \includegraphics[scale = 0.5]{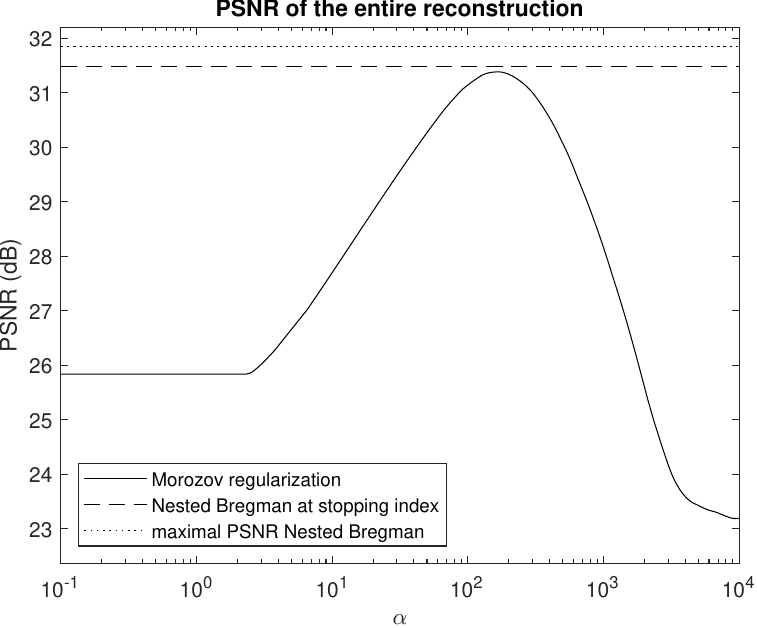}
    \end{minipage}
        \begin{minipage}[c]{0.49\textwidth}
        \includegraphics[scale = 0.5]{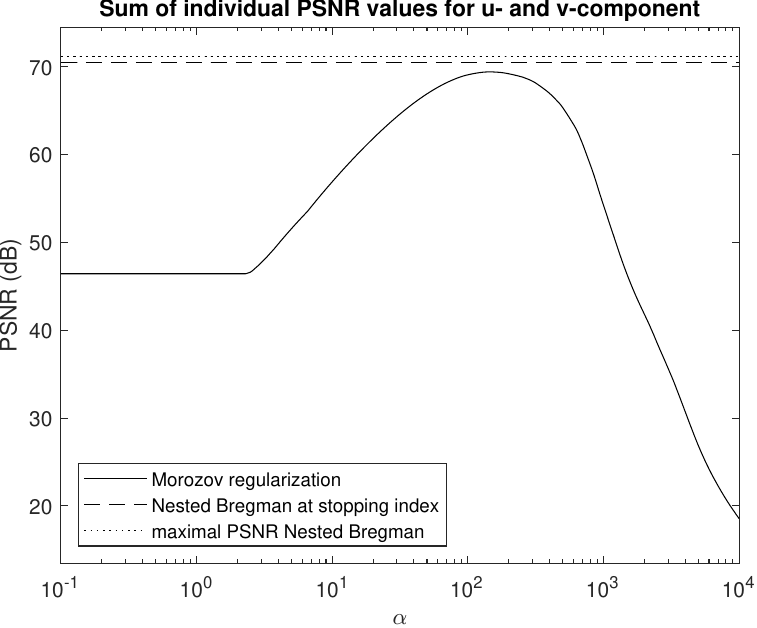}
    \end{minipage}
    \caption{{PSNR values for Morozov regularization for different values of $\alpha$, Nested Bregman iterations with weight $\alpha = 1000$ at the stopping index and maximal PSNR from Nested Bregman iterations. From top left to bottom right: $u-$component, $v-$component, full reconstruction, sum of $u-$component and $v-$component PSNR values.}}
    \label{fig:comp_Morozov_l1_l2}
\end{figure}

\begin{figure}[H]
    \centering
    \begin{minipage}[c]{0.49\textwidth}
    \includegraphics[scale = 0.5]{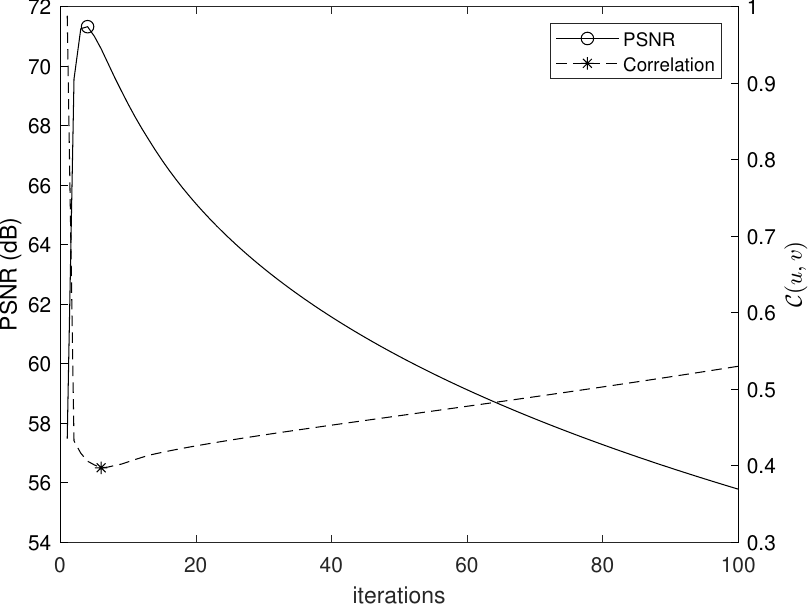}
\end{minipage}
\begin{minipage}[c]{0.49\textwidth}
    \includegraphics[scale = 0.5]{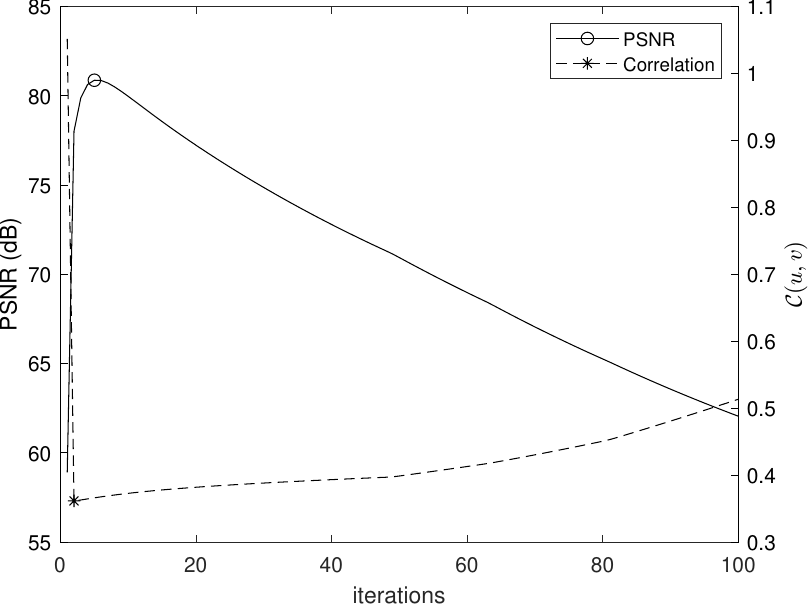}
\end{minipage}
\caption{Sum of the componentwise PSNR values obtained by Algorithm \ref{algo:NestedBregman_Noisy_Morozov} (left) and Algorithm \ref{algo:NestedBregman_Noisy} (right). The circle and the asterisk mark the maximum of the PSNR and the  first local minimum of the scalar normalized cross-correlation, respectively.}
\label{fig:PSNR_nested_l1_l2}
\end{figure}

Since the data are artificially created, we know the exact noise level $\delta = \norm{\eta^\delta}_2$. This means, in each step of the algorithm we solve \begin{equation}\label{eq:l1_l2_numrics_functional_morozov}
    \begin{cases}
        &\min \frac{\alpha}{2} \norm{Du}_2^2 + \norm{v}_1 - \inner{p_l^\delta}{u},\\
        &s.t \quad \norm{u+v{-f^\delta}}_2\le \delta.
    \end{cases}
\end{equation}
We employ the Matlab \textit{CVX} package for the numerical solution of this problem. Since $g$ is differentiable, we can update the subgradient in accordance with the second statement in Remark \ref{rem:update_differrentiable} by \begin{equation*}
    p_l^\delta = \alpha D^{*}Du_l^\delta.
\end{equation*} Figure \ref{fig:PSNR_nested_l1_l2} shows the sum of the PSNR values of the two components reconstructed using Algorithm \ref{algo:NestedBregman_Noisy_Morozov} (left) and Algorithm \ref{algo:NestedBregman_Noisy} (right) with $100$ outer iterations. The asterisk marks the iterate at which the first local minimum of $\mathcal{C}(u_l^\delta,v_l^\delta)$ defined as in \eqref{eq:scalar_correlation} is obtained. Note that, in both cases, this minimum is obtained closely to the maximum of the PSNR values, which is indicated by the circle in Figure \ref{fig:PSNR_nested_l1_l2}.  Furthermore, we observe that Algorithm \ref{algo:NestedBregman_Noisy} produces components with higher cumulated PSNR values. Figures \ref{fig:nested_morozov_l1_l2} and \ref{fig:nested_bregman_l1_l2} represent, from top to bottom, the true components and the noisy observation, the first iterates, the iterates with minimal scalar cross-correlation, the iterates with maximal cumulated PSNR values, and the  iterates after the final outer iteration obtained by Algorithm \ref{algo:NestedBregman_Noisy_Morozov} (Figure \ref{fig:nested_morozov_l1_l2}) and Algorithm \ref{algo:NestedBregman_Noisy} (Figure \ref{fig:nested_bregman_l1_l2}). In the initial step for both algorithms, the component $v_1^\delta$ identifies all peaks from the true component $v^\dagger$, but also includes some smaller peaks, that do not belong to  $v^\dagger$. This happens because $\alpha$ is chosen too small, thus resulting in $u_1^\delta$ being slightly over-regularized. When stopped at the first local minimum of the scalar cross-correlation  ($l = 6$), the iterates obtained by using Morozov regularization in the inner loop closely resemble the true components visually. Notably, the $v$ component with this stopping rule contains fewer small artifacts than the ones with the best PSNR values at $l = 4$. On the other side, the $u$ component of the maximal sum of PSNR values appears marginally smoother. When using Bregman iterations as a regularization method in the inner loop, the $v$ components with the best PSNR ($l = 5)$ and the minimal scalar cross-correlation ($l = 3$) are visually indistinguishable, while the $u$ component with the best correlation appears smoother. However, in general, one can not simply compare Algorithm \ref{algo:NestedBregman_Noisy_Morozov} and Algorithm \ref{algo:NestedBregman_Noisy}, since the iteration in Algorithm \ref{algo:NestedBregman_Noisy} depends on the initial weighting parameters in the infimal convolution, while the iteration in Algorithm \ref{algo:NestedBregman_Noisy_Morozov} only depends on the ratio of the weights.  Thus, for a fair comparison, it would be necessary to compute $\alpha$ and $\beta$ such that the first minimization problem in Algorithm \ref{algo:NestedBregman_Noisy} computes the same pair of minimizers as the first minimization in \ref{algo:NestedBregman_Noisy_Morozov}. For both algorithms, we also observe that the $u$ components become eventually overrepresented along the iterations, resulting in worse decompositions. This again illustrates the need for stopping the procedure according to some quality indicator of the decomposition. \\
{We test the performance of the proposed method in contrast to approaches using a single variational problem, such as bilevel methods with respect to some quality indicator or a grid search in parameter space. To this end, we compare the results from Algorithm \ref{algo:NestedBregman_Noisy_Morozov} to Morozov regularization for various parameters. Note that the minimizing pair in \eqref{eq:Morozov_decomposition} does not change when dividing the objective function by $\beta$. Thus fixing $\beta =1$ and varying $\alpha$ is sufficient to compute the reconstructions that are obtainable with Morozov regularization. For this we choose a logarithmically spaced grid of $1000$ values for $\alpha$ ranging from $10^{-1}$ to $10^4$, and perform the regularization with $J = \frac{\alpha}{2}\norm{D\cdot}_2 \square \norm{\cdot}_1$.  The results are illustrated in Figure \ref{fig:comp_Morozov_l1_l2} by solid lines for the PSNR values in the following cases: $u-$component, $v-$component, full reconstruction, as well as the sum of the PSNR values of the $u-$ and $v-$ component. The dashed and dotted lines represent the PSNR of Algorithm \eqref{algo:NestedBregman_Noisy_Morozov} at the stopping index and the maximal PSNR, respectively. Notably, the results from Algorithm \eqref{algo:NestedBregman_Noisy_Morozov} with the cross-correlation-based stopping criterion yield better results in all cases apart from the $v-$component, when the PSNR value is relatively large nonetheless. Moreover, the best PSNR values obtained with the Nested Bregman iterations are larger in all cases. We also observe that the proposed method behaves stably with respect to the weighting choice for the infimally convoluted functionals. As long as the weight in the Nested Bregman iterations is not such that the $v-$component is under-regularized in the first iteration of the method, repeating the experiment with different choices yields similar results. In summary, the proposed method is stable with respect to the weighting parameter and produces reconstructions that are more accurate than those that a grid search with an optimally chosen grid or any bilevel method could obtain since both approaches can only produce solutions of Morozov regularization. Furthermore, this indicates that similar to the improved contrast obtained by Bregman iterations compared to Tikhonov regularization with total variation regularizer (cf.\cite[Section 6]{Benning-Burger}), Nested Bregman iterations mitigate biases introduced by the choices of $g$ and $h$ in this case.}

 {Additionally}, we were able to numerically confirm the bound for the values $h(v_l^\delta)$. This is exhibited in Figure \ref{fig:nested_morozov_l1_l2_bound}, where the  {dashed and dotted lines show} the values of $h(v_l^\delta)$ {in Algorithms \ref{algo:NestedBregman_Noisy_Morozov} and \ref{algo:NestedBregman_Noisy} respectively}, while the {full} line marks the theoretical upper-bound $\frac{\alpha\norm{D (u^\dagger+v^\dagger)}_2^2}{l}$ from Theorem \ref{thm:Convergence_Nested_Bregman-Morozov}. Clearly, the value $h(u_l^\delta)$ lies significantly below the theoretical bound, as observed also in the other numerical experiments.

{
\begin{figure}[H] 
\centering
\hspace{-1cm}\begin{minipage}[c]{0.32\textwidth}
    \scalebox{.4}{\includegraphics{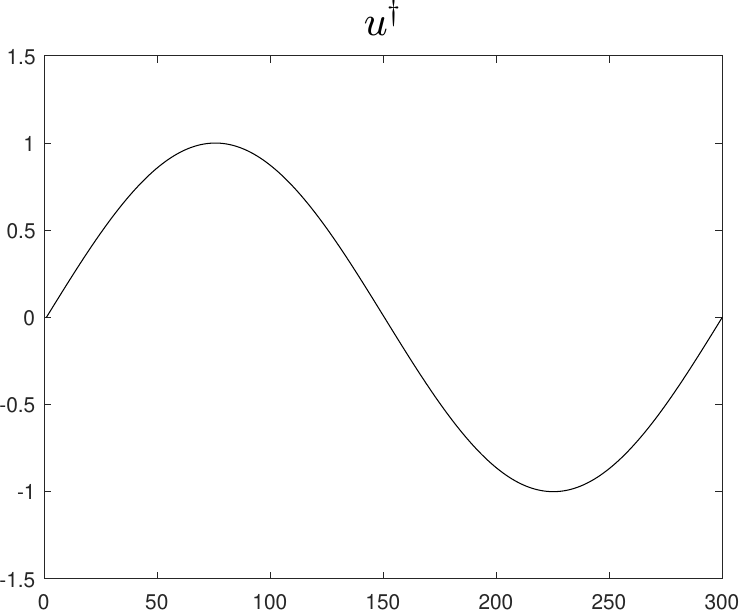}}
\end{minipage}
\begin{minipage}[c]{0.32\textwidth}
    \scalebox{.4}{\includegraphics{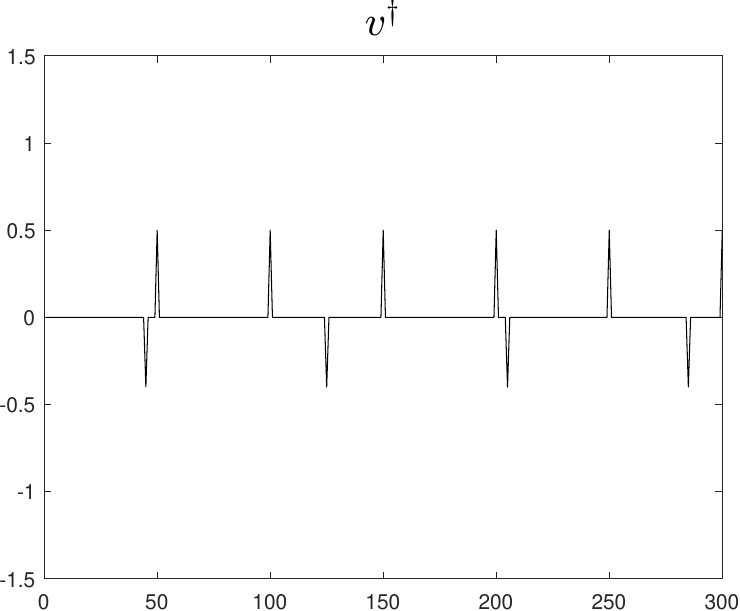}}
\end{minipage}
\begin{minipage}[c]{0.32\textwidth}
   \scalebox{.4}{\includegraphics{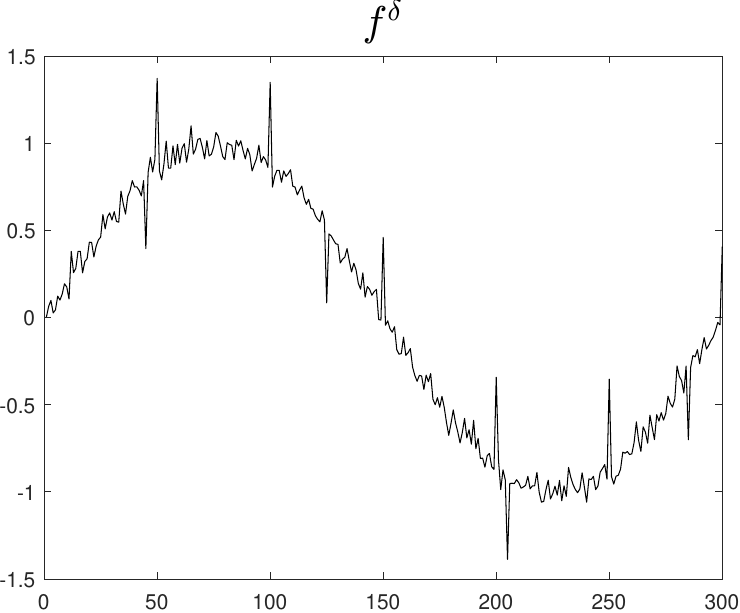}}
\end{minipage}

\vspace{0.5cm}
\centering
\hspace{-1cm}\begin{minipage}[c]{0.32\textwidth}
    \scalebox{.4}{\includegraphics{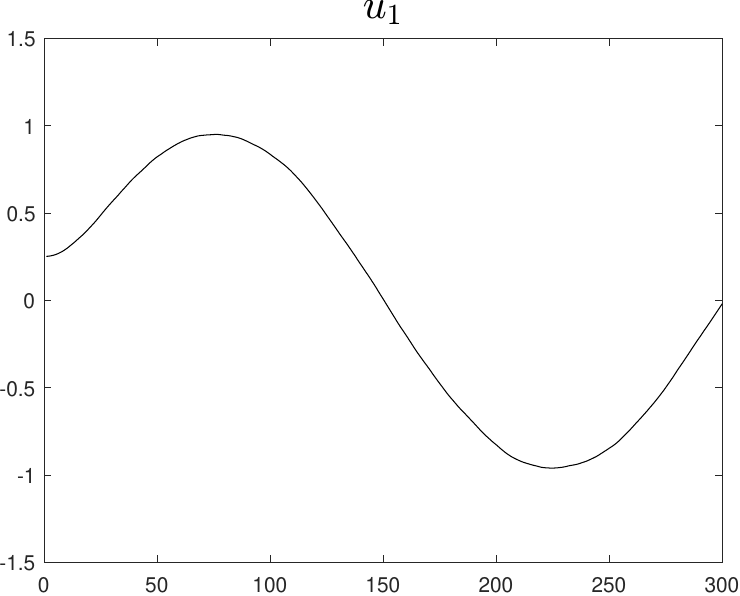}}
\end{minipage}
\begin{minipage}[c]{0.32\textwidth}
    \scalebox{.4}{\includegraphics{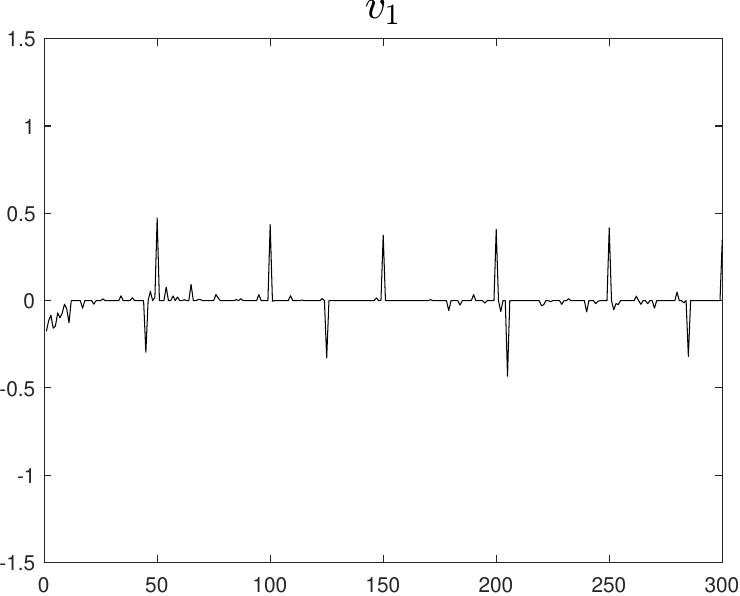}}
\end{minipage}
\begin{minipage}[c]{0.32\textwidth}
    \scalebox{.4}{\includegraphics{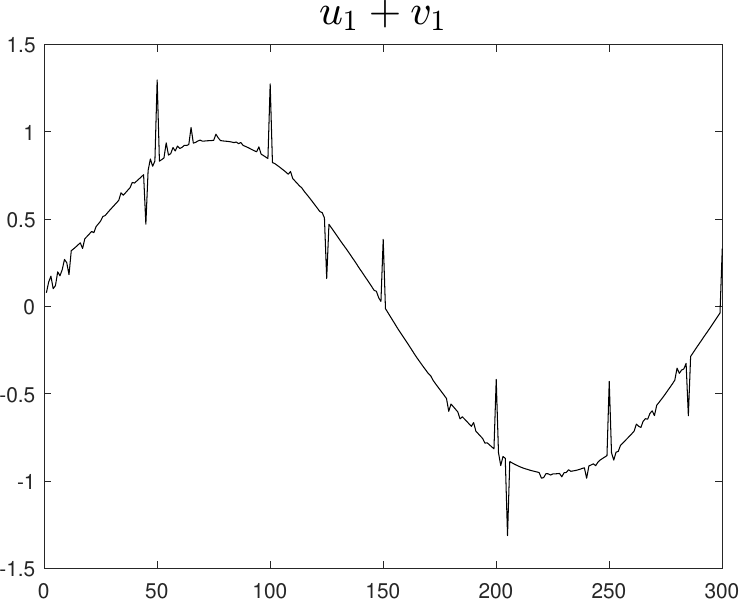}}
\end{minipage}

\vspace{0.5cm}
\centering
\hspace{-1cm}\begin{minipage}[c]{0.32\textwidth}
    \scalebox{.4}{\includegraphics{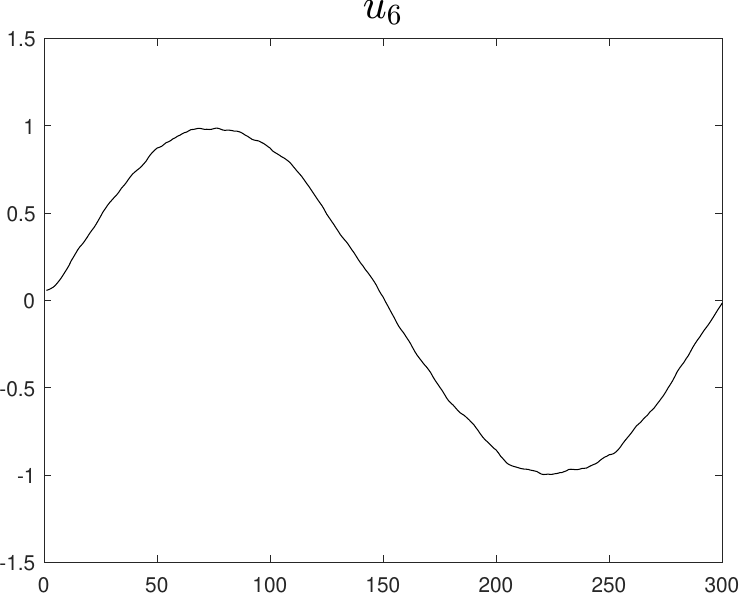}}
\end{minipage}
\begin{minipage}[c]{0.32\textwidth}
    \scalebox{.4}{\includegraphics{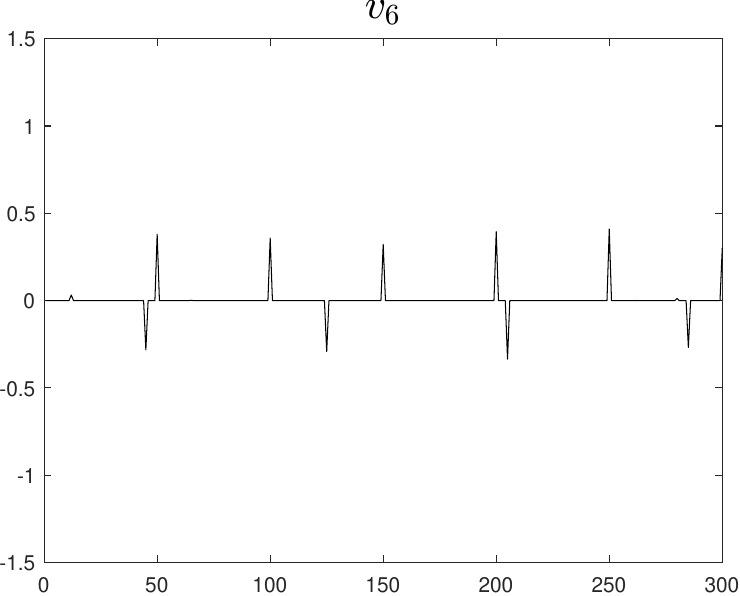}}
\end{minipage}
\begin{minipage}[c]{0.32\textwidth}
    \scalebox{.4}{\includegraphics{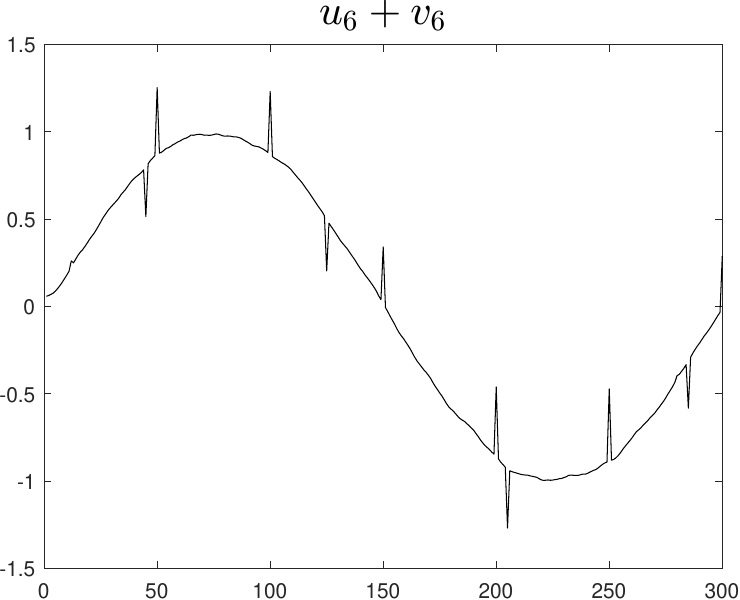}}
\end{minipage}

\vspace{0.5cm}
\centering
\hspace{-1cm}\begin{minipage}[c]{0.32\textwidth}
    \scalebox{.4}{\includegraphics{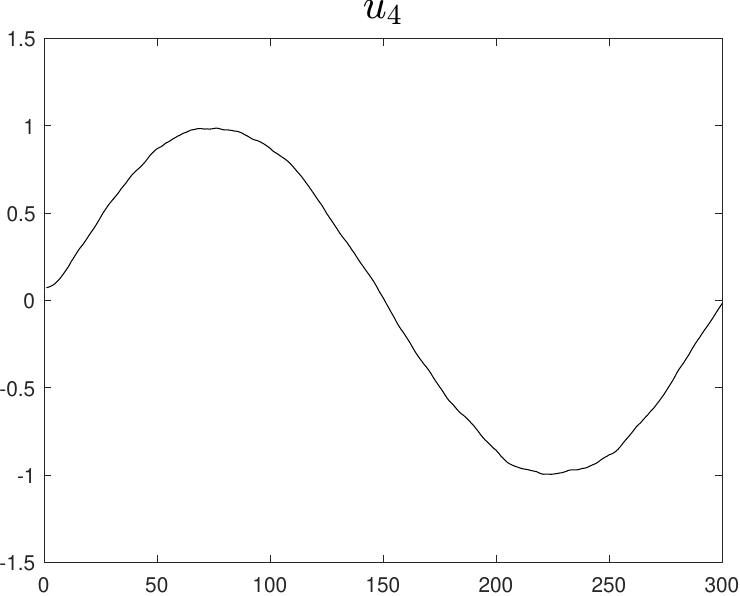}}
\end{minipage}
\begin{minipage}[c]{0.32\textwidth}
    \scalebox{.4}{\includegraphics{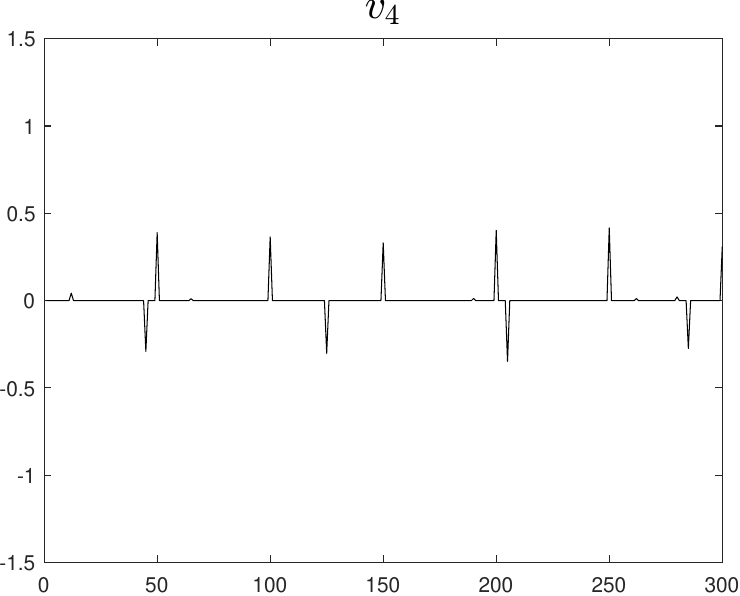}}
\end{minipage}
\begin{minipage}[c]{0.32\textwidth}
    \scalebox{.4}{\includegraphics{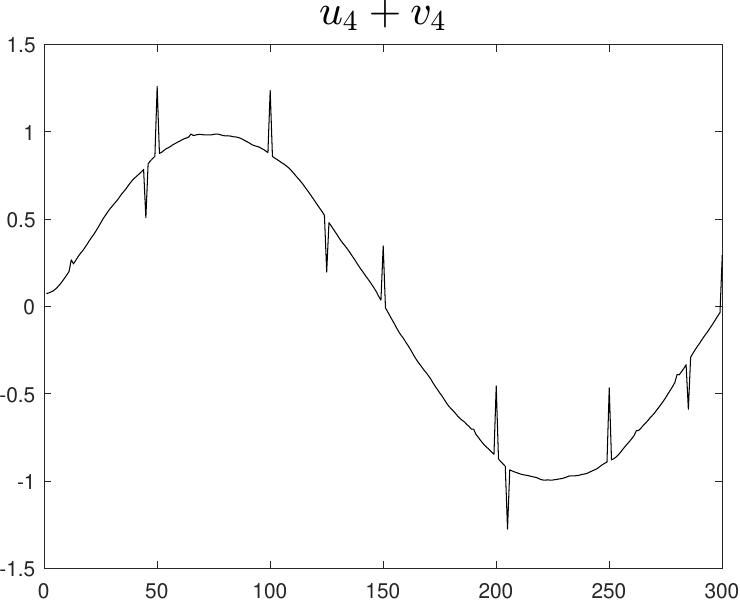}}
\end{minipage}

\vspace{0.5cm}
\centering
\hspace{-1cm}\begin{minipage}[c]{0.32\textwidth}
    \scalebox{.4}{\includegraphics{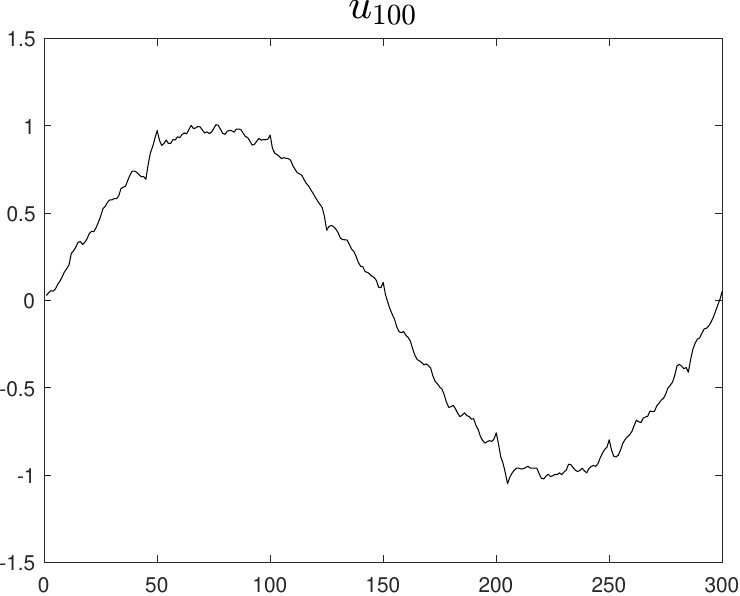}}
\end{minipage}
\begin{minipage}[c]{0.32\textwidth}
    \scalebox{.4}{\includegraphics{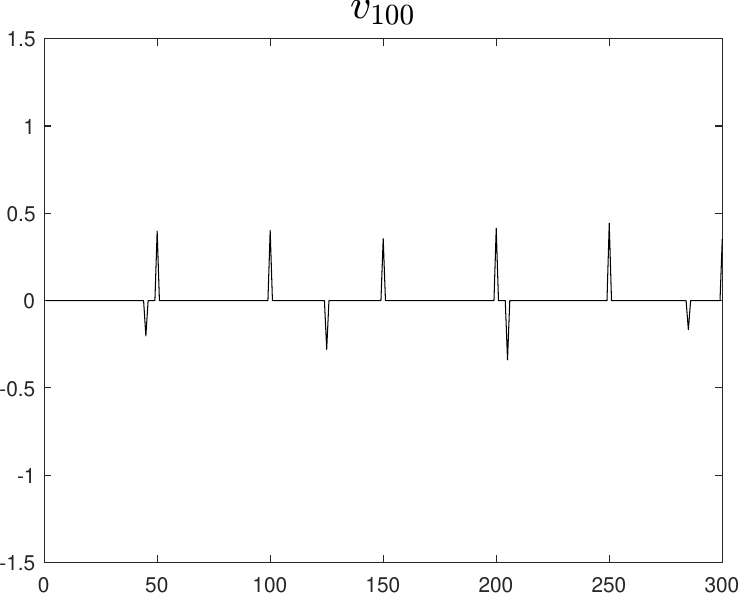}}
\end{minipage}
\begin{minipage}[c]{0.32\textwidth}
    \scalebox{.4}{\includegraphics{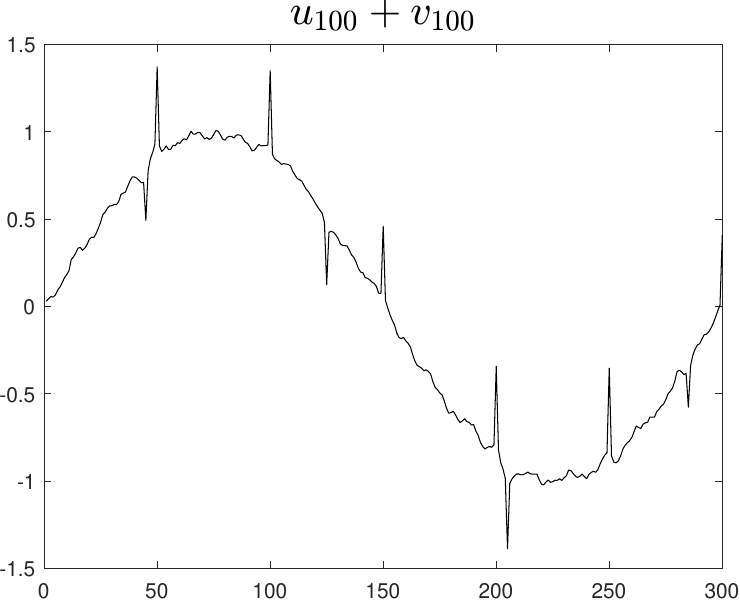}}
\end{minipage}

\caption{Results obtained by Algorithm \ref{algo:NestedBregman_Noisy_Morozov} with $\alpha = {1000}$ and $g,h$ as in \eqref{eq:l1_l2_numrics_functional_morozov}. First line: true components and noisy observation. Lines $2-4$ (top to bottom) reconstructed components for the iteration steps $l = 1,{6}\text{ (minimal cross-correlation)},{4}\text{ (best PSNR)},100$}

\label{fig:nested_morozov_l1_l2}

\end{figure}

\begin{figure}[H] 
\centering
\hspace{-1cm}\begin{minipage}[c]{0.32\textwidth}
    \scalebox{.4}{\includegraphics{Images/Images_H1_l1/u_real.pdf}}
\end{minipage}
\begin{minipage}[c]{0.32\textwidth}
    \scalebox{.4}{\includegraphics{Images/Images_H1_l1/v_real.pdf}}
\end{minipage}
\begin{minipage}[c]{0.32\textwidth}
   \scalebox{.4}{\includegraphics{Images/Images_H1_l1/f.pdf}}
\end{minipage}

\vspace{0.5cm}
\centering
\hspace{-1cm}\begin{minipage}[c]{0.32\textwidth}
    \scalebox{.4}{\includegraphics{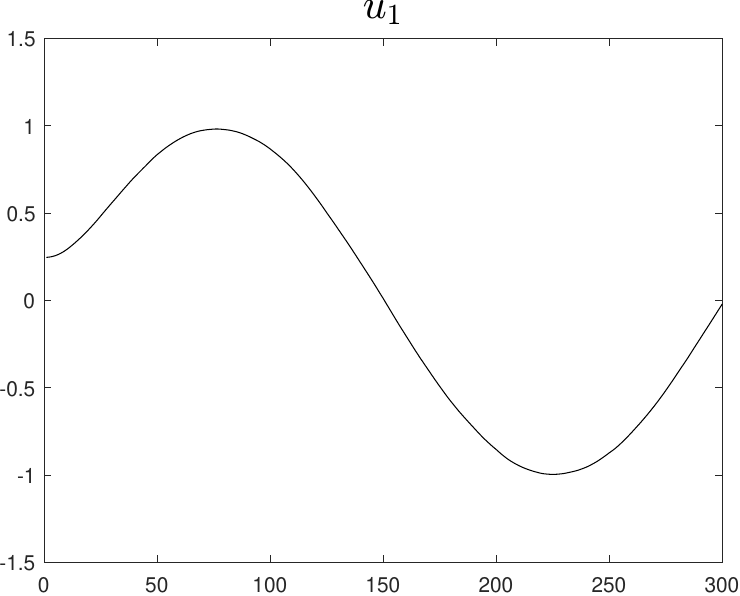}}
\end{minipage}
\begin{minipage}[c]{0.32\textwidth}
    \scalebox{.4}{\includegraphics{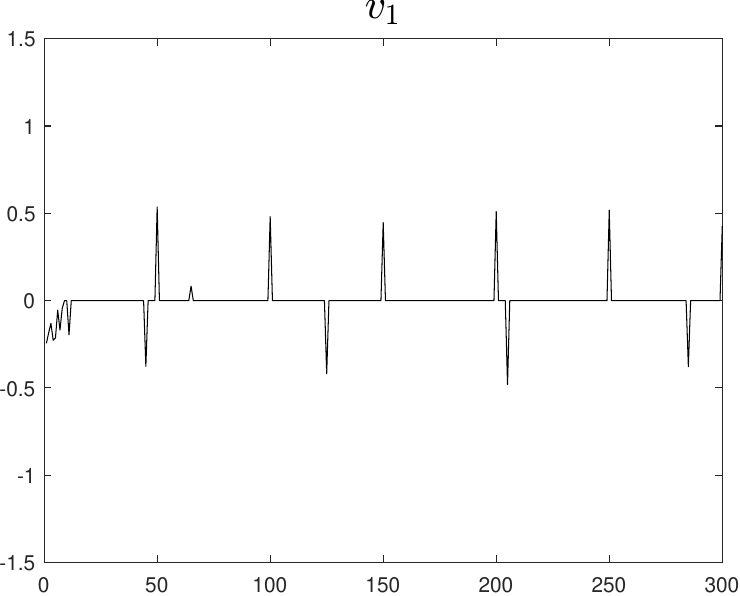}}
\end{minipage}
\begin{minipage}[c]{0.32\textwidth}
    \scalebox{.4}{\includegraphics{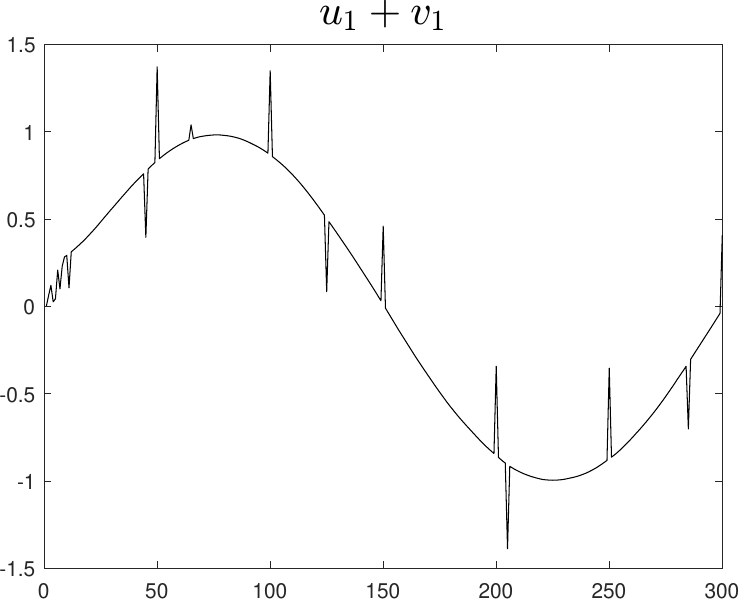}}
\end{minipage}

\vspace{0.5cm}
\centering
\hspace{-1cm}\begin{minipage}[c]{0.32\textwidth}
    \scalebox{.4}{\includegraphics{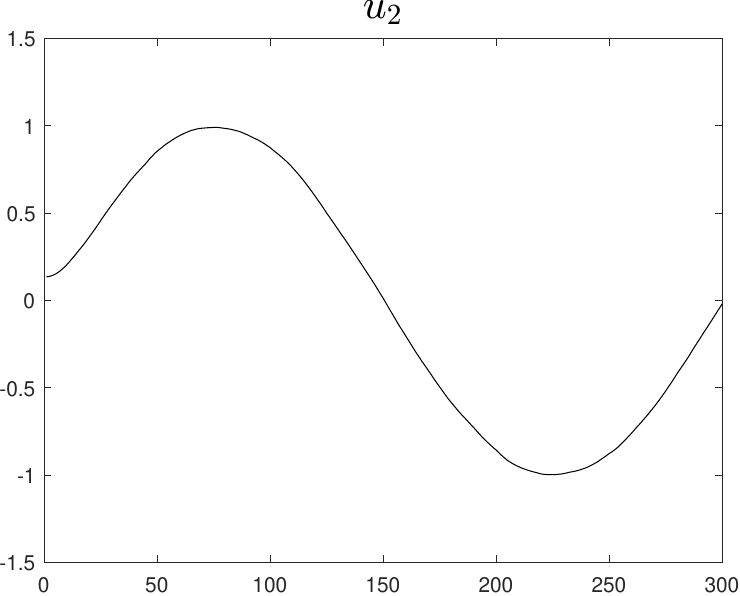}}
\end{minipage}
\begin{minipage}[c]{0.32\textwidth}
    \scalebox{.4}{\includegraphics{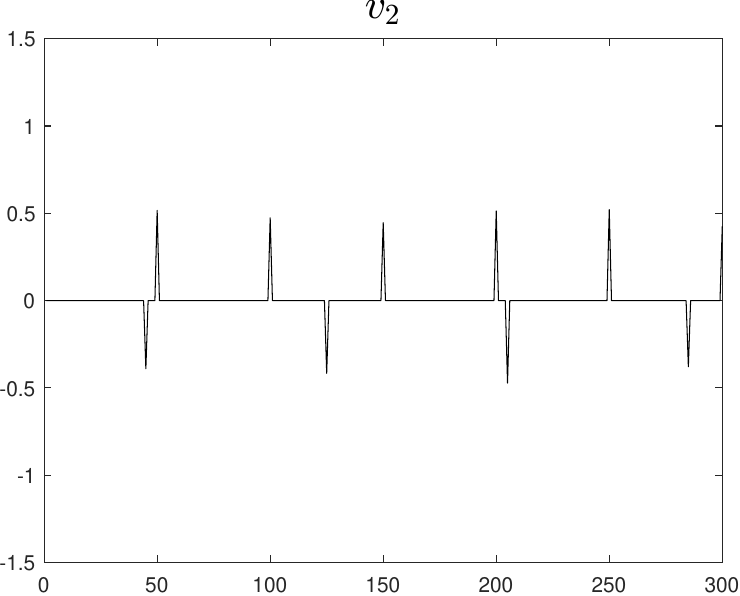}}
\end{minipage}
\begin{minipage}[c]{0.32\textwidth}
    \scalebox{.4}{\includegraphics{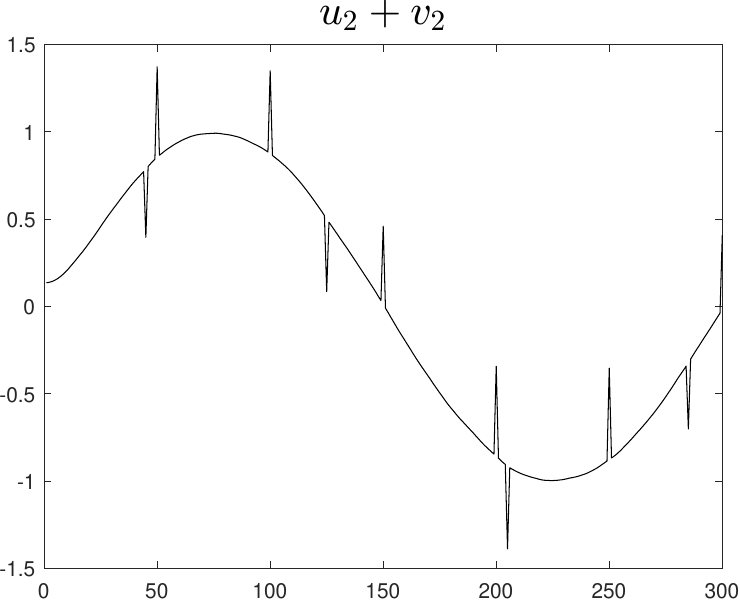}}
\end{minipage}

\vspace{0.5cm}
\centering
\hspace{-1cm}\begin{minipage}[c]{0.32\textwidth}
    \scalebox{.4}{\includegraphics{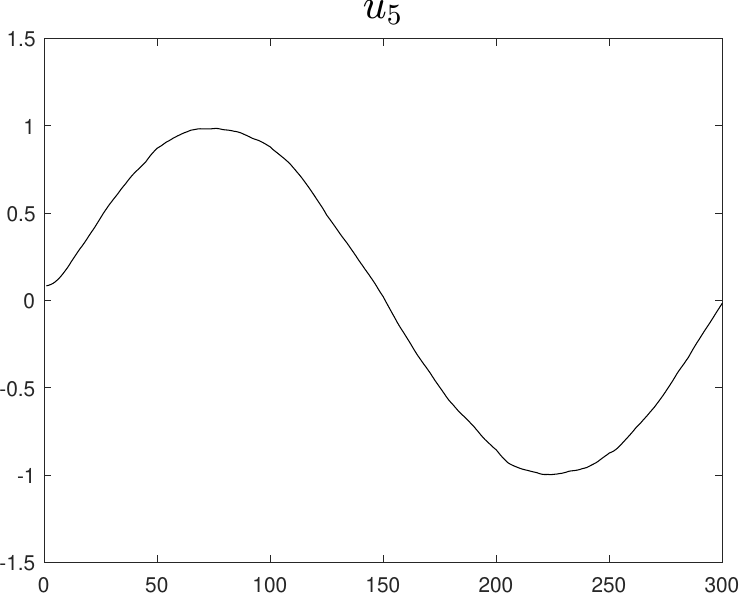}}
\end{minipage}
\begin{minipage}[c]{0.32\textwidth}
    \scalebox{.4}{\includegraphics{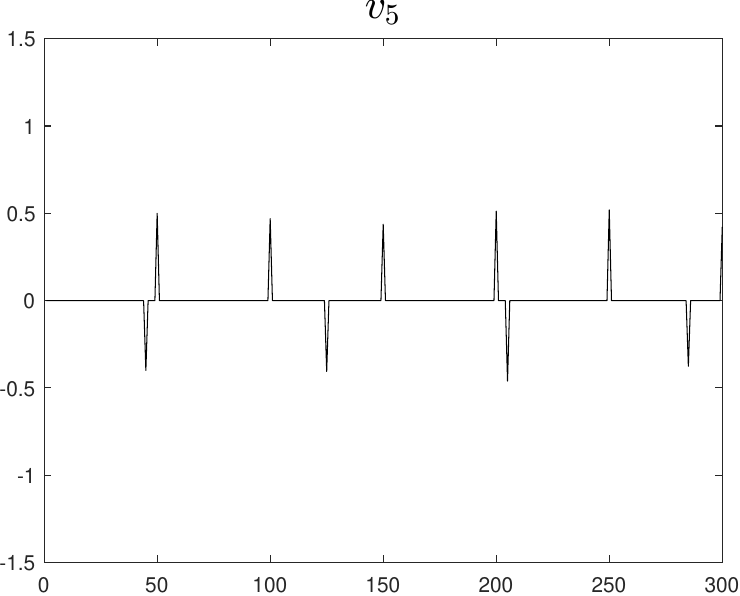}}
\end{minipage}
\begin{minipage}[c]{0.32\textwidth}
    \scalebox{.4}{\includegraphics{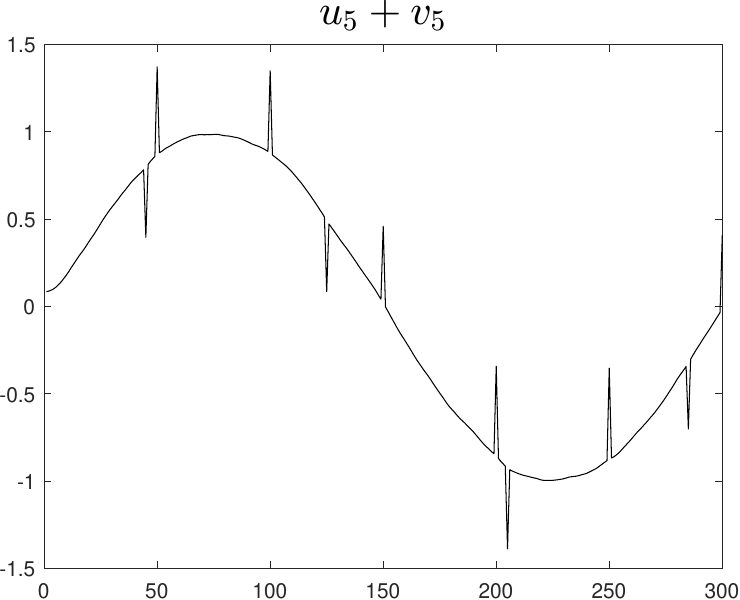}}
\end{minipage}

\vspace{0.5cm}
\centering
\hspace{-1cm}\begin{minipage}[c]{0.32\textwidth}
    \scalebox{.4}{\includegraphics{Images/Images_H1_l1/Images_Nested/u_100.pdf}}
\end{minipage}
\begin{minipage}[c]{0.32\textwidth}
    \scalebox{.4}{\includegraphics{Images/Images_H1_l1/Images_Nested/v_100.pdf}}
\end{minipage}
\begin{minipage}[c]{0.32\textwidth}
    \scalebox{.4}{\includegraphics{Images/Images_H1_l1/Images_Nested/u_100+v_100.pdf}}
\end{minipage}
\caption{Results obtained by Algorithm \ref{algo:NestedBregman_Noisy} with $\alpha ={1000}$ and $g,h$ as in \eqref{eq:l1_l2_numrics_functional_morozov}. First line: true components and noisy observation. Lines $2-4$ (top to bottom) reconstructed components for the iteration steps $l = 1,2\text{ (minimal cross-correlation)},{5}\text{ (best PSNR)},100$}

\label{fig:nested_bregman_l1_l2}

\end{figure}

\begin{figure}[H]
    \centering
    \scalebox{.5}{\includegraphics{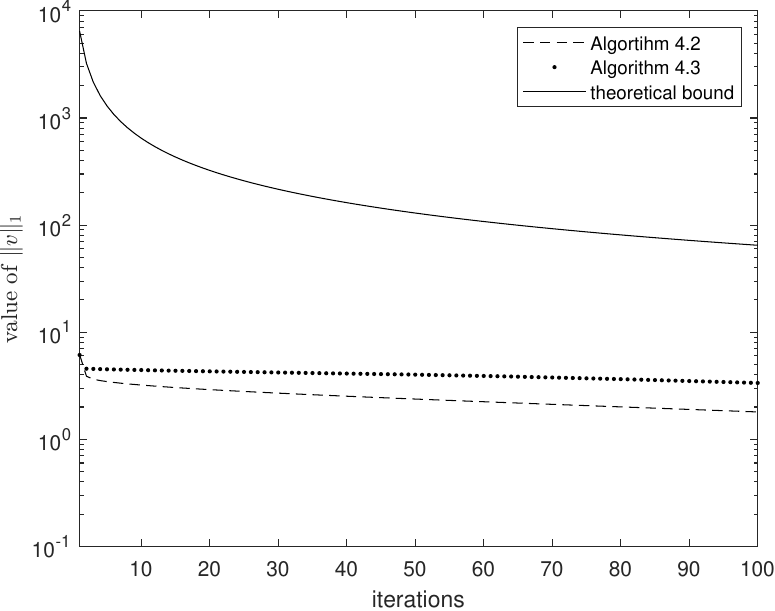}}
    \caption{Values of $h(v_l^\delta)$ in Algorithms \ref{algo:NestedBregman_Noisy_Morozov} and \ref{algo:NestedBregman_Noisy} vs the theoretical upper-bound from \eqref{eq:Convergence_h_noisy}.}
    \label{fig:nested_morozov_l1_l2_bound}
\end{figure}
}{}

\subsection{\texorpdfstring{Experiment 2: $TV-H^1$-decomposition}{TEXT}}
\begin{figure}[H]
    \centering
    \begin{minipage}[c]{0.49\textwidth}
        \includegraphics[scale = 0.5]{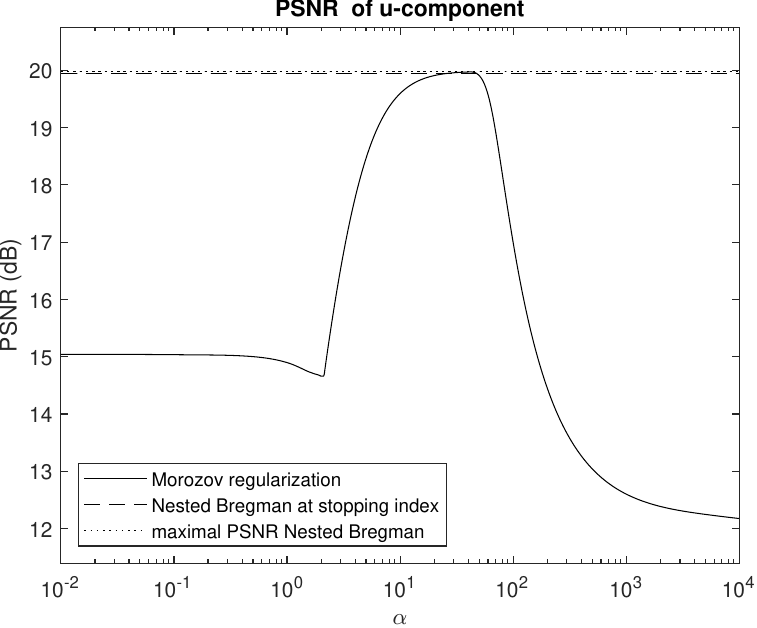}
    \end{minipage}
        \begin{minipage}[c]{0.49\textwidth}
        \includegraphics[scale = 0.5]{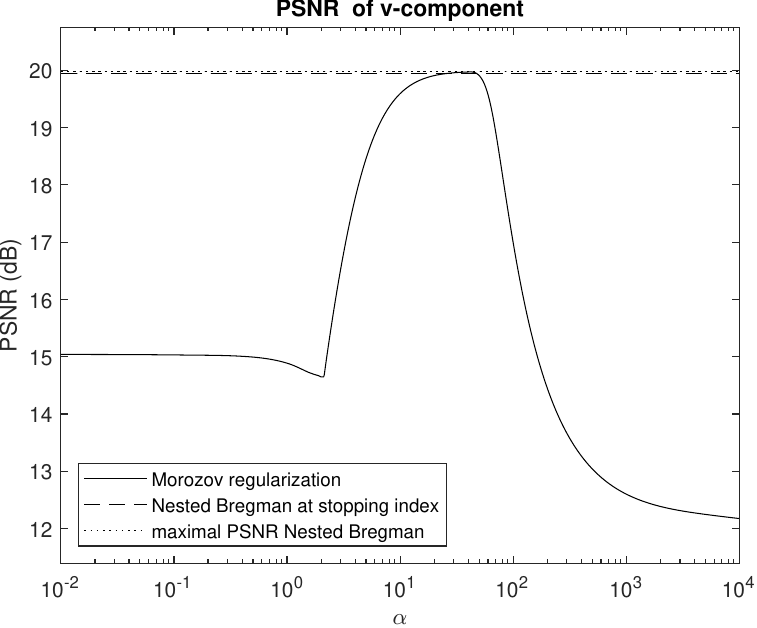}
    \end{minipage}
\vspace{0.7em}

        \begin{minipage}[c]{0.49\textwidth}
        \includegraphics[scale = 0.5]{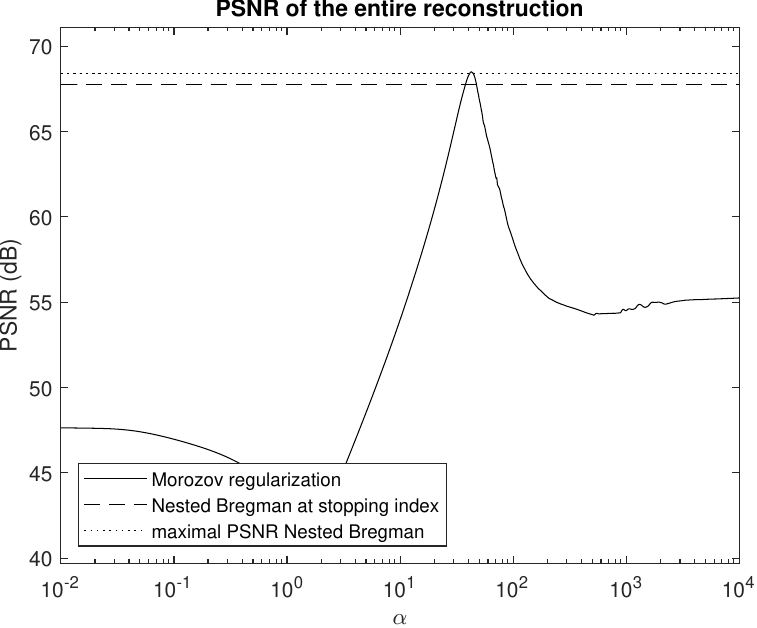}
    \end{minipage}
        \begin{minipage}[c]{0.49\textwidth}
        \includegraphics[scale = 0.5]{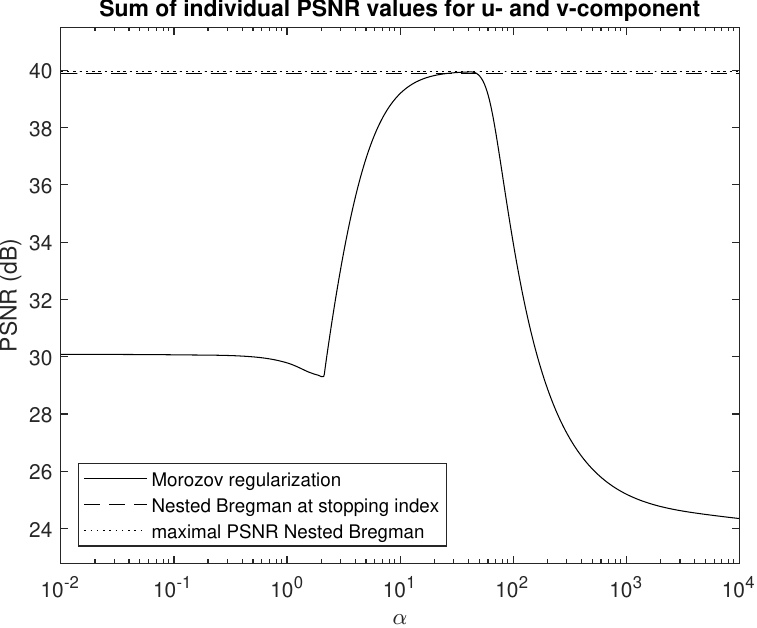}
    \end{minipage}
    \caption{{PSNR values for Morozov regularization for different values of $\alpha$, Nested Bregman iterations with weight $\alpha = 1000$ at the stopping index and maximal PSNR from Nested Bregman iterations. From top left to bottom right: $u-$component, $v-$component, full reconstruction , sum of $u-$component and $v-$component PSNR values.}}
    \label{fig:comp_Morozov_H1_TV}
\end{figure}

We revisit the example from Section \ref{sec:regularization}. For the ratio $\frac{\alpha}{\beta} = 1000$ (Line 2 in Figures \ref{fig:variational_deblurring} and \ref{fig:bregman_deblurring}), we observed that the total variation component $v$ was overrepresented in both the variational regularization and Bregman iterations. We therefore use Algorithm \ref{algo:NestedBregman_Noisy_Morozov} and Algorithm \ref{algo:NestedBregman_Noisy} to decrease its share in the decomposition. In the latter case, we let  $\alpha = 465$ and $\beta = 0.465$, to be consistent with the results presented in Figure \ref{fig:bregman_deblurring}. Figure \ref{fig:PSNR_H1-TV} shows the sum of the PSNR values for the individual components obtained by Algorithm \ref{algo:NestedBregman_Noisy_Morozov} (left) and Algorithm \ref{algo:NestedBregman_Noisy} (right) with the first local minimum of the scalar cross-correlation marked by an asterisk. When considering Morozov regularization in the inner loop, the cross-correlation based stopping rule actually suggests terminating the iteration at the exact maximum of the cumulative PSNR. When using Bregman iterations as an inner loop, the suggested stopping index is  $l = 21$, while the maximal PSNR is attained at $l = 22$. Therefore,  both experiments employing this stopping rule yield a very good approximation of the true decomposition. The corresponding behavior can be seen in Figures \ref{fig:Nested_Morozov_ratio=0.001} and \ref{fig:Nested_Bregman_ratio=0.001} (from top to bottom), where the true component and the observation, the first, the minimal scalar cross-correlation, the maximal cumulative PSNR {iterate, respectively}, {as well as the one} at $l=50$ are depicted for both algorithms. As in the case of $L^1-H^1$-decompositions, we note that the iterates obtained by Algorithm \ref{algo:NestedBregman_Noisy} yield better approximations $x_l$ of the true phenomenon with respect to the PSNR. {We observe similar results, when employing Algorithm \ref{algo:NestedBregman_Noiseless} for noise-free data.}\\
{As in the previous part, we also compare the results obtained from Algorithm \ref{algo:NestedBregman_Noisy_Morozov} to the ones obtained by single-step Morozov regularization with penalty term $J(x) = \frac{\alpha}{2}\norm{\nabla \cdot }_{L^2}^2\square \tv{\cdot}$ for $1000$ logarithmically spaced values of $\alpha$ ranging from $10^{-2}$ to $10^4$. The results of these experiments can be seen in Figure \ref{fig:comp_Morozov_H1_TV}. We observe that the results of our method are only marginally smaller than the optimal ones from the single-step regularization. Furthermore, they are outperformed only within a relatively small range for the weight $\alpha$. The precise values of $\alpha$ and corresponding PSNR values can be seen in Table \ref{tab:Nested_vs_Morozov_H1_TV}. This means the proposed method produces reconstructions that are comparable to what can be obtained by Morozov regularization with optimally chosen parameters. Hence, the reconstructions via Nested Bregman iterations perform as well as an optimal bilevel method. We also argue that the Nested Bregman iterations are preferable to a grid search with respect to the weighting parameter, as such a method would require choosing a grid containing those values of $\alpha$ for which the single-step methods yield accurate reconstructions. However, this is rather challenging without any further analysis of the decomposition problem.}

\begin{table}
    \centering
    \begin{tabular}{|c|c|c|c|c|c|}\hline
         & {$\text{PSNR}_{si}$}  & {$\text{PSNR}_{si}\le \text{PSNR}_{ss}$} & $\text{PSNR}_{ne}^{\max}$ & {$\text{PSNR}_{ne}\le \text{PSNR}_{ss}$} & {$\text{PSNR}_{ss}^{\max}$}\\  \noalign{\hrule height 1.5pt}
        $u$ & $19.944$ & $26.147\le \alpha \le 47.389$  & $19.971$ & $41.268\le \alpha \le 41.843$&$19.971$\\ \hline
        $v$& $19.944$ & $26.147\le \alpha \le 47.389$ & $19.971$ &$41.268\le \alpha \le 41.843$ &$19.971$\\ \hline
       $x$ & $67.770$ & $38.511\le \alpha \le 47.389$  & $68.378$ &$41.268\le \alpha \le 44.843$ &$68.488$\\ \hline
         $u+v$& $39.889$ & $26.147\le \alpha \le 47.389$ & $39.942$ &$41.268\le \alpha \le 41.843$ &$39.942$\\ \hline
    \end{tabular}
    \caption{{Left to right: PSNR with Nested Bregman iteration at stopping index (si), weighting choices for which single-step (ss) outperforms Nested Bregman at stopping index, maximal PSNR of Nested Bregman iteration (ne), weighting choices for which single-step outperforms optimal Nested Bregman, maximal PSNR of single-step. \\
    top to bottom: $u-$component, $v-$component, full reconstruction, sum of the PSNR values of the $u-$ and $v-$ component}}
    \label{tab:Nested_vs_Morozov_H1_TV}
\end{table}

{\begin{figure}[H]
    \centering
  \begin{minipage}[c]{0.49\textwidth}
    \includegraphics[scale = 0.5]{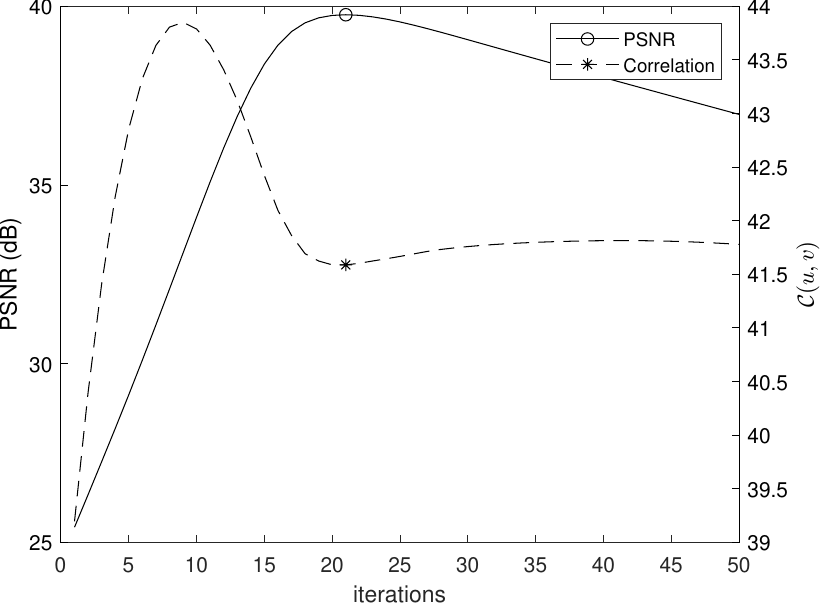}
\end{minipage}
\begin{minipage}[c]{0.49\textwidth}
    \includegraphics[scale = 0.5]{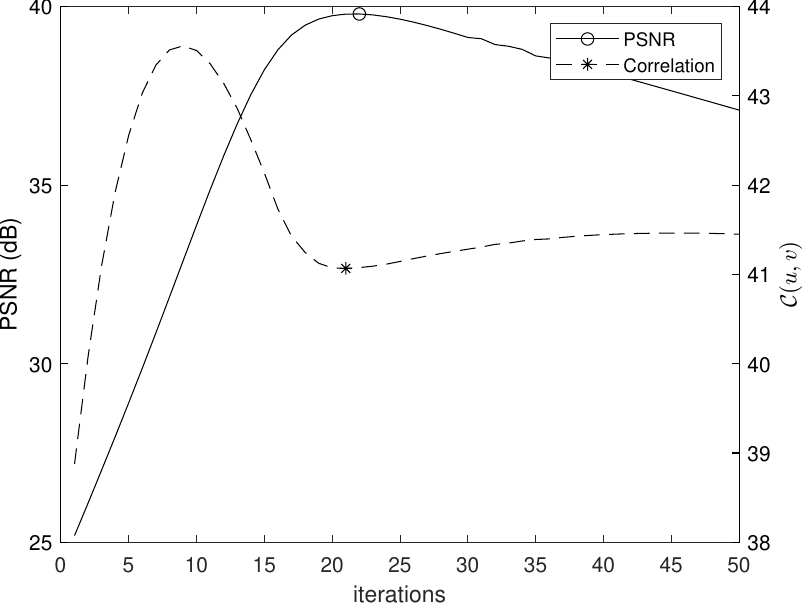}
\end{minipage}
\caption{Sum of the componentwise PSNR values obtained by Algorithm \ref{algo:NestedBregman_Noisy_Morozov} (left) and Algorithm \ref{algo:NestedBregman_Noisy} (right). The circle and the asterisk mark the maximum of the PSNR and the  first local minimum of the scalar normalized cross-correlation, respectively.}
\label{fig:PSNR_H1-TV}

\end{figure}

\begin{figure}[H]
\centering
\hspace{-1cm}\begin{minipage}[c]{0.24\textwidth}
    \includegraphics[scale = 0.4]{Images/Images_H1_TV/Images_Regularization/u.pdf}
\end{minipage}
\begin{minipage}[c]{0.24\textwidth}
    \includegraphics[scale = 0.4]{Images/Images_H1_TV/Images_Regularization/v.pdf}
\end{minipage}
\begin{minipage}[c]{0.24\textwidth}
    \includegraphics[scale = 0.4]{Images/Images_H1_TV/Images_Regularization/u+v.pdf}
\end{minipage}
\begin{minipage}[c]{0.24\textwidth}
    \includegraphics[scale = 0.4]{Images/Images_H1_TV/Images_Regularization/f.pdf}
\end{minipage}

\vspace{0.5cm}
\centering
\hspace{-1cm}\begin{minipage}[c]{0.24\textwidth}
    \includegraphics[scale = 0.4]{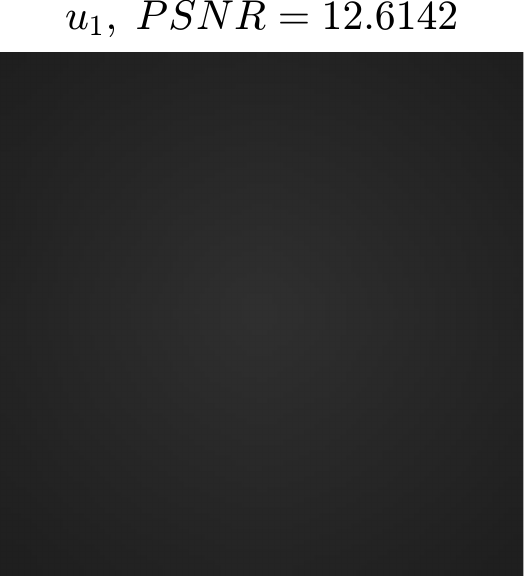}
\end{minipage}
\begin{minipage}[c]{0.24\textwidth}
    \includegraphics[scale = 0.4]{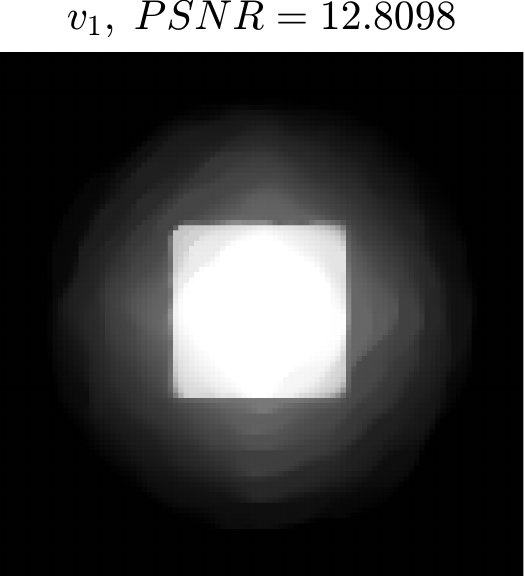}
\end{minipage}
\begin{minipage}[c]{0.24\textwidth}
    \includegraphics[scale = 0.4]{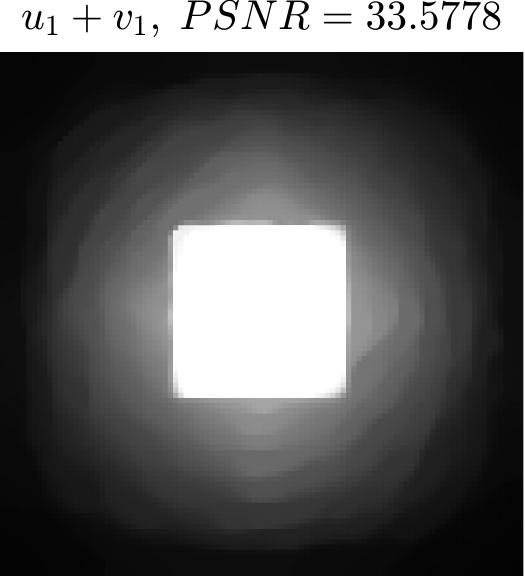}
\end{minipage}
\begin{minipage}[c]{0.24\textwidth}
    \includegraphics[scale = 0.4]{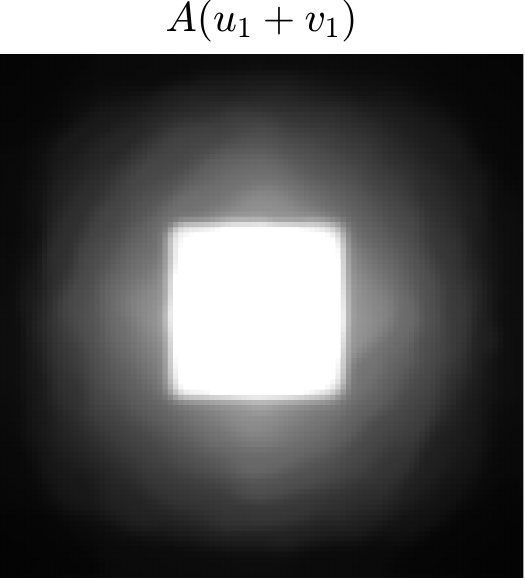}
\end{minipage}

\vspace{0.5cm}
\centering
\hspace{-1cm}\begin{minipage}[c]{0.24\textwidth}
    \includegraphics[scale = 0.4]{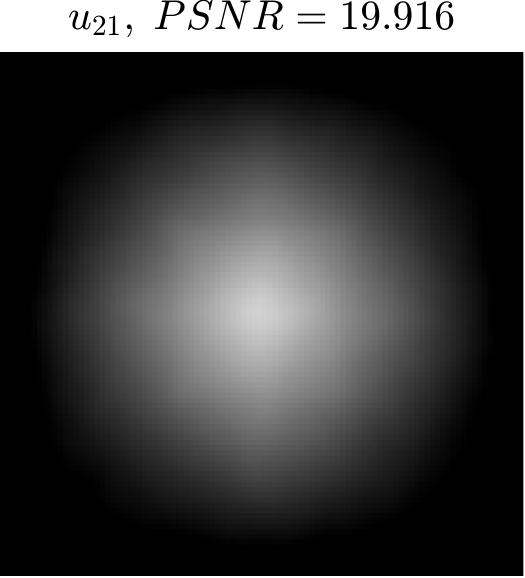}
\end{minipage}
\begin{minipage}[c]{0.24\textwidth}
    \includegraphics[scale = 0.4]{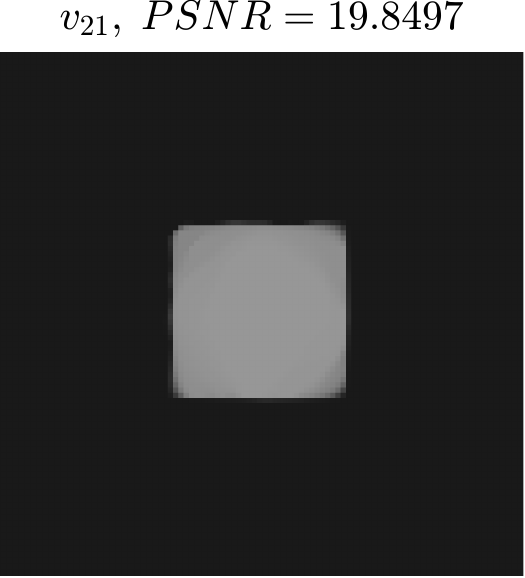}
\end{minipage}
\begin{minipage}[c]{0.24\textwidth}
    \includegraphics[scale = 0.4]{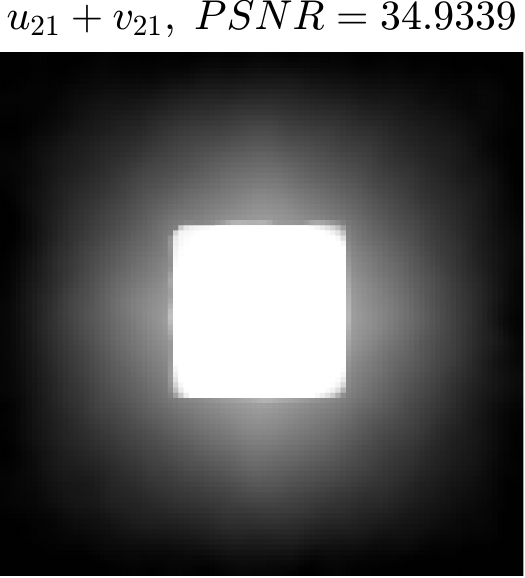}
\end{minipage}
\begin{minipage}[c]{0.24\textwidth}
    \includegraphics[scale = 0.4]{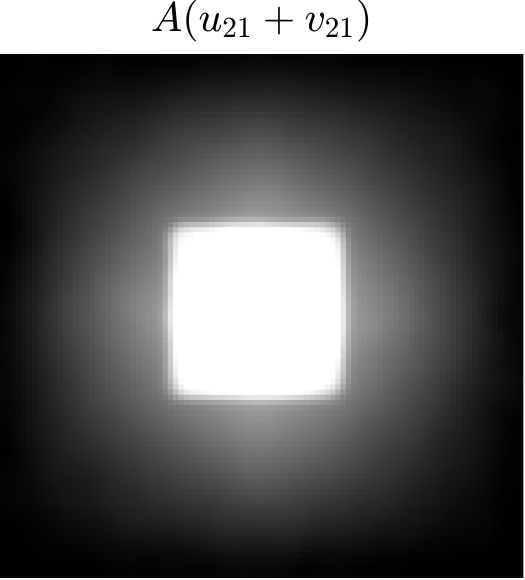}
\end{minipage}

\vspace{0.5cm}
\centering
\hspace{-1cm}\begin{minipage}[c]{0.24\textwidth}
    \includegraphics[scale = 0.4]{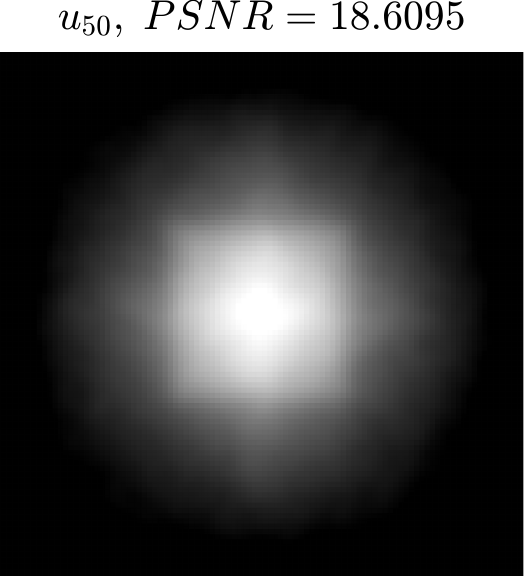}
\end{minipage}
\begin{minipage}[c]{0.24\textwidth}
    \includegraphics[scale = 0.4]{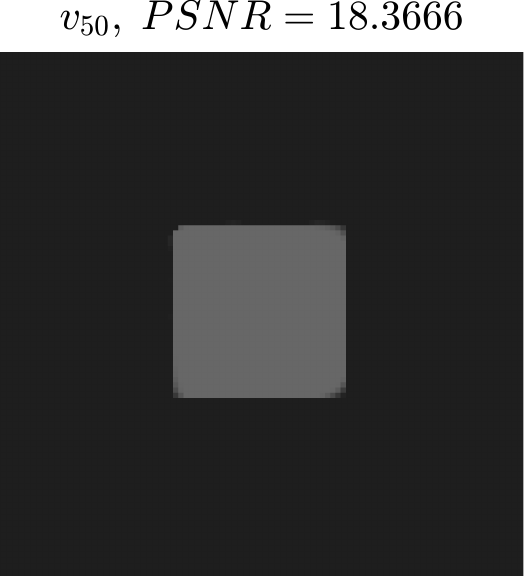}
\end{minipage}
\begin{minipage}[c]{0.24\textwidth}
    \includegraphics[scale = 0.4]{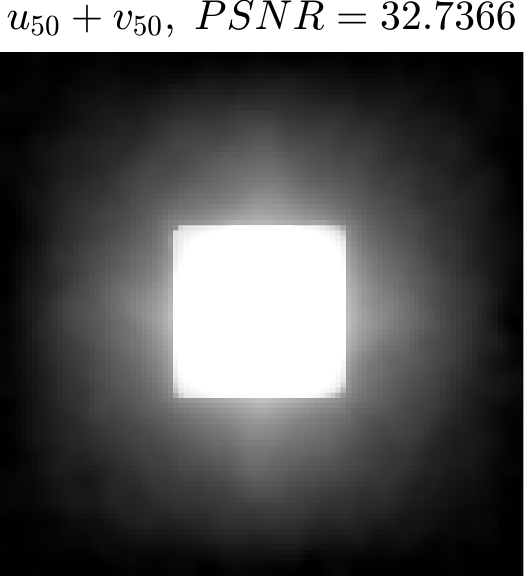}
\end{minipage}
\begin{minipage}[c]{0.24\textwidth}
    \includegraphics[scale = 0.4]{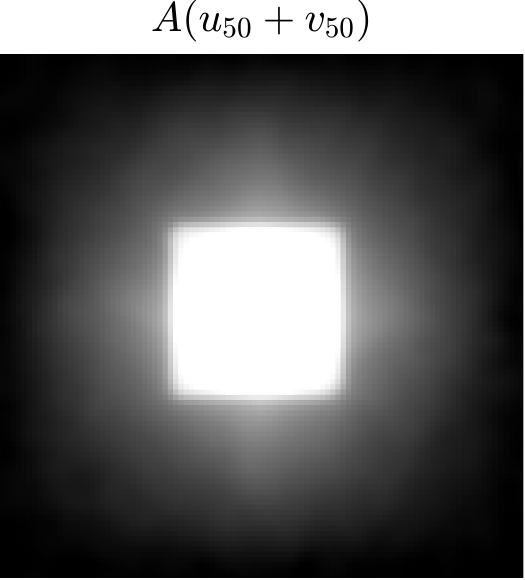}
\end{minipage}

\caption{ Results obtained  by Algorithm \ref{algo:NestedBregman_Noisy_Morozov} with $\alpha =1000$, $\beta = 1$ and $g,h$ as in \eqref{eq:functionals_TV_H^1}. First line: true components and noisy observation. Lines $2-4$ (top to bottom) reconstructed components for the iteration steps $l = 1,21\text{ (minimal cross-correlation and best PSNR)},50$}

\label{fig:Nested_Morozov_ratio=0.001}

\end{figure}
%-------------------------------------------------now Bregman------------------------------

\begin{figure}[H]
\centering
\hspace{-1cm}\begin{minipage}[c]{0.24\textwidth}
    \includegraphics[scale = 0.4]{Images/Images_H1_TV/Images_Regularization/u.pdf}
\end{minipage}
\begin{minipage}[c]{0.24\textwidth}
    \includegraphics[scale = 0.4]{Images/Images_H1_TV/Images_Regularization/v.pdf}
\end{minipage}
\begin{minipage}[c]{0.24\textwidth}
    \includegraphics[scale = 0.4]{Images/Images_H1_TV/Images_Regularization/u+v.pdf}
\end{minipage}
\begin{minipage}[c]{0.24\textwidth}
    \includegraphics[scale = 0.4]{Images/Images_H1_TV/Images_Regularization/f.pdf}
\end{minipage}

\vspace{0.5cm}
\centering
\hspace{-1cm}\begin{minipage}[c]{0.24\textwidth}
    \includegraphics[scale = 0.4]{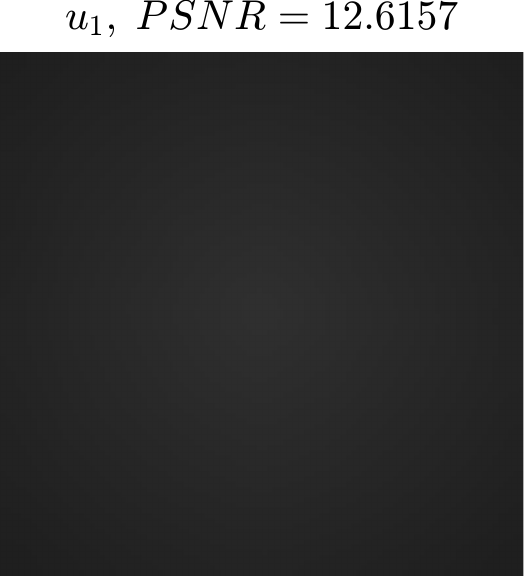}
\end{minipage}
\begin{minipage}[c]{0.24\textwidth}
    \includegraphics[scale = 0.4]{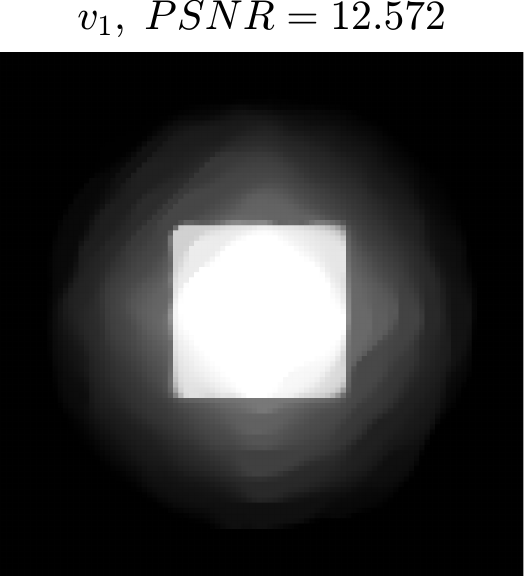}
\end{minipage}
\begin{minipage}[c]{0.24\textwidth}
    \includegraphics[scale = 0.4]{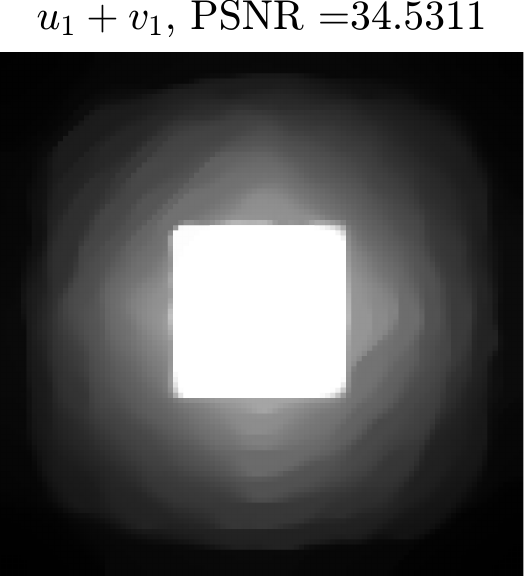}
\end{minipage}
\begin{minipage}[c]{0.24\textwidth}
    \includegraphics[scale = 0.4]{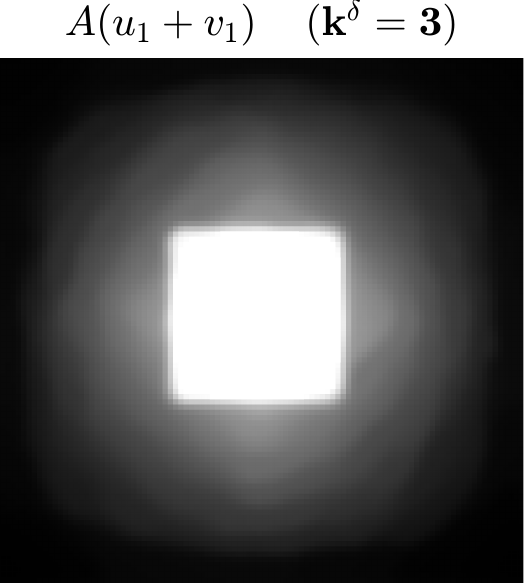}
\end{minipage}

\vspace{0.5cm}
\centering
\hspace{-1cm}\begin{minipage}[c]{0.24\textwidth}
    \includegraphics[scale = 0.4]{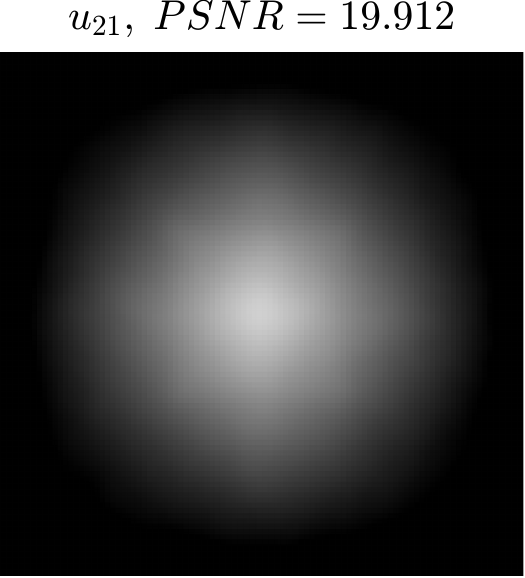}
\end{minipage}
\begin{minipage}[c]{0.24\textwidth}
    \includegraphics[scale = 0.4]{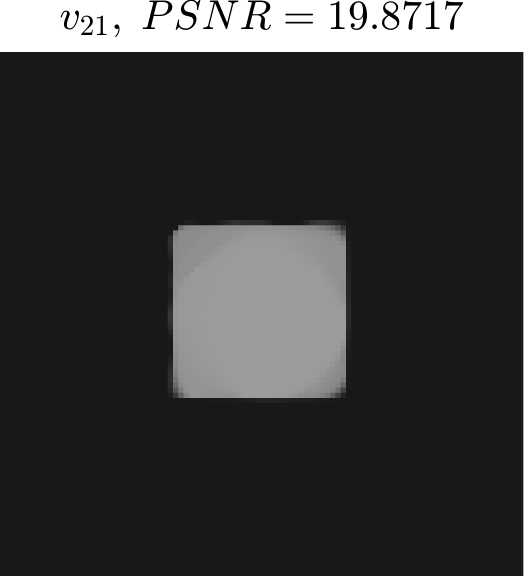}
\end{minipage}
\begin{minipage}[c]{0.24\textwidth}
    \includegraphics[scale = 0.4]{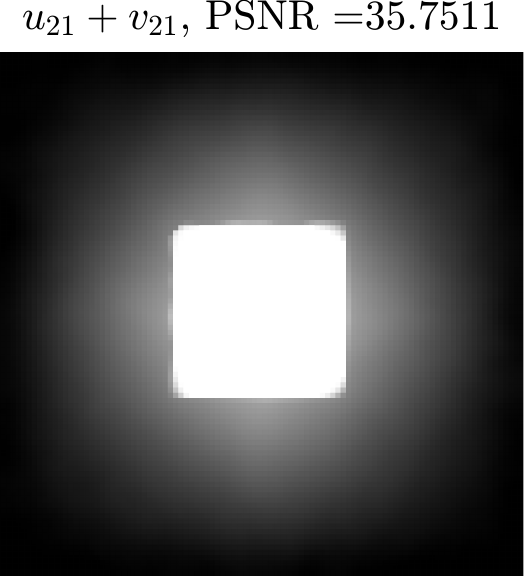}
\end{minipage}
\begin{minipage}[c]{0.24\textwidth}
    \includegraphics[scale = 0.4]{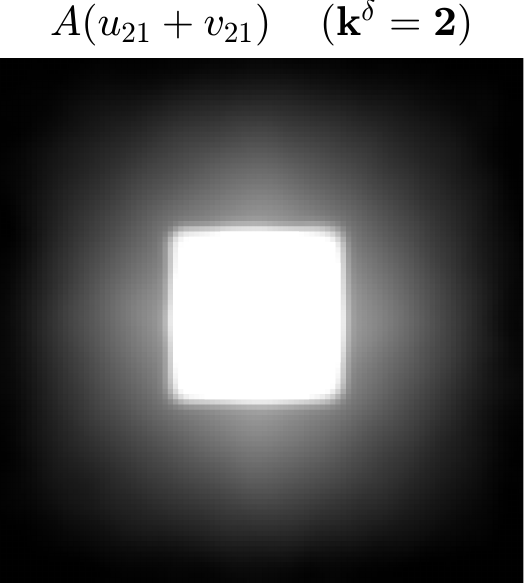}
\end{minipage}

\vspace{0.5cm}
\centering
\hspace{-1cm}\begin{minipage}[c]{0.24\textwidth}
    \includegraphics[scale = 0.4]{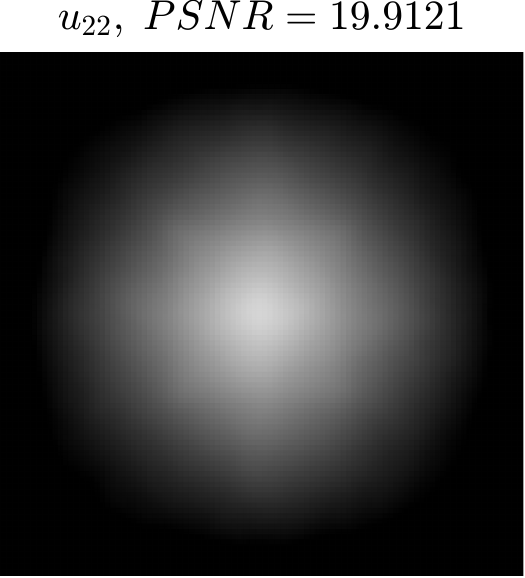}
\end{minipage}
\begin{minipage}[c]{0.24\textwidth}
    \includegraphics[scale = 0.4]{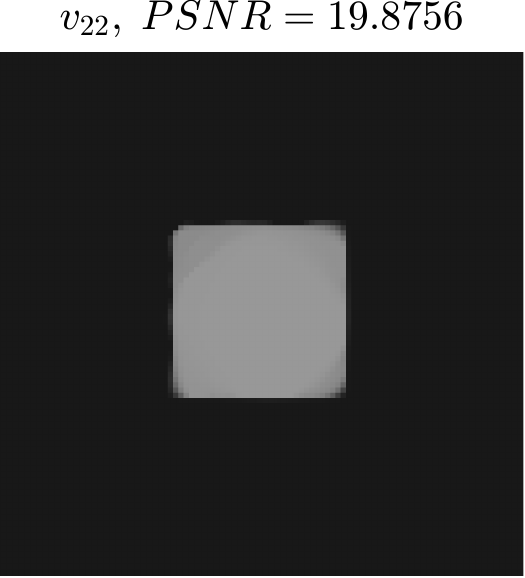}
\end{minipage}
\begin{minipage}[c]{0.24\textwidth}
    \includegraphics[scale = 0.4]{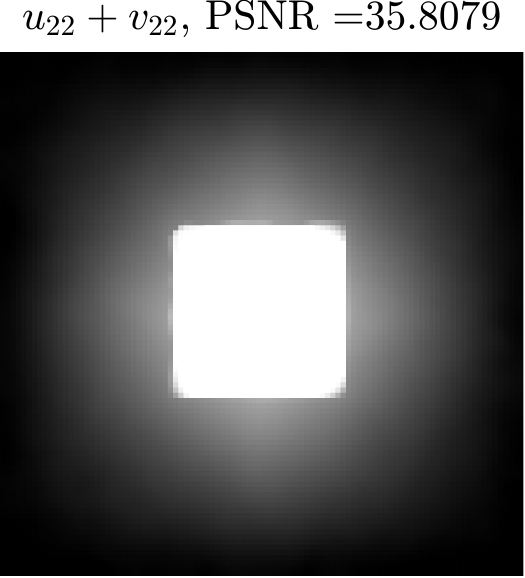}
\end{minipage}
\begin{minipage}[c]{0.24\textwidth}
    \includegraphics[scale = 0.4]{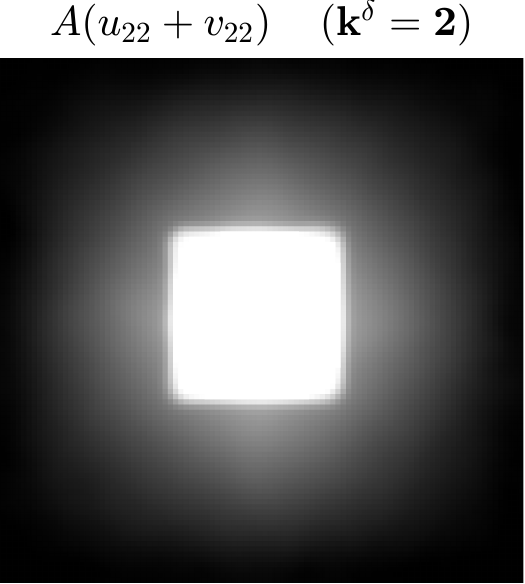}
\end{minipage}

\vspace{0.5cm}
\centering
\hspace{-1cm}\begin{minipage}[c]{0.24\textwidth}
    \includegraphics[scale = 0.4]{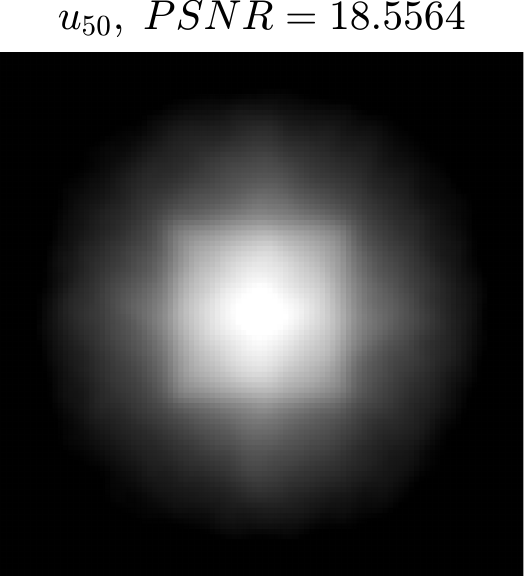}
\end{minipage}
\begin{minipage}[c]{0.24\textwidth}
    \includegraphics[scale = 0.4]{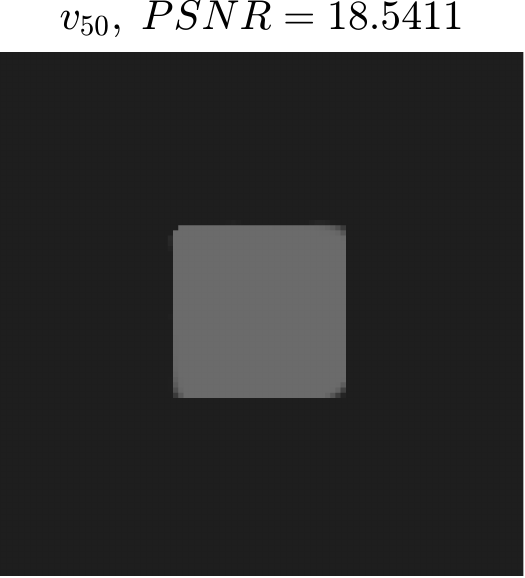}
\end{minipage}
\begin{minipage}[c]{0.24\textwidth}
    \includegraphics[scale = 0.4]{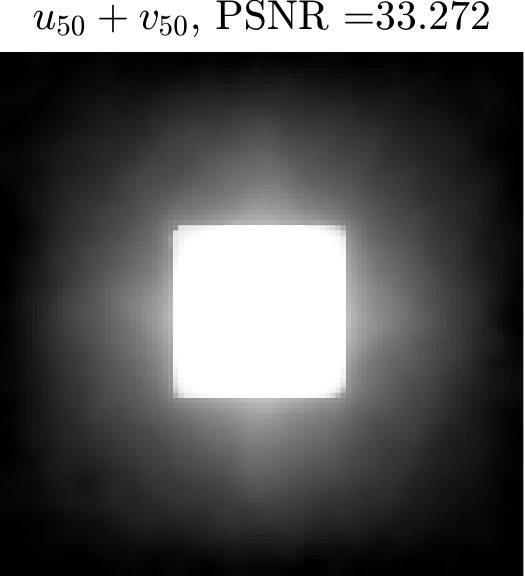}
\end{minipage}
\begin{minipage}[c]{0.24\textwidth}
    \includegraphics[scale = 0.4]{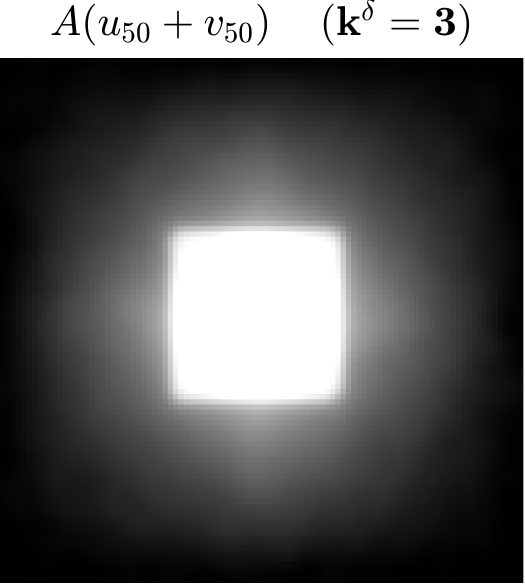}
\end{minipage}

\caption{Results obtained by Algorithm \ref{algo:NestedBregman_Noisy} with $\alpha =464$, $\beta = 465$ and $g,h$ as in \eqref{eq:functionals_TV_H^1}. First line: true components and noisy observation. Lines $2-4$ (top to bottom) reconstructed components for the iteration steps $l = 1,21\text{ (minimal cross-correlation)},22\text{ (best PSNR)},50$}

\label{fig:Nested_Bregman_ratio=0.001}

\end{figure}

}{}

\subsection{Experiment 3: $TGV-TGV^{osci}$-decompositions}
In order to show the effectiveness of our method for more complex regularizers, we conclude the numerical examples with the infimal convolution of the total generalized variation and its oscillatory version. That is, we choose $g$ and $h$ in the Nested Bregman algorithms as in \eqref{eq:TGV} and \eqref{eq:TGV_osci}. Since even in the discretized setting, none of these functionals is differentiable, we replace the $1-$norm by its Moreau-envelope.  As a test image, we consider the sum of two components $x^\dagger = u^\dagger +v^\dagger$, where $u^\dagger$ is a square with piecewise affine texture on an empty background and $v^\dagger$is an oscillating pattern with an empty square in its center (see line 1 in Figure \ref{fig:TGV_osci_Nested_Morozov}). To construct this image, we consider two components consisting of nine sub-squares with $75\times 75$ pixels each. In the first component $u^\dagger$, only the central square is described by an affine function, while the pixel values of the other squares are $0$. In the other components, the central square is void, while the other squares contain some oscillating texture, whose pixel values follow the definition  $\cos(\omega\cdot x) + \sin(\omega\cdot x)$ with $x  \in [0,75]^2$ for $\omega = (0.25,0.5)^T$. The observation $f^\delta$ was created by adding Gaussian noise with mean $\mu  = 0$ and variance $\sigma = 0.05$ to the sum $u^\dagger+v^\dagger$. We  run Algorithm \ref{algo:NestedBregman_Noisy_Morozov} with $g(u) = TGV_{\alpha_1,\beta_1}(u)$ and $h(v) = TGV^{osci}_{\alpha_2,\beta_2,C}(v)$ with $C = \omega^T \omega$. As initial guess for the parameters, we consider $\alpha_1 = 5$ and $\alpha_2 = 1$. Following the numerical experiments used to obtain Figure 5 in \cite{oscICTV}, we choose $\beta_i = \alpha_i$ for $i = 1,2$. As in the previous experiments, the sum of the componentwise PSNR and the stopping index proposed by the scalar cross-correlation can be seen in Figure \ref{fig:PSNR_TGV_osci}. The minimal cross-correlation was achieved at iteration step $l = 3$ and the best PSNR at step $l = 4$. The iterates at steps $l = 1,3,4,10$ can be seen in lines $2-5$ of Figure \ref{fig:TGV_osci_Nested_Morozov}. We observe that the piecewise affine component in the first iteration is over-regularized,  some of its features being contained in the oscillatory component. The iterates performing best with respect to cross-correlation (line 3) and as well as those with the maximal PSNR (line 4) yield good reconstructions of the true components, with only minor inaccuracies around the edges of the middle square. However, continuing the iteration further leads to an over-regularization of the oscillatory component and, consequently, to oscillatory features occurring in the piecewise affine component (line 5).

{

\begin{figure}[H][H]
    \centering
    \includegraphics[scale = 0.6]{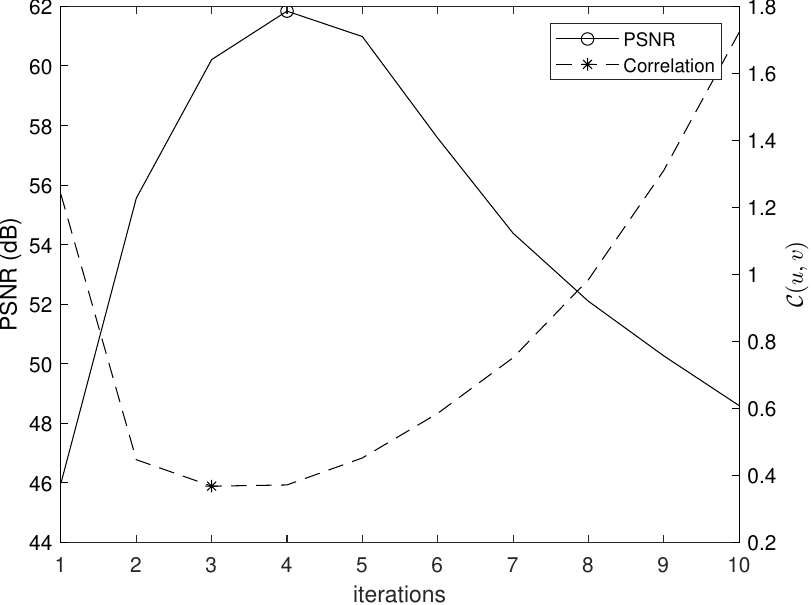}
    \caption{Sum of the componentwise PSNR values obtained by Algorithm \ref{algo:NestedBregman_Noisy_Morozov}. The circle and the asterisk mark the maximum of the PSNR and the  first local minimum of the scalar normalized cross-correlation, respectively.}
    \label{fig:PSNR_TGV_osci}
\end{figure}

\begin{figure}[H]
\centering
\hspace{-1cm}\begin{minipage}[c]{0.24\textwidth}
    \includegraphics[scale = 0.4]{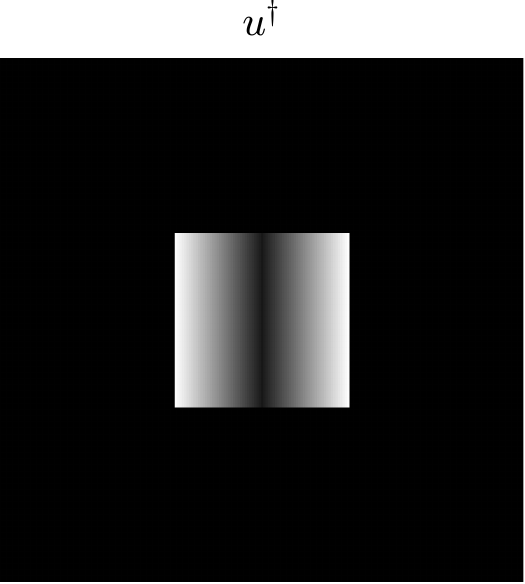}
\end{minipage}
\begin{minipage}[c]{0.24\textwidth}
    \includegraphics[scale = 0.4]{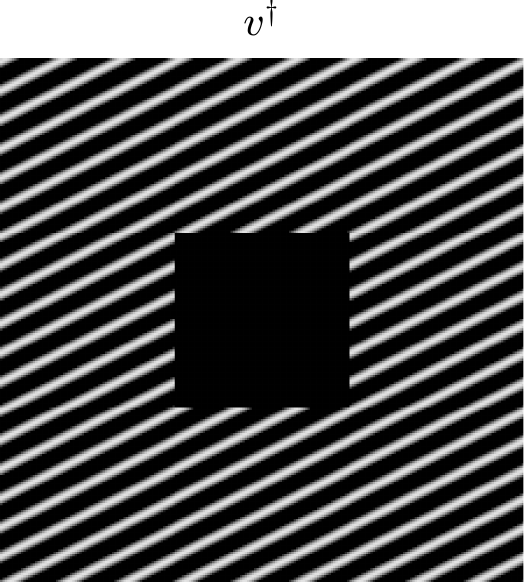}
\end{minipage}
\begin{minipage}[c]{0.24\textwidth}
    \includegraphics[scale = 0.4]{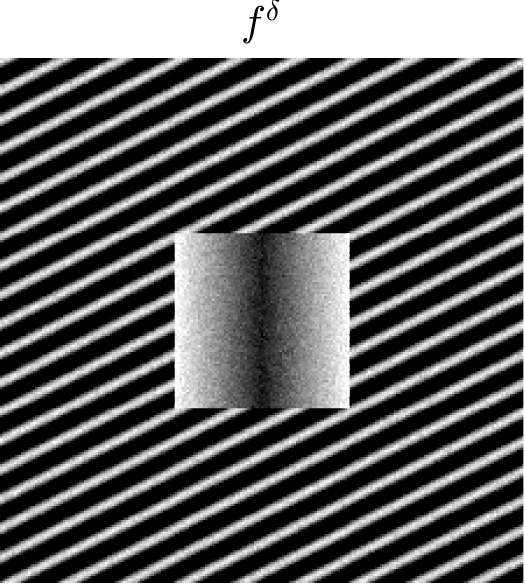}
\end{minipage}

\vspace{0.5cm}
\centering
\hspace{-1cm}\begin{minipage}[c]{0.24\textwidth}
    \includegraphics[scale = 0.4]{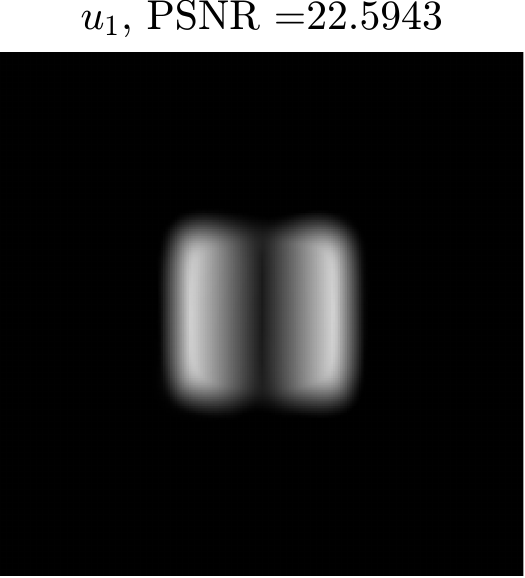}
\end{minipage}
\begin{minipage}[c]{0.24\textwidth}
    \includegraphics[scale = 0.4]{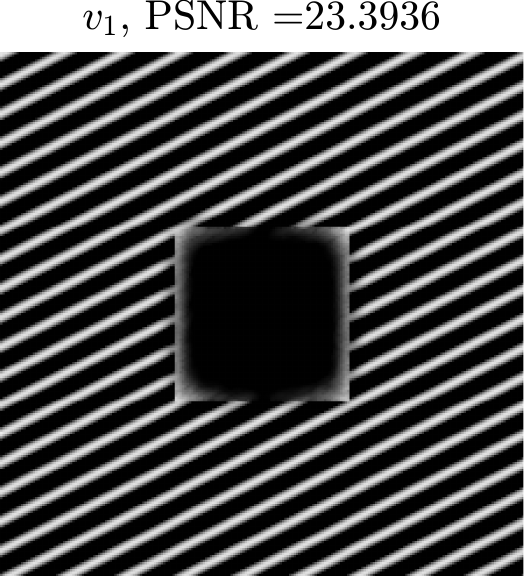}
\end{minipage}
\begin{minipage}[c]{0.24\textwidth}
    \includegraphics[scale = 0.4]{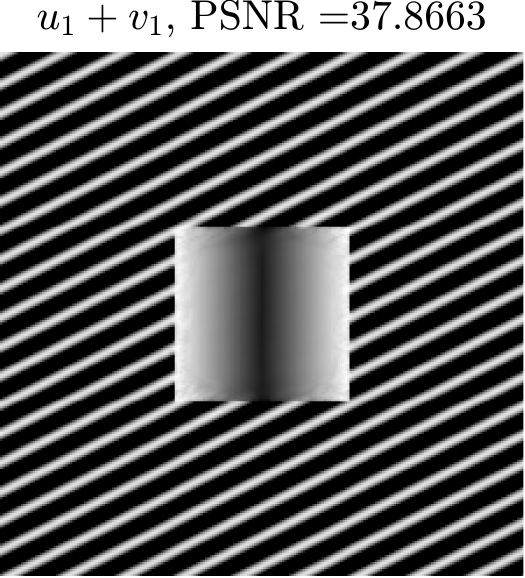}
\end{minipage}

\vspace{0.5cm}
\centering
\hspace{-1cm}\begin{minipage}[c]{0.24\textwidth}
    \includegraphics[scale = 0.4]{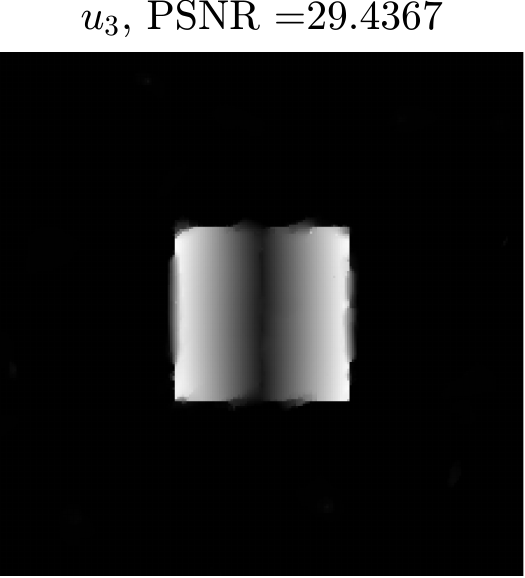}
\end{minipage}
\begin{minipage}[c]{0.24\textwidth}
    \includegraphics[scale = 0.4]{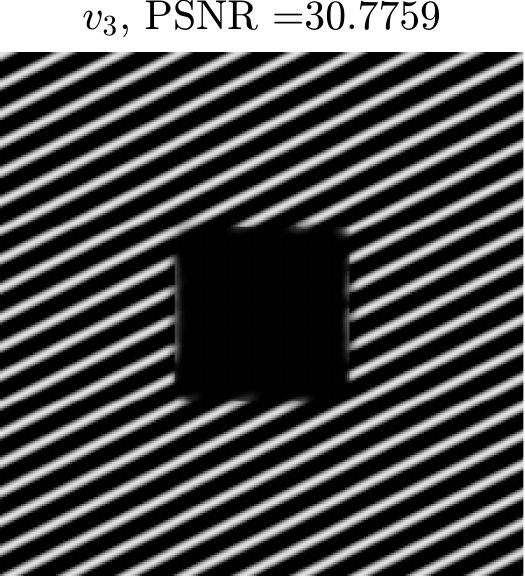}
\end{minipage}
\begin{minipage}[c]{0.24\textwidth}
    \includegraphics[scale = 0.4]{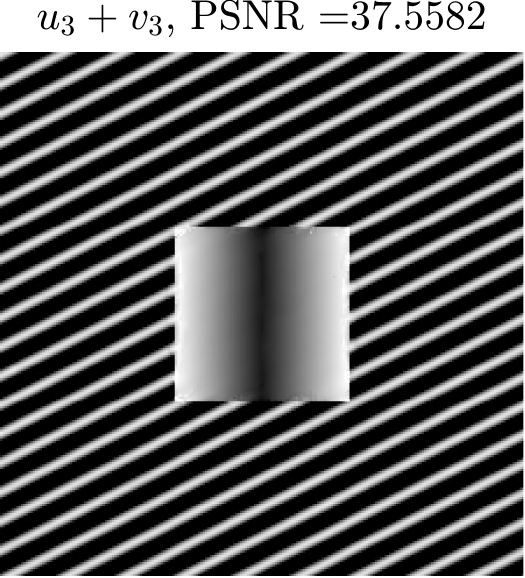}
\end{minipage}

\vspace{0.5cm}
\centering
\hspace{-1cm}\begin{minipage}[c]{0.24\textwidth}
    \includegraphics[scale = 0.4]{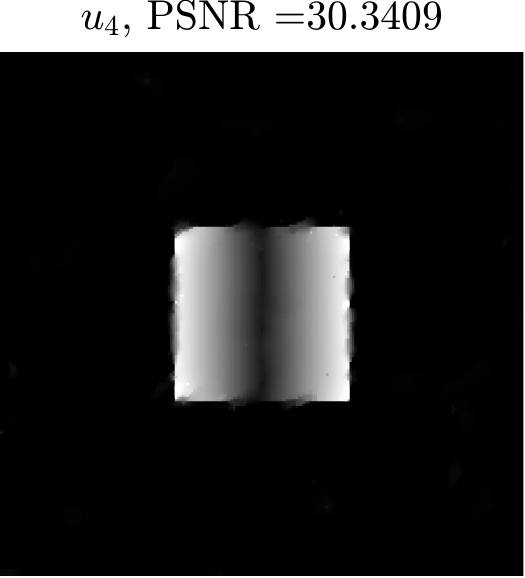}
\end{minipage}
\begin{minipage}[c]{0.24\textwidth}
    \includegraphics[scale = 0.4]{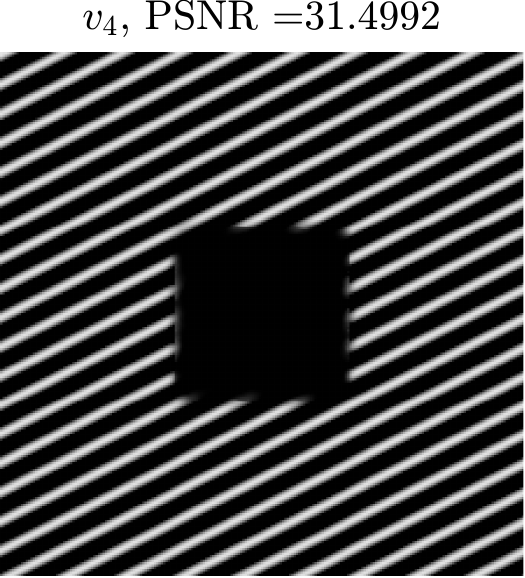}
\end{minipage}
\begin{minipage}[c]{0.24\textwidth}
    \includegraphics[scale = 0.4]{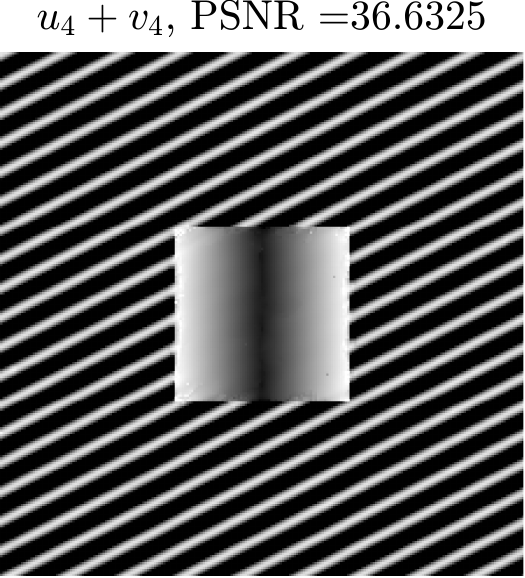}
\end{minipage}

\vspace{0.5cm}
\centering
\hspace{-1cm}\begin{minipage}[c]{0.24\textwidth}
    \includegraphics[scale = 0.4]{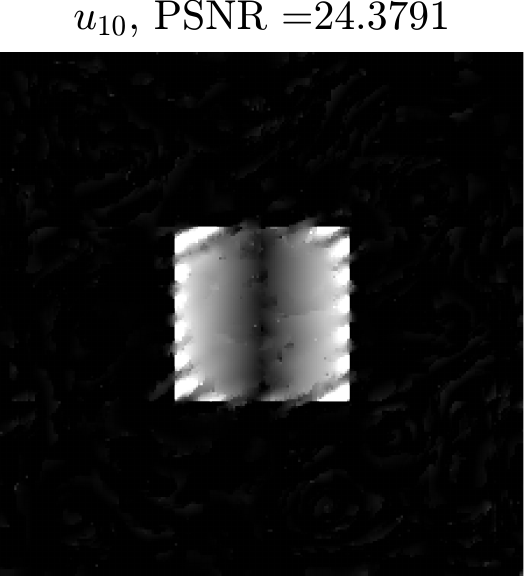}
\end{minipage}
\begin{minipage}[c]{0.24\textwidth}
    \includegraphics[scale = 0.4]{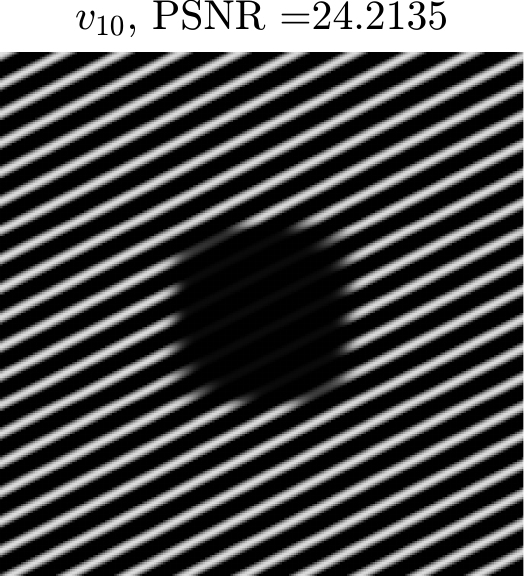}
\end{minipage}
\begin{minipage}[c]{0.24\textwidth}
    \includegraphics[scale = 0.4]{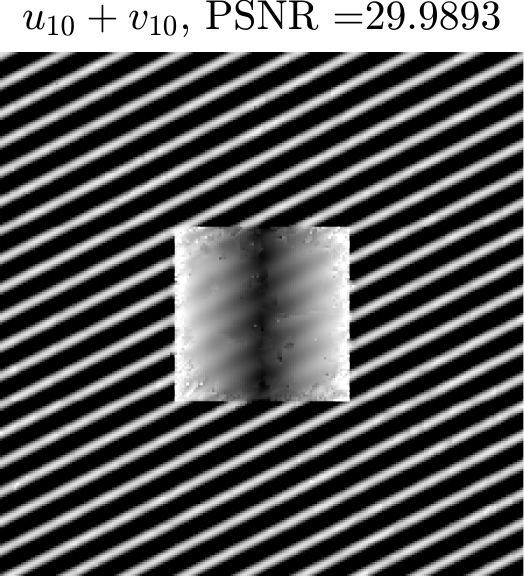}
\end{minipage}

    \caption{ Results obtained by Algorithm \ref{algo:NestedBregman_Noisy_Morozov} with $\alpha_1 = \beta_1 = 5$, $\alpha_2 = \beta_2  = 1$ and $g,h$ as in \eqref{eq:TGV}, \eqref{eq:TGV_osci}. First line: true components and noisy observation. Lines $2-4$ (top to bottom) reconstructed components for the iteration steps $l = 1,3\text{ (minimal cross-correlation)},4\text{ (best PSNR)},10$}
    \label{fig:TGV_osci_Nested_Morozov}
\end{figure}}{}

\section{Conclusion}
We introduce the method of Nested Bregman iterations for decomposition problems, and analyze the method with respect to its well-definedness and convergence behavior. In particular, we make use of the infimal convolution of regularizers to reconstruct solutions for noise corrupted, linear, ill-posed problems, consisting of structurally different components. We illustrate the strength of the method numerically on synthetic test images, where we use a simple cross-correlation based stopping principle. {In our experiments, the results with the proposed method also perform equally or better than the best possible results obtained from a single-step variational approach with suitable parameter choice.} Hence, Nested Bregman iterations {appear superior over a grid search among Morozov-regularized reconstructions, as those can only compete with the proposed method for suitably chosen \textit{predetermined} grids. Furthermore, they} provide a flexible alternative to bilevel methods for image and signal decomposition. However, some theoretical questions on the convergence behavior of the proposed method, where Bregman iterations are used in the inner loop, remain open. In future research, we aim to adapt the Nested Bregman iteration for denoising problems with mixed noise.

\section{Acknowledgements}
We thank Luca Calatroni (Laboratoire d'Informatique, Signaux et Systèmes de Sophia-Antipolis) for his valuable literature suggestions. This research was funded in part by the Austrian Science Fund (FWF) [10.55776/DOC78]. For open access purposes, the authors have applied a CC BY public copyright license to any author-accepted manuscript version arising from this submission.

\printbibliography

@article{Dem,
	Author = {Demengel, F. and Temam, R.},
	Coden = {IUMJAB},
	Fjournal = {Indiana University Mathematics Journal},
	Issn = {0022-2518},
	Issue = 5,
	Journal = {Indiana Univ. Math. J. },
	Pages = {673--709},
	Title = {Convex Functions of a Measure and Applications},
	Volume = 33,
	Year = 1984}

@book{AmbrosioBV,
	Author = {Ambrosio, L. and Fusco, N. and Pallara, D.},
	Isbn = {0198502451},
	Publisher = {Oxford University Press, USA},
	Title = {{Functions of bounded variation and free discontinuity problems}},
	Year = {2000}}

@incollection{handbook_2019,
	title = {Generating structured nonsmooth priors and associated primal-dual methods},
	volume = {20},
	isbn = {978-0-444-64140-3},
	language = {en},
	booktitle = {Handbook of {Numerical} {Analysis}},
	publisher = {Elsevier},
	author = {Hintermüller, Michael and Papafitsoros, Kostas},
	year = {2019},
	doi = {10.1016/bs.hna.2019.08.001},
	pages = {437--502}
}

@Incollection{DelosReyes2021,
author="De los Reyes, Juan Carlos
and Villac{\'i}s, David",
editor="Chen, Ke
and Sch{\"o}nlieb, Carola-Bibiane
and Tai, Xue-Cheng
and Younces, Laurent",
title="Bilevel Optimization Methods in Imaging",
bookTitle="Handbook of Mathematical Models and Algorithms in Computer Vision and Imaging: Mathematical Imaging and Vision",
year="2021",
publisher="Springer International Publishing",
address="Cham",
pages="1--34"
}

@incollection{bilevellearning,
	Author = {Calatroni, L. and Chung, C. and De Los Reyes, J.C. and Sch\"onlieb, C.B. and Valkonen, T.},
	Booktitle = {RADON book Series on Computational and Applied Mathematics, vol. 18},
	Publisher = {Berlin, Boston: De Gruyter},
	Title = {Bilevel approaches for learning of variational imaging models},
	Year = {2017}}

@article{Bredies_2022,
doi = {10.1088/1361-6420/ac668d},
year = {2022},
month = {9},
publisher = {IOP Publishing},
volume = {38},
number = {10},
pages = {105006},
author = {Kristian Bredies and Marcello Carioni and Martin Holler},
title = {Regularization graphs—a unified framework for variational regularization of inverse problems},
journal = {Inverse Problems}
}

@article{HollerKunischIC2014,
	Author = {Holler, M. and Kunisch, K.},
	Doi = {10.1137/130948793},
	Issn = {1936-4954},
	Journal = {SIAM J. Imaging Sci.},
	Number = {4},
	Pages = {2258--2300},
	Title = {On infimal convolution of {TV}-type functionals and applications to video and image reconstruction},
	Volume = {7},
	Year = {2014},
}

@article{TGV_dynamic_MRI,
author = {Schloegl, M. and Holler, M. and Schwarzl, A. and Bredies, K. and Stollberger, R.},
title = {Infimal convolution of total generalized variation functionals for dynamic {MRI}},
journal = {Magnetic Resonance in Medicine},
volume = {78},
number = {1},
pages = {142--155},
doi = {10.1002/mrm.26352},
year={2017}
}

@article{huber1973robust,
  title={Robust regression: asymptotics, conjectures and Monte Carlo},
  author={Huber, Peter J},
  journal={Ann. Statist.},
  pages={799--821},
  year={1973},
  publisher={JSTOR}
}

@InProceedings{tv_linf,
author="Burger, Martin
and Papafitsoros, Konstantinos
and Papoutsellis, Evangelos
and Sch{\"o}nlieb, Carola-Bibiane",
editor="Bociu, Lorena
and D{\'e}sid{\'e}ri, Jean-Antoine
and Habbal, Abderrahmane",
title="Infimal Convolution Regularisation Functionals of BV and $\mathrm{L}^{p}$ Spaces. The Case $p=\infty$",
booktitle="System Modeling and Optimization",
year="2016",
pages="169--179"
}

@article{journal_tvlp,
	Author = {Burger, M. and Papafitsoros, K. and Papoutsellis, E. and Sch\"onlieb, C.B.},
	Journal = {J. Math. Imaging Vision},
	doi = {10.1007/s10851-015-0624-6},
	Number = {3},
	Pages = {343--369},
	Title = {Infimal convolution regularisation functionals of {BV} and $\mathrm{L}^{p}$ spaces. {P}art {I}: The finite $p$ case},
	Volume = {55},
	Year = {2016}}

@article{TGV_learning2,
	Author = {De Los Reyes, J.C. and Sch\"{o}nlieb, C.B. and Valkonen, T.},
	Journal = {J. Math. Imaging Vision},
	doi = {10.1007/s10851-016-0662-8},
	Number = {1},
	Pages = {1-25},
	Title = {Bilevel Parameter learning for higher-order {T}otal {V}ariation regularisation models},
	Volume = {57},
	Year = {2017}}

@article{papafitsoros2013study,
	Author = {Papafitsoros, K. and Bredies, K.},
	Journal = {Inverse Probl. Imaging},
	Number = {2},
	Title = {A study of the one dimensional total generalised variation regularisation problem},
	Volume = {9},
	Year = {2015}}

@article{vese2003modeling,
  title={Modeling textures with total variation minimization and oscillating patterns in image processing},
  author={Vese, Luminita A and Osher, Stanley J},
  journal={J. Sci. Comput.},
  volume={19},
  pages={553--572},
  year={2003},
  publisher={Springer}
}

@book{meyer2001oscillating,
	Author = {Meyer, Y.},
	Publisher = {American Mathematical Society},
	Title = {Oscillating patterns in image processing and nonlinear evolution equations: the fifteenth {D}ean {J}acqueline {B. L}ewis memorial lectures},
	Volume = {22},
	Year = {2001}}

@article{Osher_Sole_Vese_2003,
author = {Osher, Stanley and Sol\'{e}, Andr\'{e}s and Vese, Luminita},
title = {Image Decomposition and Restoration Using Total Variation Minimization and the H1},
journal = {Multiscale Model. Simul.},
volume = {1},
number = {3},
pages = {349-370},
year = {2003},
doi = {10.1137/S1540345902416247},
}

@article{nikolova2002minimizers,
	Author = {Nikolova, M.},
	Journal = {SIAM J. Numer. Anal.},
	doi = {10.1137/S0036142901389165},
	Number = {3},
	Pages = {965--994},
	Publisher = {SIAM},
	Title = {Minimizers of cost-functions involving nonsmooth data-fidelity terms. {A}pplication to the processing of outliers},
	Volume = {40},
	Year = {2002}}

@article{TGV,
	Author = {Bredies, K. and Kunisch, K. and Pock, T.},
	Journal = {SIAM J. Imaging Sci.},
	doi = {10.1137/090769521},
	Number = {3},
	Pages = {492--526},
	Publisher = {SIAM},
	Title = {Total Generalized Variation},
	Volume = {3},
	Year = {2010}}

@article{ChambolleLions,
	Affiliation = {Ceremade (CNRS, URA 749), Universit√{\copyright} de Paris-Dauphine, F-75775 Paris Cedex 16, France FR},
	Author = {Chambolle, A. and Lions, P.L.},
	Issue = {2},
	Journal = {Numer. Math.},
	doi = {10.1007/s002110050258},
	Pages = {167-188},
	Publisher = {Springer Berlin / Heidelberg},
	Title = {Image recovery via total variation minimization and related problems},
	Volume = {76},
	Year = {1997}}

@article{rudin1992nonlinear,
	Author = {Rudin, L.I. and Osher, S. and Fatemi, E.},
	Journal = {Physica D: Nonlinear Phenomena},
	doi = {10.1016/0167-2789(92)90242-F},
	Number = {1-4},
	Pages = {259--268},
	Title = {Nonlinear total variation based noise removal algorithms},
	Volume = {60},
	Year = {1992}}

@inproceedings{Tikhonov,
	Author = {Tikhonov, A.},
	Booktitle = {Soviet Math. Dokl.},
	Pages = {1035},
	Title = {Solution of incorrectly formulated problems and the regularization method},
	Volume = {5},
	Year = {1963}}

@book{BauschkeCombettes,
author = {Bauschke, Heinz H. and Combettes, Patrick L.},
 title = {Convex Analysis and Monotone Operator Theory in Hilbert Spaces},
 year = {2017},
 edition = {2nd},
 publisher = {Springer Publishing Company, Incorporated},

}

@article{bur-osh,
	Author = {Martin Burger and Stanley Osher},
	Doi = {10.1088/0266-5611/20/5/005},
	Journal = {Inverse Problems},
	Number = {5},
	Pages = {1411--1421},
	Title = {Convergence rates of convex variational regularization},
	Volume = {20},
	Year = 2004}

@article{oscICTV,
author = {Yiming Gao and Kristian Bredies},
title = {Infimal Convolution of Oscillation Total Generalized Variation for the Recovery of Images with Structured Texture},
journal = {SIAM J. Imaging Sci.},
volume = {11},
number = {3},
pages = {2021--2063},
year = {2018}}

@article{fri_sch,
  author     = {Frick, K. and Scherzer, O.},
  title      = {Regularization of ill-posed linear equations by the non-stationary augmented {L}agrangian method},
  doi        = {10.1216/JIE-2010-22-2-217},
  journal    = {J. Integral Equations Appl.},
  fjournal   = {J. Integral Equations Appl.},
  month      = jun,
  number     = {2},
  pages      = {217--257},
  volume     = {22},
  year       = {2010},
}

@Article{fri_lor_res,
  author =       "Klaus Frick and Dirk A. Lorenz and Elena Resmerita",
  title =        "{Morozov's Principle} for the Augmented {Lagrangian}
                 Method Applied to Linear Inverse Problems",
  journal={Multiscale Model. Simul.},
  volume =       "9",
  number =       "4",
  pages =        "1528--1548",
  year =         "2011",
  DOI =          "10.1137/100812835",
}

@article{orig_breg,
title = {An Iterative Regularization Method for Total Variation-Based Image Restoration},
author = {Stanley Osher and Martin Burger and Donald Goldfarb and Jinjun Xu and Wotao Yin},
journal = {Multiscale Model. Simul.},
year = {2005}}

@article{li_et_al,
title = {Compressed image sensing with components regularization based on Bregman iteration},
author = {Li, X., Wei Z. and Xiao L.},
journal = {The Journal of China Universities of Posts and Telecommunications},
volume = {18},
number = {2},
pages = {114-119},
year = {2011}}

@phdthesis{stromberg,
  author       = {Str\"omberg, Thomas}, 
  title        = {A study of the operation of infimal convolution},
  school       = {Lulea University of Technology},
  year         = 1994,
  address      = {},
  note         = {}
}

@book{zalinescu,
author = {Zalinescu, Constantin},
 title = {Convex analysis in general vector spaces},
 year = {2002},
 isbn = {},
 edition = {},
 publisher = {World Scientific},
}

@article{Benning-Burger,
    author ={Martin Benning and Martin Burger} ,
    title ={Modern regularization methods for
inverse problems} ,
    journal ={Acta Numrica} ,
    year = {2018},
pages={1 -111}
}

@book{barbu2012convexity,
  title={Convexity and Optimization in Banach Spaces},
  author={Barbu, V. and Precupanu, T.},
  series={Springer Monographs in Mathematics},
  year={2012},
  publisher={Springer Netherlands}
}

@article{Kaltenbacher_Morozov,
year = {2018},
month = {4},
publisher = {IOP Publishing},
volume = {34},
number = {5},
pages = {055008},
author = {Barbara Kaltenbacher and Andrej Klassen},
title = {On convergence and convergence rates for Ivanov and Morozov regularization and application to some parameter identification problems in elliptic PDEs},
journal = {Inverse Problems}
}

@article{Lorenz_2013,
year = {2013},
month = {6},
publisher = {IOP Publishing},
volume = {29},
number = {7},
pages = {075016},
author = {Dirk Lorenz and Nadja Worliczek},
title = {Necessary conditions for variational regularization schemes},
journal = {Inverse Problems},
}

@book{IvanovVasinTanana,
title = {Theory of Linear Ill-Posed Problems and its Applications},
author = {Valentin K. Ivanov and Vladimir V. Vasin and Vitalii P. Tanana},
publisher = {De Gruyter},
address = {Berlin, Boston},
isbn = {9783110944822},
year = {2002},
}

@article{GrasmairHaltmeierScherzer,
title = {The residual method for regularizing ill-posed problems},
journal = {Appl. Math. Comput.},
volume = {218},
number = {6},
pages = {2693-2710},
year = {2011},
issn = {0096-3003},
author = {Markus Grasmair and Markus Haltmeier and Otmar Scherzer},
}

@article{BilevelHuberTV,
author = {Pagliari, Valerio and Papafitsoros, Kostas and Raib\textcommabelow{t}\u{a}, Bogdan and Vikelis, Andreas},
title = {Bilevel Training Schemes in Imaging for Total Variation--Type Functionals with Convex Integrands},
journal = {SIAM J. Imaging Sci.},
volume = {15},
number = {4},
pages = {1690-1728},
year = {2022}
}

@article{HintermuellerInfConv,
author = {Hinterm\"{u}ller, M. and Stadler, G.},
title = {An Infeasible Primal-Dual Algorithm for Total Bounded Variation--Based Inf-Convolution-Type Image Restoration},
journal = {SIAM J. Sci. Comput.},
volume = {28},
number = {1},
pages = {1-23},
year = {2006}
}

@book{BookSpectral,
  title={Mathematical Analysis and Numerical Methods for Science and Technology: Volume 3 Spectral Theory and Applications},
  author={Dautray, R. and Artola, M. and Amson, J.C. and Cessenat, M. and Lions, J.L.},
  year={2012},
  publisher={Springer Berlin Heidelberg}
}

@article{Vese_negative_sobolev2008,
author = {Lieu, Linh and Vese, Luminita},
year = {2008},
month = {10},
pages = {167-193},
title = {Image Restoration and Decomposition via Bounded Total Variation and Negative Hilbert-Sobolev Spaces},
volume = {58},
journal = {Appl. Math. Optim.},
doi = {10.1007/s00245-008-9047-8}
}

@article{Vese_divBMO2005,
author = {Le, Triet and Vese, Luminita},
year = {2005},
month = {01},
pages = {390-423},
title = {Image Decomposition Using Total Variation and div( BMO )},
volume = {4},
journal={Multiscale Model. Simul.},
doi = {10.1137/040610052}
}

@article{AABC_decomposition2005,
author = {Aujol, Jean-François and Aubert, Gilles and Blanc-Féraud, Laure and Chambolle, Antonin},
year = {2005},
month = {01},
pages = {71-88},
title = {Image Decomposition into a Bounded Variation Component and an Oscillating Component},
volume = {22},
journal = {J. Math. Imaging Vision},
doi = {10.1007/s10851-005-4783-8}
}

@article{HuskaDecompositionOnSurfaces_2019,
doi = {10.1088/1361-6420/ab2d44},
year = {2019},
month = {11},
publisher = {IOP Publishing},
volume = {35},
number = {12},
pages = {124008},
author = {Martin Huska and Alessandro Lanza and Serena Morigi and Ivan Selesnick},
title = {A convex-nonconvex variational method for the additive decomposition of functions on surfaces},
journal = {Inverse Problems}
}

@article{HuskaKangMorigiLanza2021,
author = {Huska, Martin and Kang, Sung and Lanza, Alessandro and Morigi, Serena},
year = {2021},
month = {11},
pages = {1749-1789},
title = {A Variational Approach to Additive Image Decomposition into Structure, Harmonic, and Oscillatory Components},
volume = {14},
journal = {SIAM J. Imaging Sci.},
doi = {10.1137/20M1355987}
}

@article{AujaolGilboaChanOsher2006,
author = {Aujol, Jean-François and Gilboa, Guy and Chan, Tony and Osher, Stanley},
year = {2006},
month = {04},
pages = {111-136},
title = {Structure-Texture Image Decomposition—Modeling, Algorithms, and Parameter Selection},
volume = {67},
journal = {International Journal of Computer Vision},
doi = {10.1007/s11263-006-4331-z}
}

@article{TernaryImageDecomposition,
author = {Girometti, Laura and Lanza, Alessandro and Morigi, Serena},
year = {2022},
month = {12},
pages = {},
title = {Ternary image decomposition with automatic parameter selection via auto- and cross-correlation},
volume = {49},
journal = {Adv. Comput. Math.},
doi = {10.1007/s10444-022-10000-4}
}

@inproceedings{QuaternaryImageDecomposition,
author = {Girometti, Laura and Huska, Martin and Lanza, Alessandro and Morigi, Serena},
title = {Quaternary Image Decomposition with Cross-Correlation-Based Multi-parameter Selection},
year = {2023},
isbn = {978-3-031-31974-7},
publisher = {Springer-Verlag},
address = {Berlin, Heidelberg},
doi = {10.1007/978-3-031-31975-4_10},
booktitle = {Scale Space and Variational Methods in Computer Vision: 9th International Conference, SSVM 2023, Santa Margherita Di Pula, Italy, May 21–25, 2023, Proceedings},
pages = {120–133},
numpages = {14},
}

@article{SignalDecomposition2023,
author = {Huska, Martin and Cicone, Antonio and Kang, Sung and Morigi, Serena},
year = {2023},
month = {05},
pages = {153-166},
title = {A Two-stage Signal Decomposition into Jump, Oscillation and Trend using ADMM},
volume = {13},
journal = {Image Processing On Line},
doi = {10.5201/ipol.2023.417}
}
\end{document}